\documentclass{lms} 
\usepackage{tikz-cd}
\usepackage{tikz}
\usepackage[utf8]{inputenc}
\usepackage{calculator}

\usepackage{graphicx}
\usepackage{amsfonts}
\usepackage{amsmath}
\usepackage{amssymb}
\usepackage{hyperref}

\numberwithin{equation}{section}
\newcommand{\nn}{\nonumber}
\newcommand{\Zig}{{\rm Zig}}
\newcommand{\Zag}{{\rm Zag}}
\newcommand{\Exp}{{\rm Exp}}
\newcommand{\wt}{{\rm wt}}

\newcommand{\Tr}{{\rm Tr}}
\newcommand{\Id}{{\rm Id}}
\newcommand{\Hom}{{\rm Hom}}

\newcommand{\Sym}{{\rm Sym}}
\newcommand{\Supp}{{\rm Supp}}

\newcommand{\Index}{{\rm Ind}}

\newcommand{\ie}{i.e. }
\ifx\noTheorems\undefined
\newtheorem{theorem}{Theorem}[section]
\newtheorem{proposition}[theorem]{Proposition}
\newtheorem{lemma}[theorem]{Lemma}
\newtheorem{corollary}[theorem]{Corollary}

\newtheorem{conjecture}[theorem]{Conjecture}
\newtheorem{remark}[theorem]{Remark}
\newtheorem{definition}[theorem]{Definition}
\newtheorem{example}[theorem]{Example}

\title{Cohomological DT invariants from localization}

\author{Pierre Descombes}
\classno{14N35, 16G20}
\extraline{The author is hosted at the Laboratoire de Physique Théorique et Hautes Energies at Sorbonne Université.}
\begin{document}

\maketitle

\begin{abstract}
    Given a quiver with potential associated to a toric Calabi-Yau threefold, the numerical Donaldson-Thomas invariants for the moduli space of framed representations can be computed by using toric localization, which reduces the problem to the enumeration of pyramid partitions. We provide a refinement of this localization procedure, which allows to compute cohomological Donaldson-Thomas invariants. Using this approach, we prove a universal formula which gives the BPS invariants of any toric quiver, up to undetermined contributions which are invariant under Poincaré duality. When the toric Calabi-Yau threefold has compact divisors, these self-Poincaré dual contributions have a complicated dependence on the stability parameters, but explicit computations suggest that they drastically simplify for the self-stability condition (also called attractor chamber). We conjecture a universal formula for the attractor invariants, which applies to any toric Calabi-Yau singularity with compact divisors.
\end{abstract}

\tableofcontents

\section{Introduction}
Donaldson-Thomas (DT) invariants are the mathematical counterpart of the BPS invariants counting supersymmetric bound states in type II string compactifications. On a non-compact toric Calabi-Yau threefold $X$, the study of DT invariants can be translated into a representation-theoretic problem using an equivalence between the bounded derived category of coherent sheaves on $X$ and the bounded derived category of representations of a quiver with potential $(Q,W)$, encoded in a brane tiling. We denote by $Q_0$ (resp $Q_1$) the set of nodes (resp arrows) of the quiver. We will consider cohomological DT invariants as defined by Davison in \cite{Dav13} and Davison-Meinhardt in \cite{DavMein} from idea of Kontsevich-Soibelman in \cite{KonSol10}. The prime object of interest is the generating series of the cohomological DT invariants $\mathcal{A}(x)$, or BPS monodromy, first defined in \cite{KonSol10}: it is a generating series in the Grothendieck ring of monodromic mixed Hodge structures (MMHS). The cohomological BPS invariants $\Omega_{\theta,d}$ defined in \cite[Theo A]{DavMein}, for dimension vector $d\in\mathbb{N}^{Q_0}$ and generic King stability parameter $\theta\in\mathbb{R}^{Q_0}$ are valued in the Grothendieck ring of MMHS, and gives a virtual version of the cohomology with compact support of the moduli space of $\theta$-semistable $d$-dimensional representation. The Harder-Narasimhan decomposition, expressing a general quiver representation as an extension of semistable representation with increasing slope $\mu=\theta.d/\sum_{i\in Q_0}d_i$, and then the Jordan-Hölder filtration, expressing a semistable representation as an extension of stable objects with the same slope, can be expressed by the formula \cite[eq 7]{DavMein}:
\begin{align}\label{framBPS}
    \mathcal{A}(x)=\prod^{\curvearrowright}_l
    \Exp\left(\sum_{d\in l}\frac{\Omega_{\theta,d}}{\mathbb{L}^{1/2}-\mathbb{L}^{-1/2}}x^d\right)
\end{align}
Here $\mathbb{L}^{1/2}$ denotes the square root of the Tate motive, $Exp$ denotes the plethystic exponential defined in \cite[eq 6]{DavMein}, and the product ranges over rays $l$ with increasing slope. Denoting by $\langle\cdot,\cdot\rangle$ the anti-symmetrized Euler form of the quiver introduced below, the attractor invariants $\Omega_{\ast,d}$ defined in \cite[sec. 3.6]{MP20} are special instances of $\Omega_{\theta,d}$ for $\theta$ a small generic deformation of the self-stability (or attractor) condition $\langle \cdot,d\rangle$, subject to the constraint $\theta(d)=0$. The attractor invariants correspond with initial data of the stability scattering diagram introduced in \cite{bridgeland2016scattering}, and one can extract from them all the DT invariants using the recently proven attractor and flow tree formulas, see \cite{Mozgovoy:2021iwz} and \cite{Arguz:2021zpx}. We denote by $\Omega_\theta(x):=\sum_d\Omega_{\theta,d}x^d$ and $\Omega_\ast(x):=\sum_d\Omega_{\ast,d}x^d$ the corresponding  generating series. These series are in general hard to compute, and there is to our knowledge no general closed formula unless $X$ has no compact divisors.\medskip

One way to compute these BPS invariants is to consider $i$-cyclic representations, \ie representations with a vector generating the whole representation at the node $i$. Equivalently, one considers DT invariants for the framed quiver with potential $(Q_i,W)$   (with a single framing node $\infty$ and a single framing arrow $f:\infty\to i$) in the non-commutative stability chamber. Defining the automorphisms $S_{\pm i}(x^d)=\mathbb{L}^{\pm d_i/2}x^d$, the generating series of cohomological framed invariants $Z_i(x)$ is related to the generating series of unframed invariants $\mathcal{A}(x)$ by a wall crossing formula \cite{Mor11,Moz11},
\begin{align}\label{framunframed}
    Z_i(x)=S_i(\mathcal{A}(x))S_{-i}(\mathcal{A}(x)^{-1})
\end{align}
For $D$ a non-compact divisor of $X$, corresponding to a corner of the toric diagram, one can also  consider $D$-cyclic representations, as defined in \cite[sec 3.2]{TYY}.
The corresponding framed quiver $(Q_D,W_D)$ has a single framing node $\infty$ and a pair of arrows $\infty\to i$ and $j\to \infty$ with an additional potential term (see  section \ref{D4def} below for details). We denote by $Z_D(x)$ their cohomological generating series. An $i$-cyclic (resp. $D$-cyclic) representation can be viewed as a noncommutative analogue of a sheaf with a map from the structure sheaf $\mathcal{O}_X$ (resp. $\mathcal{O}_D$). In physics, framed DT invariants count framed BPS states with a D6-brane or non-compact D4-brane charge. Accordingly, we shall refer to the two types of framings as D6- and D4-brane framing, respectively.\medskip

The moduli space of $i$-cyclic (resp. $D$ cyclic) representations admits a maximal torus action rescaling the arrows of $Q_i$ (resp. $Q_f$), leaving the potential $W$ (resp. $W_D$) equivariant, \ie invariant up to a scalar: we denote by $\Lambda$ the character lattice of the torus. We further denote by $\Delta_i$ (resp. $\Delta_D$) the subset of $\Lambda$ (called the Empty Room Configuration, or ERC) given by weights of paths starting at the framing node which are non vanishing in an $i$-cyclic (resp $D$-cyclic) representation of $(Q_i,W)$ (resp $Q_D,W_D$) : $\Delta_i$ can be interpreted a pyramid with an atom of type $i$ on the top, whose facets are given by $\Delta_D$, for $D$ running over the corners of the toric diagram.\medskip

In Lemma \ref{lemfix} we show that the $i$-cyclic (respectively, $D$-cyclic) representations which are fixed under the maximal torus leaving the potential invariant are in bijection with the set $\Pi_i$ of sub-pyramids of $\Delta_i$ (respectively, the set $\Pi_D$ of sub-facets of $\Delta_D$). This allows to translate the computation of the numerical limit of the generating series $Z_i(x)$ (resp. $Z_D(x)$) into a purely combinatoric problem, as proven in \cite[Cor 5.7]{MR}. The formalism of $K$-theoretic localization, developed in \cite{NekOk} allows to compute by toric localization a refinement of the numerical Donaldson-Thomas invariants, known as the $K$-theoretic DT invariants (which are expected to agree with the $\chi_y$ genus evaluation of the cohomological DT invariants), provided the moduli space of framed representations is compact. This formalism therefore applies when the ERC is finite (see for example in \cite{Cir19}). In our situation, the moduli space is non-compact, and the invariants obtained naively by applying the K-theoretic localization formula in the non-compact setting differ from the the cohomological invariants. It can be seen by comparing the computations for the Hilbert scheme of points on $\mathbb{C}^3$ in the $K$-theoretic setting in \cite[sec 8.3]{NekOk}, and in the cohomological setting in \cite{BBS}.\medskip

For a one dimensional torus $\mathbb{C}^\ast$ acting on a smooth scheme, the Bialynicki-Birula decomposition allows to express the cohomology of the attracting variety, \ie the sub-variety of points flowing onto a fixed point when $t\to 0$, as a sum of the cohomology of the fixed points components, shifted by the number of contracting weights in the $\mathbb{C}^\ast$-equivariant tangent space of the fixed locus. The moduli space of cyclic representations of a framed quiver with potential $(Q_f,W_f)$ is not smooth: it is the critical locus of the functional $\Tr(W_f)$, but the general philosophy of derived geometry allows to think about it as a smooth scheme, provided that one replaces the tangent space by the full tangent-obstruction complex. We establish then a derived version of Bialynicki-Birula decomposition. Namely, consider a moduli space $M$ which is the critical locus of a potential on a smooth ambient scheme with a $\mathbb{C}^\ast$-action leaving the potential invariant. A choice of such a $\mathbb{C}^\ast$-action is called a choice of slope, and is denoted by $s$. Denotes by $M^+$ the attracting variety, and by $M^0_\pi$ for $\pi\in\Pi$ the fixed components of the torus action. For $\pi\in\Pi$, denotes by $\Index^s_pi$ the signed number of contracting weight in the restriction to $M^0_\pi$ of the tangent-obstruction complex of $M$. Then we proves \eqref{locnonproj}:
\begin{align}
    [M^+]^{vir}=\sum_{\pi\in\Pi}\mathbb{L}^{\Index^s_\pi/2}[M^0_\pi]^{vir}
\end{align}
This formula holds also when $M$ is a $[-1]$-shifted symplectic scheme or stack, \ie can only be described locally as the critical locus of a potential, as proven in \cite{Descombes:2022cpc}. It explains the observed discrepancy between K-theoretic and cohomological/motivic computations: the K-theoretic computations provides only the refined invariants of the attracting variety of the $\mathbb{C}^\ast$-action given by the slope $s$.\medskip

We apply this result to D6 and D4 brane framings. The fixed components are then isolated points corresponding to pyramids $\pi\in\Pi_i$ (resp $\pi\in\Pi_D$), hence $[M^0_\pi]^{vir}=1$. A choice of slope is then equivalent to a choice of a generic line separating the brane tiling lattice $L$, into two half planes $L^{>0}$ and $L^{<0}$, corresponding to contracting (resp. repelling) weights in the $t\to0$ limit. According to Lemma \ref{lemattr}, the attracting variety of the moduli space of framed representations is then given by representations in which the cycles with repelling weights are nilpotent. To a side $z$ of the toric diagram one associates a vector $l_z\in L$, given by the outward normal to one subdivision of this side, which corresponds to the $L$-weight of a particular cycle of $(Q,W)$ denoted by $v^z$. Those cycles generate all the cycles of $(Q,W)$ (precisely, for $w$ a cycle of $Q$, one has a power $n\in\mathbb{N}$ such that $w^n$ can be written as a product of $v^z$), and correspond to the toric coordinates on $X$ when one views the Jacobian algebra of $(Q,W)$ as a noncommutative crepant resolution of the coordinate ring of $X$. The attracting variety is then the set of framed representations $z$ of the toric diagram such that for $l_z\in L^{<0}$, $v^z$ is nilpotent.\medskip

Imposing nilpotency and invertibility of various cycles of $Q$ amounts to restricting to a Serre sub-category of the category of critical representations of the quiver. Consequently, the formalism of cohomological Hall algebra and wall crossing still applies. For  two disjoint sets of sides of the toric diagram $Z_I$ and $Z_N$, we use the superscript $Z_I:I,Z_N:N$ to denotes the invariants computed by restricting to the representations such that for $z\in Z_I$ (resp. $z\in Z_N$), $v^z$ is invertible (resp. nilpotent). We denotes for convenience by $[z,z']$ the set of sides of the toric diagram between $z$ and $z'$ in the clockwise order, and use the superscript $I$ (resp $N$) to denotes fully invertible (resp fully nilpotents) invariants, \ie invariants counted by considering only representations where all the cycles are invertible (resp nilpotent). For D4 brane framing, we can choose a generic slope $s$ such that for $z,z'$ the sides of the toric diagram adjacent to the corner corresponding to $D$, $l_z,l_{z'}\in L^{>0}$. For a D6 brane framing, there must be always some cycles $v^z$ with repelling weights, hence for a generic slope $s$ we denotes by $[z,z']$ the set of sides $\tilde{z}$ of the toric diagram such that $\tilde{z}\in L^{<0}$. We obtains then:
\begin{align}
    Z_D(x)&=\sum_{\pi\in\Pi_D}\mathbb{L}^{\Index^s_\pi/2}x^{d_\pi}\nn\\
    Z_i^{[z,z']:N}(x)&=\sum_{\pi\in\Pi_i} \mathbb{L}^{\Index^s_\pi/2}x^{d_\pi}
\end{align}
We must then relate the generating series $ Z_i^{[z,z']:N}(x)$ of partially nilpotent $i$-cyclic representations with the full generating series $ Z_i(x)$. It is done using an invertible/nilpotent decomposition of BPS invariants, namely from Proposition \ref{nilinvdec}:
\begin{align}
    \Omega_{\theta}(x)=\Omega^{z:I}(x)+\Omega^{z:N}_\theta(x)
\end{align}
and the fact that from Lemma \ref{lemcenter}, partially invertible representations exists for dimensions vectors in the kernel of the anti-symmetrized Euler form, hence $\Omega^{z:I}(x)$ is in the center of the quantum affine plane and is insensitive to wall crossing. When some cycles are invertible in a representation, we can use the isomorphisms given by the arrows of the cycle to identify the nodes of the cycle, and we obtain a representations of a reduced quiver. Partially invertible BPS invariants of a quiver can then be expressed as BPS invariants of a simpler quiver, and we obtain then universal formulas for them in Section \ref{partinvsec}. We provide some notations for dimensions vectors supporting invertible BPS invariants: to a side $z$ of the toric diagram with $K_z$ subdivisions, one associates zig-zag paths, which are special paths on the brane tiling dividing the torus into $K_z$ parallel strips: for $k\neq k'\in\mathbb{Z}/K_z\mathbb{Z}$, we denote by $\alpha^z_k$ the dimension vector with $1$ on nodes of $Q$ inside the $k$-th strip of the torus, $\alpha^z_{[k,k'[}=\alpha^z_k+\alpha^z_{k+1}+...+\alpha^z_{k'-1}$, and $\delta$ the dimension vector with $1$ on each node. We use these expressions and invertible/nilpotent decompositions to express $\Omega_\theta(x)$ in terms of $\Omega_\theta^{[z,z']:N}(x)$ in Proposition \ref{nilpart}. Using the formula \eqref{framunframed} and \eqref{framBPS} relating framed invariants and BPS invariants, we obtain then:\medskip

\begin{theorem}(Theorem \ref{theolocD})\medskip

    $i)$ For $D$ a non-compact divisor of $X$, corresponding to the corner $p$ of the toric diagram lying between the two sides $z,z'$, and a generic slope $s$ such that $l_z,l_{z'}\in L^{>0}$ (such slopes always exist, because the angle between $l_z$ and $l_{z'}$ is smaller than $\pi$), we have:
    \begin{align}
        Z_D(x)=\sum_{\pi\in\Pi_D}\mathbb{L}^{\Index^s_\pi/2}x^{d_\pi}
    \end{align}\\
    $ii)$ For a generic slope $s$ such that $l_{\tilde{z}}\in L^{<0}\iff\tilde{z}\in[z,z']$, one has:
    \begin{align}
    Z_i(x)&=S_{-i}[\Exp\left(\textstyle\sum_d\Delta^{s}\Omega_d\frac{\mathbb{L}^{d_i}-1}{\mathbb{L}^{1/2}-\mathbb{L}^{-1/2}}x^d\right)]\sum_{\pi\in\Pi_i}\mathbb{L}^{\Index^s_\pi/2}x^{d_\pi}
    \end{align}
    Using the correction term:
    \begin{align}
        \Delta^s\Omega(x)=&(\mathbb{L}^{3/2}+({\textstyle\sum_{\tilde{z}\in[z,z']}} K_{\tilde{z}}-2)\mathbb{L}^{1/2}-({\textstyle\sum_{\tilde{z}\in[z,z']}} K_{\tilde{z}}-1)\mathbb{L}^{-1/2})\sum_{n\geq 1} x^{n\delta}\nonumber\\
    &+(\mathbb{L}^{1/2}-\mathbb{L}^{-1/2})\sum_{\tilde{z}\in[z,z']}\sum_{k\neq k'}\sum_{n\geq 0} x^{n\delta+\alpha^{\tilde{z}}_{[k,k'[}}
    \end{align}
\end{theorem}

This localization formula can be easily implemented on a computer to calculate the framed cohomological DT invariants explicitly for any brane tiling and reasonably small dimension vectors.\medskip

The invertible/nilpotent decomposition allows also to gave a general result about BPS invariants of toric quiver. Namely, we relate in Proposition \ref{nilpart} $\Omega_\theta(x)$, and $\Omega^N_\theta(x)$, the fully nilpotent the BPS invariants. But the Corollary \ref{lemcenter} shows that $\Omega^N_\theta(x)$ is the Poincaré dual of $\Omega_\theta(x)$, hence we prove an universal formula for BPS invariants of toric quivers up to a self Poincaré dual contribution:\medskip

\begin{theorem}(Theorem \ref{theoasym})
\begin{align}
    \Omega_\theta(x)=(\mathbb{L}^{3/2}+(b-3+i)\mathbb{L}^{1/2}+i\mathbb{L}^{-1/2})\sum_{n\geq 1} x^{n\delta}+\mathbb{L}^{1/2}\sum_{z}\sum_{k\neq k'}\sum_{n\geq 0} x^{n\delta+\alpha^z_{[k,k'[}}+\Omega^{sym}_\theta(x)
\end{align}
with $\Omega^{sym}_\theta(x)$ self Poincaré dual, and supported on dimension vectors $d\not\in\langle\delta\rangle$. The same formula holds for attractor invariants.
\end{theorem}

For toric Calabi-Yau threefolds without compact divisors (also known as local curves, 
corresponding to toric diagrams with no internal points), the quiver $Q$ is symmetric,
and consequently  the unframed DT invariants do not exhibit wall-crossing. They are known in most cases, see \cite{BBS,MMNS,MorNag} and
\cite[sec 5]{MP20} for a review. We check that theorem \ref{theoasym} is consistent with these results: in some cases, including the simplest case of the conifold, there exists infinite towers of dimension vectors $d$ with $\Omega_\theta(d)=1$, associated to rational curves with normal bundle $\mathcal{O}(-1)+\mathcal{O}(-1)$, whose contributions are included in $\Omega^{sym}(x)$. In contrast, the dimension vectors with $\Omega_\theta(d)=\mathbb{L}^{1/2}$ appearing in Theorem \ref{theoasym} are associated to rational curves with normal bundle $\mathcal{O}(-2)+\mathcal{O}(0)$. In some cases one can find 'preferred slopes' (as shown in \cite[sec 4.3]{Arb}) where many cancellations occur in the index, and obtain a closed formula for the full BPS invariants from the cohomological localization: we check that it agrees with the cohomological computations for $\mathbb{C}^3$, the conifold and $\mathbb{C}^2/(\mathbb{Z}/2\mathbb{Z})\times\mathbb{C}$ in Section \ref{localcurves}.\medskip

For Calabi-Yau threefolds with compact divisors, corresponding to asymmetric quivers, there is no closed formula to our knowledge for numerical invariants, let alone for the cohomological ones. In particular BPS invariants depend on the King stability parameter $\theta$, and the symmetric contribution $\Omega^{sym}_\theta(x)$ is quite intricate for arbitrary $\theta$. In \cite{BMP20}, toric quivers associated to toric Fano surfaces (\ie toric diagrams with one interior point and no interior boundary points) are studied. It is conjectured in \cite[p. 21]{BMP20}, \cite[Conj 1.2]{MP20} that in this case the only attractor invariants are those supported on dimension vectors $e_i$ for $i\in Q_0$, corresponding to simple representations, and those supported on the dimensions vectors $\mathbb{N}^\ast\delta$, corresponding to D0 branes, \ie Hilbert schemes of points. In \cite{MP20} weak toric Fano surfaces (\ie toric diagrams with one interior point and interior boundary points) are considered: it is observed that there can be additional dimension vectors with non-vanishing attractor invariants, but it is conjectured (\cite[Conj 1.1]{MP20}) that they lie in the kernel of the anti-symmetrized Euler form. We shall formulate a refinement of these conjectures:\medskip

\begin{conjecture}(Conjecture \ref{conj}) 
For toric diagram with $i\geq 1$ internal lattice points, the attractor invariants are given by:
   \begin{align}
       \Omega_\ast(x)=\sum_i x_i+(\mathbb{L}^{3/2}+(b-3+i)\mathbb{L}^{1/2}+i\mathbb{L}^{-1/2})\sum_{n\geq 1} x^{n\delta}+\mathbb{L}^{1/2}\sum_{z}\sum_{k\neq k'}\sum_{n\geq 0} x^{n\delta+\alpha^z_{[k,k'[}}
    \end{align}
\end{conjecture}

The attractor invariants associated to simple representations and Hilbert scheme of points are known. When there are $K_z-1$ lattice points on a side $z$ of the toric diagram, the toric threefold $X$ exhibits a $\mathbb{C}^2/\mathbb{Z}_{K_z}\times \mathbb{C}^\ast$ singularity away from the zero locus of the toric coordinate corresponding to $z$, as recalled in the proof of Proposition \ref{1partinv}. The conjecture then predicts that the only additional attractor invariants correspond to D2-branes wrapped on rational curves in this extended singularity.\medskip

The rest of this article is organized as follows:
\begin{itemize}
    \item In section 2 we review known results on Donaldson-Thomas theory on toric threefolds, and introduce the basic definitions and notations. In section 2.1, we introduce the moduli spaces of representations associated to unframed and framed quiver, their cohomological DT invariants and generating series thereof. In section 2.2 we recall how the quiver with potential for toric Calabi-Yau threefolds can be deduced from brane tiling, and emphasize the utility of zig-zag paths. In section 2.3 we introduce the D6- and D4-framing.\medskip
    
    \item In section 3, using invertible/nilpotent decompositions of unframed representations, we relate generating series of BPS invariants with various nilpotency constraints. In section 3.1 we introduce the notion of partially invertible/nilpotent representations, and define their generating series and BPS invariants. In section 3.2 we show that the invertible/nilpotent decomposition on unframed representations implies a decomposition of BPS invariants. In section 3.3, we compute BPS invariants for partially invertible representations. In section 3.4 we orchestrate previous results, express the BPS invariants in terms of the partially nilpotent BPS invariants accessible by toric localization, and prove the Theorem \ref{theoasym}.\medskip
    
    \item In section 4, we study toric localization for framed quivers. In section 4.1 we describe the fixed locus and the attracting locus of the toric action scaling the arrows of D4- and D6-framed representations. In section 4.2, we describe the $\mathbb{C}^\ast$-equivariant tangent-obstruction complex at a $\mathbb{C}^\ast$-fixed component of the moduli space. In section 4.3, we prove the 'derived Bialynicki-Birula decomposition' for general framed quivers and provide our localization theorem \ref{theolocD} for D4- and D6-framed invariants. In section 4.4, we relate our localization result to the localization of K-theoretic DT invariants.\medskip
    
    \item  In section 5 we illustrate our formula and formulate our conjecture for the complete set of attractor invariants. In section 5.1 we compare our results with the known formulas for local curves, and explains on specific examples the discrepancy between K-theoretic and cohomological computations. In section 5.2 we compare our result \ref{theoasym} and conjecture \ref{conj} with the computations in \cite{BMP20,MP20} for toric threefolds with one compact divisor. In order to facilitate comparison with future computations, we spell out our conjecture for the canonical bundle over toric weak Fano surfaces, using the brane tilings listed in \cite{HanSeo}. 
\end{itemize}\medskip

$Acknowledgments$ I am grateful to my PhD advisors Boris Pioline and Olivier Schiffmann for useful discussions, and all their advice and suggestions during the writing of this article. I also thank Ben Davison and Sergey Mozgovoy for many useful suggestions, and the anonymous referee for precious advice on the first version of this paper.

\section{Basic notions on Donaldson-Thomas theory and toric quivers}

\subsection{Invariants of quivers with potential}

\subsubsection{Representations and cohomological DT invariants}

Consider a quiver with potential $(Q,W)$, with $Q_0$ (resp. $Q_1$) the set of nodes (resp. arrows) of $Q$, the source and target of an arrow $a$ being denoted respectively $s(a)$ and $t(a)$, and $W$ a linear combination of cycles of $Q$ (we follow the notations of \cite{MR} whenever possible). The path algebra of the quiver $Q$,  denoted by $\mathbb{C}Q$, is the free algebra generated by arrows of the quiver, such that $ba=0$ if $s(b)\neq t(a)$. A cycle is a path $w=a_1...a_n$ with $s(a_n)=t(a_1)$. The cyclic derivative is defined by 
\begin{align}
    \partial_a w=\sum_{i:a_i=a}a_{i+1}...a_na_1...a_{i-1}
\end{align}
and extended to $\mathbb{C}Q$ by linearity. The cyclic
derivatives of the potential define the ideal $(\partial W)=((\partial_aW)_{a\in Q_1})$. The Jacobian algebra is the quotient  $J_{Q,W}=\mathbb{C}Q/(\partial W)$. 
We shall usually identify a path with its image in $J_{Q,W}$, \ie paths
which differ by derivatives of the potential will be identified.

\medskip

Consider a framed quiver with potential $(Q_f,W_f)$ obtained from $(Q,W)$ by adding a single framing node $\infty$, (possibly multiple) framing arrows between the framing node and nodes of $Q$, and (when allowed)  additional cycles in the potential, corresponding to path starting and ending at the framing node. One consider the projective $\mathbb{C}Q_f$ module $\mathfrak{P}_f$ generated by paths of $Q_f$ starting at the framing node. One can also consider the  Jacobian algebra $J_{Q_f,W_f}:=\mathbb{C}Q_f/(\partial W_f)$ and the left $J_{Q_f,W_f}$ module $P_f:=\mathfrak{P}_f/((\partial W_f)\cap\mathfrak{P}_f)$.\medskip

For any dimension vector $d\in\mathbb{N}^{Q_0}$, we denote by $\mathfrak{M}_{Q,d}$ the moduli stack of $d$-dimensional representations of the unframed quiver $Q$ (without imposing the potential relations), \ie the moduli stack of left $\mathbb{C}Q$ modules, which can be expressed more explicitly by:
\begin{align}
    \mathfrak{M}_{Q,d}=\frac{\prod_{(a:i\to j)\in Q_1}\Hom(\mathbb{C}^{d_i},\mathbb{C}^{d_j})}{\prod_{i\in Q_0}GL_{d_i}}
\end{align}
Here the gauge group $G_d=\prod_{i\in Q_0} GL_{d_i}$ acts on $a\in Hom(\mathbb{C}^{d_i},\mathbb{C}^{d_j})$ by $a\mapsto g_j a g_i^{-1}$. For a stability parameter $\theta$, we denote by $\mathcal{M}^{\theta,ss}_{Q,d}$ (resp $\mathcal{M}^{\theta,s}_{Q,d}$) the moduli space of $\theta$-semistable representations (resp the smooth open subset of $\theta$-stable representations), obtained by geometric invariant theory as in \cite{king_moduli}.\medskip

Similarly, for any dimension vector $d\in\mathbb{N}^{Q_0}$ we denote by $\mathcal{M}_{Q_f,d}$ the moduli space of $f$-cyclic representations of the framed quiver $Q_f$ with dimension vector $d'=(d,1)\in \mathbb{N}^{Q_f}_0$ , \ie representations with dimension $1$ on the framing node, such that the sub-representation generated by the framing node is the whole representation:
\begin{align}
    \mathcal{M}_{Q_f,d}=\frac{(\prod_{(a:i\to j)\in (Q_f)_1}\Hom(\mathbb{C}^{d'_i},\mathbb{C}^{d'_j}))^{\rm cycl}}{\prod_{i\in Q_0}GL_{d_i}}
\end{align}
Here the subscript ${\rm cycl}$ denotes the open subset of $f$-cyclic representations. $f$-cyclic representations are $\theta$-stable representations of $Q_f$, for a stability parameter $\theta\in\mathbb{R}^{(Q_f)_0}$ such that $\theta.d=0$, $\theta_\infty>0$ and $\theta_i<0$ for $i\in Q_0$, hence from geometric invariant theory $\mathcal{M}_{Q_f,d}$ is a smooth scheme. Equivalently, $\mathcal{M}_{Q_f,d}$ is the scheme which corresponds to $d$-dimensional quotients of the module of paths $\mathfrak{P}_f$, \ie quotient by a $\mathbb{C}Q_f$ sub-module $\rho$ of codimension $d$.\medskip

We consider the functional $\Tr(W)$ on $\mathfrak{M}_{Q,d}$ and $\mathcal{M}^{\theta,ss}_{Q,d}$ (resp. $Tr(W_f)$ on $\mathcal{M}_{Q_f,d}$), and their critical locus $\mathfrak{M}_{Q,W,d}$ and $\mathcal{M}^{\theta,ss}_{Q,W,d}$ (resp. $\mathcal{M}_{Q_f,W_f,d}$). Representations in the critical locus are called critical representations, and correspond to left $J_{Q,W}$ modules (resp. quotients of $P_f$). One denotes by $\phi_W$ (resp. $\phi_{W_f}$) the vanishing cycle functor of $Tr(W)$ (resp. $\Tr(W_f)$), having support on critical representations: it is a functor with source the category of mixed Hodge modules on $\mathfrak{M}_{Q,d}$ (resp. $\mathcal{M}^{\theta,ss}_{Q,d}$, resp. on $\mathcal{M}_{Q_f,d}$), and target the category of monodromic mixed Hodge modules on $\mathfrak{M}_{Q,d}$ and $\mathcal{M}^{\theta,ss}_{Q,d}$ (resp. on $\mathcal{M}_{Q_f,d}$) with support on $\mathfrak{M}_{Q,W,d}$ and $\mathcal{M}^{\theta,ss}_{Q,W,d}$ (resp. $\mathcal{M}_{Q_f,W_f,d}$).\medskip

Consider a constructible sub-stack $\mathfrak{M}^S_{Q,d}$ of $\mathfrak{M}_{Q,d}$ (resp. a constructible subscheme $\mathcal{M}^{\theta,ss,S}_{Q,d}$ of $\mathcal{M}^{\theta,ss}_{Q,d}$, resp. $\mathcal{M}^S_{Q_f,d}$ of $\mathcal{M}_{Q_f,d}$). In this work we shall consider sub-stacks or sub-schemes that are attracting varieties of a toric action, or giving representations with nilpotency and invertibility constraints on particular cycles. Following the general formalism of cohomological Donaldson-Thomas invariants developed in \cite{KonSol10},\cite{Dav13} and \cite{DavMein}, one defines the cohomological DT invariants of critical representations in the Grothendieck group of monodromic mixed Hodge structures:
\begin{align}\label{critcoho}
    [\mathfrak{M}^S_{Q,W,d}]^{vir}=&H_c^\bullet(\mathfrak{M}^S_{Q,d},\phi_W\mathcal{IC}_{\mathfrak{M}_{Q,d}})\nn\\
    [\mathcal{M}^S_{Q_f,W_f,d}]^{vir}=&H_c^\bullet(\mathcal{M}^S_{Q_f,d},\phi_{W_f}\mathcal{IC}_{\mathcal{M}_{Q_f,d}})
\end{align}
denoting by $H^\bullet_c(M,F)$ the Grothendieck class of the cohomology with compact support of the complex of the monodromic mixed Hodge module $F$ on $M$ and by $\mathcal{IC}_M$ the intersection complex of $M$, or more precisely the corresponding mixed Hodge module. We omit the superscript $S$ when we consider the entire stack or scheme of representations.\medskip

In \cite{DavMein}, Davison and Meinhardt introduce the BPS sheaf on $\mathcal{M}^{\theta,ss}_{Q,d}$:
\begin{align}
    \mathcal{BPS}_{W,d}^\theta:=\left\{
    \begin{array}{ll}
        \phi_W\mathcal{IC}_{\mathcal{M}^{\theta,ss}_{Q,d}} & \mbox{if }  \mathcal{M}^{\theta,st}_{Q,d}\neq\emptyset\\
        0 & \mbox{otherwise}
    \end{array}
\right.
\end{align}
The BPS invariants are then defined by:
\begin{align}
     \Omega_{\theta,d}^S=H_c^\bullet(\mathcal{M}^{\theta,ss,S}_{Q,d},\mathcal{BPS}_{W,d}^\theta)
\end{align}

For a dimension vector $d$, considering a small deformation $\theta_d$ of the self stability condition $\langle -,d\rangle$, generic such that $\theta_d(d)=0$, we define the attractor invariants:
\begin{align}
    \Omega^S_{\ast,d}=\Omega^S_{\theta_d,d}
\end{align}
Then \cite[Theo 3.7]{MP20}, based on the theory of cluster scattering diagrams developed in \cite{gross_canonical}, states that attractors invariants $\Omega_{\ast,d}$ are well defined, \ie they do not depend on the small generic deformation $\theta_d$. Since the formalism of cluster scattering diagram also applied when restricting to a Serre sub-category of the category of representations of a quiver, the same arguments ensure that $\Omega^S_{\ast,d}$ are also well defined.
\medskip

\begin{remark}
We will interpret previous computations in the motivic setting as formulas in  this Grothendieck ring of MMHS, using the realization map from monodromic motives to MMHS (the compatibility between the motivic and cohomological definitions was checked in \cite[Appendix A]{davison2019refined}). We replace the multiplication by the square root $\mathbb{L}^{1/2}$ of the Tate motive by the cohomological shift $[-1]$ at the level of perverse sheaves, or by the tensor product with the MMHS given by the vanishing cycles of $z\to z^2:\mathbb{C}\to\mathbb{C}$ as in \cite[p. 19]{DavMein}, which is a square root of the MMHS of the affine line. In particular, when we computes the Hodge polynomial associated to the monodromic mixed Hodge structure, we replace $\mathbb{L}^{1/2}$ by $(-y)$, and by $-1$ in the numerical limit, in agreement with \cite{BBS}.
\end{remark}\medskip

\begin{remark}
In \cite{Dav13} and \cite{DavMein}, Davison and Meinhardt use Borel-Moore homology, \ie the dual of the cohomology with compact support, hence their invariants are the Poincaré dual of our invariants. Here we follow the convention of \cite{MP20}, which is also the convention used in the literature about motivic invariants of quivers with potential. 
\end{remark}

\subsubsection{Quantum affine space and generating series}

For $d,d'\in \mathbb{Z}^{Q_0}$, the Euler form $\chi_Q$ and its anti-symmetrized version $\langle,\rangle$ are defined by:
\begin{align}
    &\chi_Q(d,d')=\sum_{i\in Q_0}d_id'_i-\sum_{(a:i\to j)\in Q_1}d_id'_j\nn\\
    &\langle d,d'\rangle=\chi_Q(d,d')-\chi_Q(d',d)
\end{align}
The quantum affine space $\hat{\mathbb{A}}$ is the algebra generated by elements $x^d$, for $d\in\mathbb{N}^{Q_0}$, with coefficients in the Grothendieck group (having a ring structure) of monodromic mixed Hodge structures, and relations:
\begin{align}
    x^dx^{d'}=\mathbb{L}^{\langle d,d'\rangle/2}x^{d+d'}
\end{align}
We introduce the algebra automorphism $S_{\pm i}$ of the quantum affine space $\hat{\mathbb{A}}$ (denoting $P$ a class of the Grothendieck group of monodromic mixed Hodge structure):
\begin{align}
    &S_{\pm i}:Px^d\mapsto \mathbb{L}^{\pm d_i/2}Px^d\nn\\
\end{align}
Consider $S$, a Serre subcategory of the Abelian category of representations of $Q$, \ie a full subcategory such that for each exact sequence:
\begin{align}
    0\to V_1\to V\to V_2\to 0
\end{align}
$V\in S$ if and only if $V_1\in S$ and $V_2\in S$. We denote by $\mathfrak{M}^S_{Q,d}$, $\mathcal{M}^{\theta,ss,S}_{Q,d}$ and $\mathcal{M}^S_{Q_f,d}$ the sub-stack and sub-schemes of representations lying in $S$ (resp such that the induced representation of $Q$ lies in $S$). The generating series of unframed or framed invariants restricted to the Serre subcategory $S$, with values in the quantum affine space $\hat{\mathbb{A}}$, are defined by:
\begin{align}\label{defgenser}
    \mathcal{A}^S(x)=&\sum_d[\mathfrak{M}^S_{Q,W,d}]^{vir}x^d\nn\\
    Z_f^S(x)=&\sum_d[\mathcal{M}^S_{Q_f,W_f,d}]^{vir}x^d\nn\\
    \Omega^S_\theta(x)&=\sum_d\Omega^S_{\theta,d}x^d\nn\\
    \Omega^S_\ast(x)&=\sum_d\Omega^S_{\ast,d}x^d
\end{align}
As recalled in the introduction, the Harder-Narasimhan decomposition express a general quiver representation as an extension of semistable representation with increasing slope $\mu=\theta.d/\sum_{i\in Q_0}d_i$, and the Jordan-Hölder filtration, express a semistable representation as an extension of stable objects with the same slope. Consider a stability condition $\theta$ which is generic, \ie such that if $d,d'$ have the same slope then $\langle d-d',\bullet\rangle=0$. The Harder-Narasimhan and Jordan-Hölder decompositions can then be expressed by the formula \cite[eq 7]{DavMein}:
\begin{align}
    \mathcal{A}^S(x)=\prod^{\curvearrowright}_l
    \Exp\left(\sum_{d\in l}\frac{\Omega^S_{\theta,d}}{\mathbb{L}^{1/2}-\mathbb{L}^{-1/2}}x^d\right)
\end{align}
Here $Exp$ denotes the plethystic exponential defined in \cite[eq 6]{DavMein}, and the product ranges over rays $l$ with increasing slope.\medskip

\subsection{Unframed quivers associated to toric threefolds}

\subsubsection{From toric diagrams to quivers with potentials}\medskip

Let us consider a toric Calabi Yau threefold $X$. The fast inverse algorithm described in \cite[sec 5]{HanVegh} gives a brane tiling on the two- dimensional torus from the toric diagram of $X$, \ie a bipartite graph with white and black vertex and edges between a white and a black vertex. In fact, it can give different brane tilings that are related by toric mutations.
\medskip

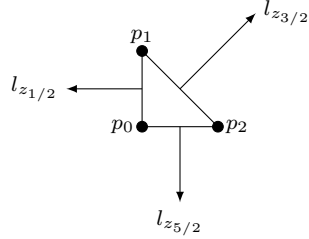
\begin{figure}\caption{Toric diagram of $\mathbb{C}^3$}\label{torc3}
\centerline{
\begin{tikzpicture} 
\filldraw [black] (0,0) circle (2pt);
\draw (0,0) node[left]{$p_0$};
\filldraw [black] (0,1) circle (2pt);
\draw (0,1) node[above]{$p_1$};
\filldraw [black] (1,0) circle (2pt);
\draw (1,0) node[right]{$p_2$};
\draw (0,0) -- (0,1); 
\draw (0,1) -- (1,0); 
\draw (1,0) -- (0,0); 
\draw[->,>=latex] (0.5,0) -- (0.5,-1);
\draw (0.5,-1) node[below]{$l_{z_{5/2}}$};
\draw[->,>=latex] (0,0.5) -- (-1,0.5);
\draw (-1,0.5) node[left]{$l_{z_{1/2}}$};
\draw[->,>=latex] (0.5,0.5) -- (1.5,1.5);
\draw (1.5,1.5) node[right]{$l_{z_{3/2}}$};
\end{tikzpicture}}
\end{figure}

We consider the toric diagram of $X$, which is a convex polygon in a two dimensional free lattice $L^{\vee}$. We denote by $n$ the number of corners of the toric diagram, and the corners themselves by $p_i$ for $i\in \mathbb{Z}/n\mathbb{Z}$ in the clockwise order. The side of the toric diagram between two adjacent corners $p_i$ and $p_{i+1}$ will be denoted $z_{i+1/2}$. We denote by $K_z$ the number of subdivisions of the edge $z$, i.e. the number of the lattice points on that edge (counting the endpoints) minus one. We denote by $l_z\in L$ the primitive vector generating the dual of the side $z$ in $L$. As an example, for $\mathbb{C}^3$, the toric diagram and vectors $l_z$ are given by Figure \ref{torc3}.

Let us now describe the fast inverse algorithm. On the real two dimensional torus obtained by dividing $\mathbb{R}^2$ by the lattice $L$, we draw for each edge $z$ of the toric diagram $K_z$ generic oriented lines directed along $l_z$, in generic position such that two lines intersect only in one point and three lines do not intersect. The different choices in the relative arrangement of lines will correspond to different quivers with potential related by toric mutations. The complement of these lines determines polygonal domains, or tiles, with oriented edges. We color those tiles in white, dark grey or light grey, according to the orientations of their edges:
\begin{itemize}
    \item If the edges of the tile are oriented in the clockwise order around the tile, we color the tile in dark grey
    \item if the edges of the tile are oriented in the counter-clockwise order around the tile, we color the tile in light grey
    \item if the orientations of the different edges of the tile do not agree, we color the tile in white
\end{itemize}

We define a brane tiling on the torus by putting a black node in each dark grey tile, a white node in each light grey tile, and connecting a black node and a white node if the corresponding tiles are connected at one of their corners. The white tiles are then in correspondence with tiles of the brane tiling.
\medskip

\begin{figure}\caption{The fast inverse algorithm for $\mathbb{C}^3$}\label{fiac3}
\centerline{\begin{tikzpicture}
\draw[gray!40] (0,0) -- (0,3) -- (3,3) -- (3,0) -- (0,0);
\fill[gray!80] (1,1) -- (2,1) -- (2,2) -- (1,1);
\fill[gray!20] (2,1) -- (3,1) -- (3,0) -- (2,0) -- (2,1);
\fill[gray!20] (2,2) -- (2,3) -- (3,3) -- (2,2);
\fill[gray!20] (0,0) -- (1,1) -- (0,1) -- (0,0);
\draw[->,>=latex] (2,3) -- (2,0);
\draw[->,>=latex] (3,1) -- (0,1);
\draw[->,>=latex] (0,0) -- (3,3);
\draw[->,>=latex] (2,2) -- (2,1);
\draw[->,>=latex] (2,1) -- (1,1);
\draw[->,>=latex] (1,1) -- (2,2);
\draw[->,>=latex] (2,3) -- (2,2);
\draw[->,>=latex] (3,1) -- (2,1);
\draw[->,>=latex] (0,0) -- (1,1);
\node[draw,circle,inner sep=2pt,fill] (B) at (1.7,1.3){};
\node[draw,circle,inner sep=2pt,fill,white] at (2.5,0.5){};
\node[draw,circle,inner sep=2pt] (W) at (2.5,0.5){};
\draw[dashed] (B) -- (W);
\draw[dashed] (B) -- (0,0.75);
\draw[dashed] (B) -- (2.25,3);
\draw[dashed] (W) -- (3,0.75);
\draw[dashed] (W) -- (2.25,0);
\draw[->,>=latex] (1.5,0.5) -- (2.5,1.5);
\draw[->,>=latex] (0.75,1.5) -- (1.25,0.5);
\draw[->,>=latex] (2.5,1.75) -- (1.5,2.25);

\begin{scope}[scale=2]
\node (A) at (3,1) {$1$};
\path[->,,>=latex] (A) edge [in=100,out=140,loop] node [left]{a} ();
\path[->,,>=latex] (A) edge [in=220,out=260,loop] node [left]{b} ();
\path[->,,>=latex] (A) edge [in=-20,out=20,loop] node [right]{c} ();
\node at (3,0.25) {$W=abc-acb$};
\end{scope}
\end{tikzpicture}}
\end{figure}
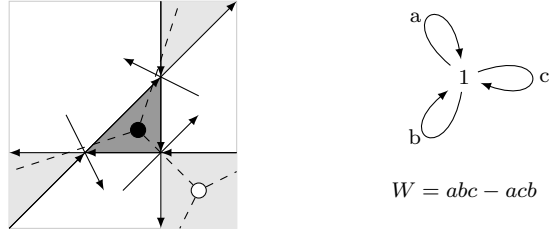

\begin{definition}
The quiver with potential $(Q,W)$ associated of a brane tiling is defined as the dual of this brane tiling, \ie:
\begin{itemize}
    \item The set of nodes $Q_0$ of the quiver is the set of tiles of the brane tiling
    \item The set of arrows $Q_1$ of the quiver is the set of edges of the brane tiling. An edge of the tiling between two tiles gives an arrow of the quiver between the two corresponding nodes, oriented such the black node is at the left of the arrow
    \item Denote by $Q_2$ the set of nodes of the brane tiling, and $Q_2^+$ (resp. $Q_2^-$) the subset of white (resp. black) nodes. To a node $F\in Q_2$ one associate the cycle $w_F$ of $Q$ composed by arrows surrounding this node. We define:
    \begin{align}
        W=\sum_{F\in Q^+_2}w_F-\sum_{F\in Q^-_2}w_F
    \end{align}
\end{itemize} 
\end{definition} 
By definition, the quiver with potential $(Q,W)$ is drawn on a torus: the unfolding of this quiver to the universal cover $\mathbb{R}^2$ of the torus is called the periodic quiver. In the case of $\mathbb{C}^3$, this procedure is described in Figure \ref{fiac3}.\medskip

\begin{definition}
    A zig-zag path of a brane tiling is a sequence of edges turning alternatively maximally right and maximally left at each node of the toric diagram.
\end{definition}

\begin{figure}\caption{Zig-zag paths and the fast inverse algorithm}\label{zig-zag}
\centerline{\begin{tikzpicture}
\fill[gray!20] (0,1) -- (0,0) -- (1,0) -- (1+0.25,1) -- (0,1);
\fill[gray!20] (2-0.1,1) -- (2,0) -- (3,0) -- (3-0.2,1) -- (2-0.1,1);
\fill[gray!20] (4,1) -- (4,0) -- (5,0) -- (5,1) -- (4,1);
\fill[gray!80] (1-0.25,-1) -- (1,0) -- (2,0) -- (2+0.1,-1) -- (1-0.25,-1);
\fill[gray!80] (3+0.2,-1) -- (3,0) -- (4,0) -- (4,-1) -- (3+0.2,-1);
\draw[->,>=latex] (0,0) -- (5.5,0);
\draw[->,>=latex] (1-0.25,-1) -- (1+0.25,1);
\draw[->,>=latex] (2-0.1,1) -- (2+0.1,-1);
\draw[->,>=latex] (3+0.2,-1) -- (3-0.2,1);
\draw[->,>=latex] (4,1) -- (4,-1);
\node[draw,circle,inner sep=2pt,fill,white] at (0.5,0.5){};
\node[draw,circle,inner sep=2pt] (W1) at (0.5,0.5){};
\node[draw,circle,inner sep=2pt,fill,white] at (2.5,0.5){};
\node[draw,circle,inner sep=2pt] (W3) at (2.5,0.5){};
\node[draw,circle,inner sep=2pt,fill,white] at (4.5,0.5){};
\node[draw,circle,inner sep=2pt] (W5) at (4.5,0.5){};
\node[draw,circle,inner sep=2pt,fill] (B2) at (1.5,-0.5){};
\node[draw,circle,inner sep=2pt,fill] (B4) at (3.5,-0.5){};
\draw[dashed] (W1) -- (B2);
\draw[dashed] (W3) -- (B2);
\draw[dashed] (W3) -- (B4);
\draw[dashed] (W5) -- (B4);
\draw[->,>=latex] (1.3,0.3) -- (0.7,-0.3);
\draw[->,>=latex] (2.3,-0.3) -- (1.7,0.3);
\draw[->,>=latex] (3.3,0.3) -- (2.7,-0.3);
\draw[->,>=latex] (4.3,-0.3) -- (3.7,0.3);
\node at (5.5,0.5) {strip $k$};
\node at (5.5,-0.5) {strip $k+1$};
\end{tikzpicture}}
\end{figure}
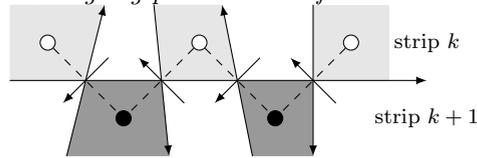

The set of edges intersecting one of the $K_z$ lines with direction $l_z$ forms a zig-zag path, following the general picture of Figure \ref{zig-zag}. These zig-zag paths divide then the torus into $K_z$ parallel strips.\medskip

As in the literature about toric quivers, we will consider general brane tiling on the torus with a consistency conditions. The consistency condition can then be expressed as the existence of an R-charge as in \cite[Def 2.4]{MozI}, or equivalently by the following conditions \cite[Def 5.1]{Ishii2009DimerMA} on zig-zag paths:
\begin{itemize}
    \item There is no homologically trivial zigzag path.
    \item No zigzag path has a self-intersection on the universal cover.
    \item No pair of zigzag paths on the universal cover intersect each other in the same direction more than once.
\end{itemize}
Not every consistent brane tiling come from the fast inverse algorithm, but one can associate to a such a brane tiling a toric diagram by considering its perfect matching as described in the next subsection. The Jacobian algebra of the quiver with potential associated to this brane tiling gives then a noncommutative crepant resolution of the corresponding toric Calabi-Yau threefold, and all the brane tiling associated to a toric diagram are related by toric mutations. There are then still $K_z$ parallel zig-zag paths with homology $l_z$ associated to a side $z$ of the toric diagram, dividing the torus into $K_z$ strips, which we can label by $k\in\mathbb{Z}/K_z\mathbb{Z}$.

The cyclic ordering is given by the orientation in the Figure \ref{zig-zag}, \ie the $k$-th strip lies to the left of the zig-zag path and the $k+1$-th strip lies to the right. We call $\Zig_k$ (resp. $\Zag_k$) the set of arrows crossing the zig-zag path, going from the $k$-th strip to the $k+1$-th strip (resp. from $k+1$-th strip to the $k$-th strip). We denote by $\alpha^z_k$ the dimension vector with component $1$ on the nodes inside the $k$-th strip, and $0$ on the other nodes. We further define
\begin{align}
\label{defalphakk}
    \alpha^z_{[k,k'[}=\alpha^z_k+\alpha^z_{k+1}+...+\alpha^z_{k'-1}
\end{align}
keeping in mind that the index $k$ lives in $\mathbb{Z}/K_z\mathbb{Z}$. In particular, $\alpha^z_{[k,k[}=\delta$ is the dimension vector with entries $1$ on each node of $Q_0$, associated to points on $X$.

\subsubsection{Perfect matchings and lattices of paths\label{sec_FFA}}\medskip

Following \cite[sec 2.2]{MozI}, consider the complex of Abelian groups:
\begin{align}
    \mathbb{Z}^{Q_2}\overset{d_2}{\to}\mathbb{Z}^{Q_1}\overset{d_1}{\to}\mathbb{Z}^{Q_0}
\end{align}
such that $d_2(F)=\sum_{a\in F}a$ and $d_1(a)=t(a)-s(a)$. We define:
\begin{align}
    \Lambda=\mathbb{Z}^{Q_1}/\langle d_2(F)-d_2(G)|F,G\in Q_2\rangle
\end{align}
and denotes by $\kappa\in\Lambda$ the image of $d_2(F)$ in the quotient $\Lambda$ for any $F\in Q_2$. The lattice $\Lambda$ (resp. its quotient $\Lambda/\mathbb{Z}\kappa=\mathbb{Z}^{Q_1}/d_2(\mathbb{Z}^{Q_2})$) is then the character lattice of the maximal torus scaling the arrows of $Q$ by leaving the potential $W$ equivariant (resp. invariant), and $\kappa$ is the $\Lambda$ weight of the potential. According to \cite[Prop 4.8]{MR}, two paths with the same source agree in $J_{Q,W}$ if and only if they have the same $\Lambda$-weight.\medskip

The map $d_1$ descends to a map $d_1:\Lambda\to\mathbb{Z}^{Q_0}$, and we define $M=ker(d)$: $M$ (resp. $M/\mathbb{Z}\kappa$) is the sub-lattice of $\Lambda$ (resp. $\Lambda/\mathbb{Z}\kappa$) giving the weights of cycles of $Q$, \ie $M$ (resp. $M/\mathbb{Z}\kappa$) gives the weight lattice of the quotient of the maximal torus scaling the arrows of the quiver leaving the potential equivariant (resp. invariant) by the gauge torus $T_G$ scaling the nodes of the quiver.\medskip

\begin{definition}
    A perfect matching is a subset $I$ of the edges of the brane tiling such that each node of the brane tiling is adjacent to exactly one edge of $I$. By duality, a perfect matching is equivalent to a cut $I$ of the quiver with potential $(Q,W)$, \ie a subset of $Q_1$ such that each cycle $w_F$ of the potential $W$ contains exactly one arrow of $I$.
\end{definition}

We define the linear map $\chi_I:\mathbb{Z}^{Q_1}\to\mathbb{Z}$ sending $a\in Q_1$ to $1$ if $a\in I$ and $0$ either. Since 
  $\chi_I(d_2(F))=1$ for $F\in Q_2$ by definition of a perfect matching, $\chi_I$ descends to a map $\chi_I:\Lambda\to\mathbb{Z}$ such that $\chi_I(\kappa)=1$, and restricts to $\bar{\chi_I}\in M^{\vee}$. Let $\sigma\in M^{\vee}_\mathbb{Q}$ be the cone generated by the $\bar{\chi_I}$. According to \cite[Remark 4.16]{MozI}, $\sigma$ gives then the fan of $X$, and the intersection of $\sigma$ with the hyperplane $\{f\in M^\vee_\mathbb{Q}|f(\kappa)=1\}$ gives the toric diagram of $X$: in particular, the lattice $L^{\vee}$ of the toric diagram is identified with $(M/\mathbb{Z}\kappa)^{\vee}$. The lattice $L$ of the brane tiling torus can then be identified as $L=M/\mathbb{Z}\kappa$.\medskip

As was first noticed in \cite[sec 4.2]{HanVegh}, $\bar{\chi_I}$ gives a node of the toric diagram: the map sending a perfect matching to the corresponding node of the toric diagram is surjective but not injective in general. However, there is a unique perfect matching associated to any corner of the toric diagram. We shall consider only such perfect matchings, and denote by $I_i$ the cut associated to the corner $p_i$.\medskip

When the two perfect matchings correspond to two adjacent corners $p_i$ and $p_{i+1}$ that are endpoints of the same side $z=z_{i+1/2}$, their union gives the zig-zag paths with direction $l_z\in L$. Removing the arrows of the two cuts $I_i,I_{i+1}$, one obtains a quiver which is a union of connected parts supported on the $K_z$ strip separated by the zig-zag paths: we denote by $Q^k$ the quiver supported on the $k$-th strip. 
We can then distinguish four types of arrows:
\begin{itemize}
    \item Arrows that are not in any cuts $I_{i}$, $I_{i+1}$ are the arrows of one connected part $Q^k$ of the remaining quiver for a $k\in\mathbb{Z}/K_z\mathbb{Z}$.
    \item Arrows that are in the intersection of the two cuts lie outside the zig-zag paths, \ie they connect nodes inside the same connected component $Q^k$, for a $k\in\mathbb{Z}/K_z\mathbb{Z}$; we denote the set of those arrows by $J_k$.
    \item Arrows in $I_{i}-I_{i+1}$ lie inside zig zag paths. With our conventions, they go from $Q^{k}$ to $Q^{k+1}$, for a $k\in\mathbb{Z}/K_z\mathbb{Z}$, \ie they forms to the above defined set $\Zig_k$
    \item Arrows in $I_{i+1}-I_{i}$ are in zig zag paths. With our conventions, they go from $Q^{k+1}$ to $Q^k$ for a $k\in\mathbb{Z}/K_z\mathbb{Z}$; \ie they forms the above defined set $\Zag_k$.
\end{itemize}

Let us denote by $M^+\subset M$ the semigroup generated by weights of cycles of $Q$. According to \cite[Cor 3.3, Cor 3.6]{MozI}, $M^+_\mathbb{Q}$ is a cone which is the dual cone of $\sigma$, and $M^+$ is saturated, \ie:
\begin{align}
    M^+=\{\lambda\in M|\chi_I(\lambda)\geq 0\quad\forall\; I\}
\end{align}
We denote by $\mathbb{C}[M^+]$ the ring generated by elements $v^\lambda$ for $\lambda\in M^+$, with relations $v^{\lambda+\lambda'}=v^\lambda v^{\lambda'}$. Because $\sigma$ is the fan of $X$, on has then $X=Spec(\mathbb{C}[M^+])$. By associating to $v^\lambda\in \mathbb{C}[M^+]$ the sum over $i\in Q_0$ of the cycles $v^\lambda_i$ of weight $\lambda$ with source and target $i$ (recall that two paths of $Q$ agree in $J_{Q,W}$ if they have the same source and $\Lambda$-weight), one obtains an inclusion $\mathbb{C}[M^+]\to J$. It was then proven in \cite[Theo 1.4]{Bro} that $\mathbb{C}[M^+]$ is the center of $J_{Q,W}$. According to \cite[Porp 3.13]{MozI}, $J_{Q,W}$ provides then a noncommutative crepant resolution of the coordinate ring of $X$.\medskip

The edges of the cone $M^+_\mathbb{Q}=\sigma^{\vee}$ are dual to sides of the toric diagram. Consider a side $z_{i+1/2}$ between the corners $p_i$ and $p_{i+1}$: the corresponding edge of $M^+_\mathbb{Q}$ lies in the intersection $\chi_{I_i}^{-1}(0)\cap \chi_{I_{i+1}}^{-1}(0)$, \ie is generated by cycles of $Q$ without arrows of $I_\cup I_{i+1}$. This shows that all the indecomposable cycles of the quivers $Q^k$ have the same $M$-weight (and equivalently the same $\Lambda$ weight) denoted by $\lambda_z$. In particular, by construction, the projection of $\lambda_z$ onto $L=M/\mathbb{Z}\kappa$ is $l_z$. We use then the notation $v^z:=v^{\lambda_z}\in J$. We have then the commutation relation, for any path $(w:i\to j)\in J$
\begin{align}\label{vicom}
    wv^z_i=v^z_jw
\end{align}

\subsubsection{Examples}

We illustrate our notations on several examples:\medskip

\begin{example}[$PdP_{3a}$]

\begin{figure}
\caption{Toric diagram and brane tiling for $PdP_{3a}$}\label{PdP3a}
\centerline{\begin{tikzpicture}
\begin{scope}[scale=0.9]
\filldraw [black] (-1,0) circle (2pt);
\draw (-1,0) node[left]{$p_1$};
\filldraw [black] (0,0) circle (2pt);
\filldraw [black] (1,0) circle (2pt);
\filldraw [black] (2,0) circle (2pt);
\draw (2,0) node[below]{$p_2$};
\filldraw [black] (0,-1) circle (2pt); 
\filldraw [black] (1,-1) circle (2pt);
\filldraw [black] (0,-2) circle (2pt); 
\draw (0,-2) node[right]{$p_0$};
\draw (-1,0) -- (2,0) -- (0,-2) -- (-1,0); 
\draw[->,>=latex] (-0.5,0) -- (-0.5,1);
\node at (-0.5,1) [above] {$l_{z_{3/2}}$};
\draw[->,>=latex] (0.5,0) -- (0.5,1);
\node at (0.5,1) [above] {$l_{z_{3/2}}$};
\draw[->,>=latex] (1.5,0) -- (1.5,1);
\node at (1.5,1) [above] {$l_{z_{3/2}}$};
\draw[->,>=latex] (-0.5,-1) -- (-2.5,-2);
\node at (-2.5,-2) [below] {$l_{z_{1/2}}$};
\draw[->,>=latex] (1.5,-0.5) -- (2.5,-1.5);
\node at (2.5,-1.5) [right] {$l_{z_{5/2}}$};
\draw[->,>=latex] (0.5,-1.5) -- (1.5,-2.5);
\node at (1.5,-2.5) [right] {$l_{z_{5/2}}$};
\end{scope}

\begin{scope}[xshift=5cm, yshift=-1.7cm,scale=0.8]
 \foreach \x in {0,...,3}{
      \foreach \y in {0,...,2}{
        \node[draw,circle,inner sep=2pt,fill] at (0+1.5*\x,0+1*\y){};
        \node[draw,circle,inner sep=2pt] at (0.5+1.5*\x,0+1*\y){};
        \node[draw,circle,inner sep=2pt,fill] at (0.75+1.5*\x,0.5+1*\y){};
        \node[draw,circle,inner sep=2pt] at (-0.25+1.5*\x,0.5+1*\y){};
        \draw (0+1.5*\x,0+1*\y) -- (0.5+1.5*\x,0+1*\y)[blue];
        \draw (0.75+1.5*\x,0.5+1*\y) -- (1.25+1.5*\x,0.5+1*\y)[blue];
        \draw (0+1.5*\x,0+1*\y) -- (-0.25+1.5*\x,0.5+1*\y)[green];
        \draw (0.75+1.5*\x,0.5+1*\y) -- (0.5+1.5*\x,1+1*\y)[green];
        \draw (0.5+1.5*\x,0+1*\y) -- (0.75+1.5*\x,0.5+1*\y)[red];
        \draw (-0.25+1.5*\x,0.5+1*\y) -- (0+1.5*\x,1+1*\y)[red];
        
        \MULTIPLY{3}{\y}{\s}
        \ADD{\s}{\x}{\so};
        \ADD{\so}{2}{\soo};
        \MODULO{\so}{6}{\sol};
        \MODULO{\soo}{6}{\sool};
        \node at (0.25+1.5*\x,0.5+1*\y){\sol};
        \node at (1+1.5*\x,0+1*\y){\sool};
        }
        \draw[->,>=latex,gray!80] (-0.125+1.5*\x,-0.5) -- (-0.125+1.5*\x,3.5);
        \draw[->,>=latex,gray!80] (-0.125+0.75+1.5*\x,-0.5) -- (-0.125+0.75+1.5*\x,3.5);
        \draw[->,>=latex,gray!80] (4.375+0.75+1.5*\x-1.5,2.75+0.5) -- (-0.125-0.75+1.5*\x-1.5,-0.75);
        \draw[->,>=latex,gray!80] (-0.125-1.5+1.5*\x-1.5,2.75+0.5) -- (4.375+1.5*\x-1.5,-0.75);
        }
\end{scope}
\end{tikzpicture}}
\end{figure}

\medskip
The toric diagram and brane tiling are given by Figure \ref{PdP3a}, where we have drawn the perfect matchings corresponding to the corners $p_0,p_1,p_2$ in blue, red and green, respectively. An arrow of the cut $I_i$ with source $j$ and target $k$ will be denoted by $\Phi^i_{jk}$.
\medskip

The zig-zag paths defined by taking the union of two consecutive perfect matchings on the boundary of the toric diagram are as follows:
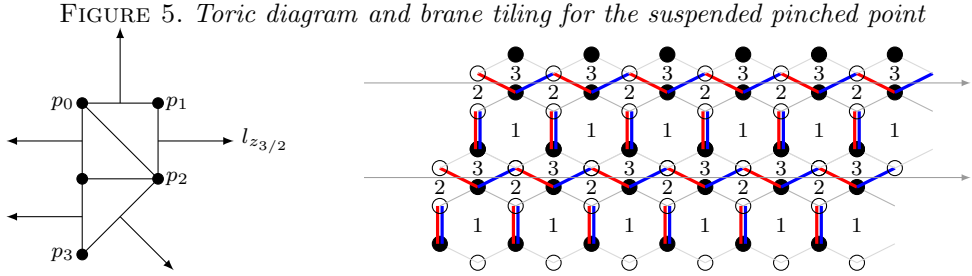
\begin{itemize}
    \item $I_0\cup I_1$: the corresponding zig-zag path corresponding to the side $z_{1/2}$ is given by the succession of blue and red edges. The remaining quiver has one connected component $Q^0$, \ie $K_z=1$ (corresponding to the fact that the associated edge of the toric diagram has one subdivision). It is a simple cyclic quiver with six nodes in the order $(0,1,2,3,4,5)$. We have then  $Q_1^0=I_2$, $J_0=\emptyset$, $\Zig_0=I_0$ and $\Zag_0=I_1$ and $\alpha^{z_{1/2}}_0=\delta$.
    
    \item $I_1\cup I_2$: the corresponding zig-zag paths corresponding to the side $z_{3/2}$ are given by the succession of red and green edges. The remaining quiver has three connected components $Q^0$, $Q^1$ and $Q^2$, \ie $K_z=3$ (corresponding to the fact that the associated edge of the toric diagram has three subdivisions). $Q^0$, $Q^1$ and $Q^2$ are simple two cycles with nodes respectively $(0,3)$, $(1,4)$ and $(2,5)$. We have:
    \begin{align}
        Q^0_0&=\{0,3\} &Q^1_0&=\{1,4\}& Q^2_0&=\{2,5\}\nn\\
        Q^0_1&=\{\Phi^0_{03},\Phi^0_{30}\}& Q^1_1&=\{\Phi^0_{14},\Phi^0_{41}\}& Q^2_1&=\{\Phi^0_{25},\Phi^0_{52}\}\nn\\
        J_0&=\emptyset &J_1&=\emptyset &J_2&=\emptyset\nn\\
        \Zig_0&=\{\Phi^1_{13},\Phi^1_{40}\}&\Zig_1&=\{\Phi^1_{24},\Phi^1_{51}\}&\Zig_2&=\{\Phi^1_{02},\Phi^1_{35}\}\nn\\
        \Zag_0&=\{\Phi^2_{01},\Phi^2_{34}\}&\Zag_1&=\{\Phi^2_{12},\Phi^2_{45}\}&\Zag_2&=\{\Phi^2_{23},\Phi^2_{50}\}\nn\\
        \alpha^{z_{3/2}}_0&=e_0+e_3&\alpha^{z_{3/2}}_1&=e_1+e_4&\alpha^{z_{3/2}}_2&=e_2+e_5
    \end{align}
    
    \item $I_2\cup I_0$: the corresponding zig-zag paths corresponding to the side $z_{5/2}$ are given by the succession of green and blue edges. The remaining quiver has two connected components $Q^0$ and $Q^1$, \ie $K_z=2$ (corresponding to the fact that the associated edge of the toric diagram has two subdivisions). $Q^0$ and $Q^1$ are respectively simple three cycles with nodes respectively $0,2,4$ and $1,3,5$. We have:
    \begin{align}
        Q^0_0&=\{0,2,4\},\quad &Q^1_0&=\{1,3,5\}\nn\\
        Q^0_1&=\{\Phi^1_{02},\Phi^1_{24},\Phi^1_{40}\}&Q^1_1&=\{\Phi^1_{13},\Phi^1_{35},\Phi^1_{51}\}\nn\\
        J_0&=\emptyset &J_1&=\emptyset\nn\\
        \Zig_0&=\{\Phi^2_{12},\Phi^2_{34},\Phi^2_{50}\}&\Zig_1&=\{\Phi^2_{01},\Phi^2_{23},\Phi^2_{45}\}\nn\\
        \Zag_0&=\{\Phi^0_{03},\Phi^0_{25},\Phi^0_{41}\},\quad&\Zag_1&=\{\Phi^0_{14},\Phi^0_{30},\Phi^0_{52}\}\nn\\
        \alpha^{z_{5/2}}_0&=e_0+e_2+e_4&\alpha^{z_{5/2}}_1&=e_1+e_3+e_5
    \end{align}
\end{itemize}
$\Box$
\end{example}\medskip

\begin{example}[Suspended pinched point]

\begin{figure}
\caption{Toric diagram and brane tiling for the suspended pinched point}\label{SPP}
\centerline{\hfill
\begin{tikzpicture} 
\begin{scope}[rotate=-90]
\filldraw [black] (0,0) circle (2pt);
\draw (0,0) node[left]{$p_0$};
\filldraw [black] (1,0) circle (2pt);
\filldraw [black] (2,0) circle (2pt);
\draw (2,0) node[left]{$p_3$}; 
\filldraw [black] (0,1) circle (2pt);
\draw (0,1) node[right]{$p_1$}; 
\filldraw [black] (1,1) circle (2pt);
\draw (1,1) node[right]{$p_2$};
\draw (0,0) -- (2,0); 
\draw (2,0) -- (1,1);
\draw (0,1) -- (1,1);
\draw (0,1) -- (0,0);
\draw (1,1) -- (0,0);
\draw (1,1) -- (1,0);
\draw[->,>=latex] (0.5,0) -- (0.5,-1);
\draw[->,>=latex] (1.5,0) -- (1.5,-1);
\SQUAREROOT{2}{\s}
\draw[->,>=latex] (1.5,0.5) -- (1.5+1/\s,0.5+1/\s);
\draw[->,>=latex] (0.5,1) -- (0.5,2);
\draw (0.5,2) node[right]{$l_{z_{3/2}}$};
\draw[->,>=latex] (0,0.5) -- (-1,0.5);
\end{scope}
\end{tikzpicture}
\hfill
\begin{tikzpicture}
 \foreach \y in {0,1}{
 \foreach \x in {0,...,5}{
    \node[draw,circle,inner sep=2pt] at (\x+0.5+0.5*\y,-0.25+1.25*\y){};
    \draw (\x+0.5+0.5*\y,-0.25+1.25*\y) -- (\x+0.5*\y,1.25*\y)[gray!30];
    \draw (\x+0.5+0.5*\y,-0.25+1.25*\y) -- (\x+1+0.5*\y,1.25*\y)[gray!30];
    \node[draw,circle,inner sep=2pt,fill] at (\x+0.5*\y,0+1.25*\y){};
    \draw (\x-0.03+0.5*\y,0+1.25*\y) -- (\x+0.5*\y-0.03,0.5+1.25*\y)[very thick,red];
    \draw (\x++0.03+0.5*\y,0+1.25*\y) -- (\x+0.5*\y+0.03,0.5+1.25*\y)[very thick,blue];
    \node[draw,circle,inner sep=2pt] at (\x+0.5*\y,0.5+1.25*\y){};
    \draw (\x+0.5*\y,0.5+1.25*\y) -- (\x+0.5+0.5*\y,0.75+1.25*\y)[gray!60];
    \draw (\x+1+0.5*\y,0.5+1.25*\y) -- (\x+0.5+0.5*\y,0.75+1.25*\y)[gray!60];
    \node[draw,circle,inner sep=2pt,fill] at (\x+0.5+0.5*\y,0.75+1.25*\y){};
    \draw (\x+0.5+0.5*\y,0.75+1.25*\y) -- (\x+0.5*\y,1+1.25*\y)[very thick,red];
    \draw (\x+0.5+0.5*\y,0.75+1.25*\y) -- (\x+1+0.5*\y,1+1.25*\y)[very thick,blue];
    \node[draw,circle,inner sep=2pt] at (\x+0.5*\y,1+1.25*\y){};
    \draw (\x+0.5*\y,1+1.25*\y) -- (\x+0.5+0.5*\y,1.25+1.25*\y)[gray!30];
    \draw (\x+1+0.5*\y,1+1.25*\y) -- (\x+0.5+0.5*\y,1.25+1.25*\y)[gray!30];
    \node[draw,circle,inner sep=2pt,fill] at (\x+0.5+0.5*\y,1.25+1.25*\y){};
    \node at (\x+0.5*\y+0.5,0.25+1.25*\y) {$1$};
    \node at (\x+0.5*\y,0.75+1.25*\y) {$2$};
    \node at (\x+0.5*\y+0.5,1+1.25*\y) {$3$};
    }}
    \draw[->,>=latex,very thin,gray!80] (-1,0.875) -- (7,0.875);
    \draw[->,>=latex,very thin,gray!80] (-1,0.875+1.25) -- (7,0.875+1.25);
\end{tikzpicture}
\hfill}
\end{figure}

For The suspended pinched point, one resolution of the toric diagram, and the corresponding brane tiling, are given by Figure \ref{SPP}. Here we have drawn  the edges of the perfect matching $p_1$ in red, and the edges of the perfect matching $p_2$ in blue. The union of these perfect matchings describes zig-zag paths corresponding to the side $z_{3/2}$ of the toric diagram. It divides the brane tiling into strips oriented from the left to the right. In particular, the quiver obtained after the two cuts has one connected component $Q^0$ (corresponding to the fact that the corresponding edge of the toric diagram has one subdivision). We have then:
\begin{align}
    &Q^0_0=\{1,2,3\},\quad Q^0_1=\{\Phi_{12},\Phi_{13},\Phi_{21},\Phi_{31}\}\nn\\
    &J_0=\{\Phi_{11}\},\quad\Zig_0=\{\Phi_{32}\},\quad\Zag_0=\{\Phi_{23}\},\quad\alpha^{z_{3/2}}_0=\delta\nn\\
    &v^z_1=\Phi_{21}\Phi_{12}=\Phi_{31}\Phi_{13},\quad v^z_2=\Phi_{12}\Phi_{21},\quad v^z_3=\Phi_{13}\Phi_{31}
\end{align}
$\Box$
\end{example}

\subsection{Framed quivers associated to toric threefolds}

\subsubsection{D6-brane framing}\label{D6def}

We introduce a first type of framed quiver built from an unframed quiver $(Q,W)$ coming from a brane tiling. Choosing $i\in Q_0$ a node of the quiver, we consider the framed quiver $Q_i$ with a framing node $\infty$, and an arrow $q:\infty\to i$, \ie $((Q_i)_0,(Q_i)_1)=(Q_0\cup\{\infty\},Q_1\cup\{q\})$. The potential is still $W$, because there is no cycle passing by the framing node. We will consider $i$-cyclic representations, \ie representations $V$ of the framed quiver $Q_i$ such that $d_\infty=1$, and the sub-representation generated by $V_\infty$ is the whole representation. We denote by $Z_i(x)$ the generating series of the cohomological DT invariants $[\mathcal{M}_{Q_i,W,d}]^{vir}$ of $i$-cyclic critical representations, following the definitions in \eqref{critcoho} and \eqref{defgenser}.
\medskip

\begin{remark}
    Such a framing corresponds to adding a D6-brane in physics terminology. Framed $i$-cyclic representations are a noncommutative analogue of sheaves with compact support on $X$ with a framing by the sheaf $\mathcal{O}_X$: such a complex is then considered as a bound state of a D6 noncompact brane (\ie a sheaf with support on the whole noncompact threefold $X$) with a D4-D2-D0 compact brane (\ie a sheaf with compact support on 2 dimensional, 1 dimensional and 0 dimensional sub-varieties).
\end{remark}

Consider a Serre subcategory $S$ of the Abelian category of representations of $Q$. There is a general formula, which is a variant of the wall crossing formula of \cite{KonSol10}, expressing the framed generating series $Z_i^S(x)$ in terms of the generating series $\mathcal{A}^S(x)$ of representations of the unframed quiver $(Q,W)$, developed in \cite{Mor11},\cite{Moz11} and \cite{MMNS}; 
\begin{align}\label{framunfram}
    Z_i^S(x)=S_i(\mathcal{A}^S(x)) S_{-i}(\mathcal{A}^S(x)^{-1})
\end{align}

\subsubsection{D4-brane framing}\label{D4def}

We consider a cut $I$ corresponding to a corner $p_i$ of the toric diagram, denoting $D$ the corresponding divisor. The divisor $D$ is in particular noncompact. We now introduce, following \cite[sec 3.2]{TYY}, a framed quiver with potential $(Q_D,W_D)$, such that $D$-cyclic representations are a noncommutative analogue of sheaves with compact support on $X$ with a framing by the sheaf $\mathcal{O}_D$. In physics terminology, 
such framed sheaves correspond to bound states of a noncompact D4 brane wrapped on $D$, together with compact D4-D2-D0 branes.\medskip

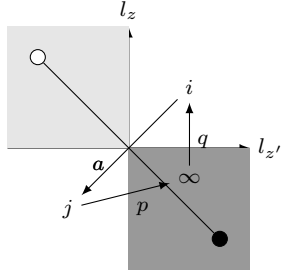
\begin{figure}
\caption{D4 brane framing and the fast inverse algorithm}\label{figd4}
\centerline{
\begin{tikzpicture}[scale=0.8]
\draw[->,>=latex] (0,-2) -- (0,2);
\draw[->,>=latex] (-2,0) -- (2,0);
\fill[gray!20] (0,0) -- (0,2) -- (-2,2) -- (-2,0) -- cycle;
\fill[gray!80] (0,0) -- (2,0) -- (2,-2) -- (0,-2) -- cycle;
\draw (0,2) node [above]{$l_z$};
\draw (2,0) node [right]{$l_{z'}$};
\node[draw,circle,inner sep=2pt,fill] (A) at (1.5,-1.5){};
\node[draw,circle,inner sep=2pt,fill,white] (B) at (-1.5,1.5){};
\node[draw,circle,inner sep=2pt] at (-1.5,1.5){};
\draw (A) -- (B); 
\node (C) at (1,1){$i$};
\node (D) at (-1,-1){$j$};
\node (E) at (1,-0.5){$\infty$};
\draw (-0.5,-0.5) node [above]{$a$};
\draw[->,>=latex] (C) -- (D);
\draw[->,>=latex] (D) -- (E);
\draw[->,>=latex] (E) -- (C);
\draw (-0.5,-0.5) node [above]{$a$};
\draw (0,-1) node [right]{$p$};
\draw (1,0.1) node [right]{$q$};
\end{tikzpicture}}
\end{figure}

The corner $p_i$ lies between the two sides $z=z_{i-1/2}$ and $z'=z_{i+1/2}$, with  $K_z$ and $K_z'$ subdivisions, respectively. We can then choose one of the intersection points of the $K_z,K_{z'}$ oriented lines on the torus with direction $l_z,l_{z'}$ (according to \cite[sec 4.4]{TYY}, different intersection points correspond to different choices for the holonomy of the gauge fields at infinity). Following the general procedure of the fast inverse algorithm, the picture at the intersection point is given by Figure \ref{figd4}.

Here we have two tiles of the brane tiling, the tile corresponding to the node $i$ of the quiver in the quadrant $(+z,+z')$, the tile corresponding to the node $j$ of the quiver in the quadrant $(-z,-z')$, and an edge  between those tiles, corresponding to an arrow $a:i\to j$ of the quiver. The corresponding framed quiver, which we denote by  $Q_D$, has one framing node $\infty$ and two framing arrows $q:\infty\to i$ and $p:j\to\infty$, \ie $((Q_{D})_0,(Q_D)_1)=(Q_0\cup\{\infty\},Q_1\cup\{p,q\}$. 
The potential for the frame quiver is obtained by adding the cycle $paq$ to the original unframed potential,
\begin{align}
    W_D=W+paq
\end{align}
We denote by $Z_D(x)$ the generating series of the cohomological DT invariants $[\mathcal{M}_{Q_D,W_D,d}]^{vir}$ of $D$-cyclic critical representations, following the definitions in \eqref{critcoho} and \eqref{defgenser}. To our knowledge, there is no known simple expression of the generating series $Z_D(x)$ in terms of the unframed generating series $\mathcal{A}(x)$ similar to the formula \eqref{framunfram} expressing $Z_i(x)$ in terms of $\mathcal{A}(x)$.
\medskip

Some general properties of $D$-cyclic critical representations are proven in \cite{TYY}. First,  in \cite[sec 3.7]{TYY} it is proven that the arrow $p$ always vanishes in such representations. Taking the partial derivative $\partial_pW_D=aq$, the arrow $p$ gives the relation $aq=0$. Second, in section 3.8 it is shown that this relation imposes that in fact all the arrows of the cut $I$ vanish. This shows that the $\mathbb{C}Q_D/(\partial W_D)$ module of paths with source at a framing node $P_D$ is generated by paths beginning by the framing node $q$, followed by a paths of the quiver with relation $(Q_I,\partial_IW)$ obtained from $Q$ by removing the arrows of $I$ and imposing the relations $\partial_aW=0$ for $a\in I$.
\medskip

In the periodic quiver plane, the paths of $(Q_I,W_I)$ beginning at the node $i$ extend in the facet between the two half lines directed by $l_{z}$, $l_{z'}$ intersecting at $i$. Indeed, the $\Zag$ arrows of the zig-zag paths associated to $z$ are in $I$, preventing paths of $Q'$ to cross the half line directed by $l_z$, and the $\Zig$ arrows of the zig-zag paths associated to $z'$ are in $I$, preventing paths of $(Q_I,W_I)$ to cross the half line directed by $l_{z'}$.

\begin{figure}
\caption{Facet corresponding to a D4 brane framing}\label{conifold}
\centerline{
\begin{tikzpicture}[scale=0.6]
 \fill[gray!40] (2.5,1) -- (-1.5,5) -- (6.5,5);
 \foreach \x in {0,...,2}{
      \foreach \y in {0,...,4}{
        \MODULO{\y}{2}{\s};
        \ADD{\y}{1}{\so}
        \MODULO{\so}{2}{\sol};
        \draw (2*\s+2*\x,\y) -- (1+2*\x,\y);
        \draw (2*\sol+2*\x,\y) -- (1+2*\x,\y)[red];
        \draw (2*\x,\y) -- (2*\x,\y+1);
        \draw (2*\x+1,\y) -- (2*\x+1,\y+1);
        \node[draw,circle,inner sep=2pt,fill] at (\s+2*\x,\y){};
        \node[draw,circle,inner sep=2pt,fill,white] at (\sol+2*\x,\y){};
        \node[draw,circle,inner sep=2pt] at (2*\x+\sol,\y){};
        \node at (0.5+2*\x,0.5+\y){\sol};
        \node at (1.5+2*\x,0.5+\y){\s};
        }}
        \draw[->,>=latex] (4,-0.5) -- (-1.5,5);
        \node at (-1.5,5) [left]{$l_z$};
        \draw[->,>=latex] (1,-0.5) -- (6.5,5);
        \node at (6.5,5) [right]{$l_{z'}$};
\end{tikzpicture}}
\end{figure}
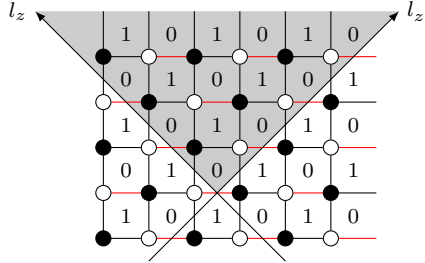
As an example, for the conifold, it gives Figure \ref{conifold}, where we have drawn  the edges of the perfect matching in red and filled in gray the nodes of the periodic quiver that are accessible from the node $0$ after having removed the arrows of the cut. \medskip

\section{Invertible and nilpotent BPS invariants}

\subsection{Definition}

Let us choose a side $z$ of the toric diagram. We have seen that the toric coordinate $v^{\lambda_z}\in\mathbb{C}[M^+]$ of $X$ is identified with the element $v^z=\sum_{i\in Q_0} v^z_i$ in the center of the Jacobian algebra $J_{Q,W}$. The noncommutative analogue of sheaves supported on the locus of $X$ where one of the toric coordinate $v^{\lambda_z}$ is non-vanishing (resp. vanishing) are critical representations $V$ where the endomorphism $v^z$ of $V$ is invertible (resp. nilpotent).\medskip

Consider a short exact sequence of $d_1$, $d=d_1+d_2$ and $d_2$ dimensional representations of $Q$:
\begin{align}
    0\to V_1\to V\to V_2\to 0
\end{align}
The operator $v^z$ is upper  triangular with respect to this block decomposition, with diagonal blocks $v^z|_{V_1}$ and $v^z|_{V_2}$, \ie $v^z$ is invertible (resp. nilpotent) in $V$ if and only it is  invertible (resp. nilpotent) in $V_1$ and $V_2$. Hence the subcategory of representations of $Q$ obtained by imposing nilpotency or invertibility to some of the $v^z$ is a Serre subcategory of the category of representations of $Q$.\medskip

For $Z_I,Z_N$ two disjoint subsets of the sides of the toric diagram, we use then the superscript $Z_I:I,Z_N:N$ to denotes the restriction to the Serre subcategories of representations where $v^z$ is invertible for $z\in Z_I$ and nilpotent for $z\in Z_N$. We use the superscript $I$ (resp. $N$) to denotes the Serre subcategory when all the cycles $v^z$ are invertible (resp nilpotent), and call the corresponding representations and generating series 'totally invertible' (resp. 'totally nilpotent'). Recalling that the sides of the toric diagram are cyclically ordered, we use also intervals notations like $[z,z']$, to denote the set the sides of the toric diagram enumerated in the clockwise order, starting at $z$ and ending at $z'$.\medskip

It is important to stress that not all results of Donaldson-Thomas theory extend to the partially invertible and nilpotent invariants. In particular, the purity result does not  hold for the BPS invariants $\Omega^{Z_I:I,Z_N:N}_{\theta,d}$, as will become apparent in the formulae of Propositions \ref{2partinv} and \ref{1partinv}.\medskip

\subsection{Invertible/Nilpotent decomposition}

Let us fix a side  $z$ of the toric diagram. Consider a representation in $\mathfrak{M}^{z:I}_{Q,W,d}$, \ie such that
the endomorphism $v^z$ is invertible. For each $k\in\mathbb{Z}/K_z\mathbb{Z}$ the connected component $Q^k$ of the quiver between the $k$-th  and the $k+1$-th zig-zag path associated to $z$ is strongly connected: for $i,j\in(Q^k)_0$, there is then paths $v:i\to j$ and $v':j\to i$ in $Q^k$, such that $v'v:i\to i$ is a cycle of $Q^k$, \ie is equal in $J_{Q,W}$ to a power of $v^z_i$. It implies that $v$ is invertible, \ie $d_i=d_j$. Then the dimension vector is constant inside a strip between two zig-zag paths associated to $z$, hence the dimension vector $d$ is in $\langle (\alpha^z_k)_k\rangle$, the linear span with positive integer coefficients of the $(\alpha^z_k)_k$. We have the useful property:\medskip

\begin{lemma}\label{lemcenter}\medskip
The dimension vectors $\alpha^z_k$ belong to the kernel of the skew-symmetrized Euler form  $\langle \cdot,\cdot\rangle$. For each $Z_I,Z_N$ such that $Z_I\neq\emptyset$, dimensions vectors supporting representations of the Serre subcategory $Z_I:I,Z_N:N$ are then in the kernel of $\langle \cdot,\cdot\rangle$, hence the corresponding BPS invariants are not subject to wall crossing, and we denote them by $\Omega^{Z_I:I,Z_N:N}(x)$, omitting the subscript $\theta$.
\end{lemma}

Proof: Consider $i\in Q^k_0$: $\langle \alpha^z_{k'},e_i\rangle$ gives the number of arrows of the quiver going from the the node $i$ to a node in $Q^{k'}_0$, minus the number of arrows of the quiver going from a node in $Q^{k'}_0$ to the node $i$. The tile $i$ is bordered by $n$ incoming and $n$ outgoing arrows. If the $k-1$-th (resp the $k$-th) zig-zag path border the tile $i$, then there is one incoming and one outgoing arrow  adjacent to a node of $Q^{k-1}_0$ (resp $Q^{k+1}_0$), and the rest of the arrows are incoming and outgoing arrows are adjacent to a node of $Q^k_0$. By disjunction of case, there is as many incoming and outgoing arrows at $i$ adjacent to a node of $Q^{k'}_0$, hence $\langle \alpha^z_{k'},e_i\rangle=0$.\medskip

Consider a representation $V\in\mathfrak{M}^{Z_I:I,Z_N:N}_{Q,d}$. Take $z\in Z_I$: $v^z$ is then invertible on $V$, \ie  $d\in ker\langle,\rangle$. In particular, the associated term $[\mathfrak{M}^{Z_I:I,Z_N:N}_{Q,W,d}]^{vir}x^d$, and then $\mathcal{A}^{Z_I:I,Z_N:N}(x)$, is in the center of the quantum affine space, \ie no wall crossing can occur, and BPS invariants do not depend on the stability parameter $\theta$.$\Box$\medskip

We now show the following Lemma, which is a direct generalization of \cite[Lem 4.1]{Davison2016TheIC} to the case of non-symmetric quivers:\medskip

\begin{lemma}\label{supportlemma}
    Consider a quiver with potential $(Q,W)$ and an element $v$ central in $J_{Q,W}$, such that representations where $v$ is invertible have a dimension vector in the kernel of the anti-symmetrized Euler form $\langle,\rangle$. Then $v$ acts as a scalar on representations in the support of $\mathcal{BPS}_{W,d}^\theta$.
\end{lemma}

Proof: We refer to the proof of the Lemma 4.1 of \cite{Davison2016TheIC} for the details of the arguments: here the major difference is that we consider a quiver which is not symmetric, and then the quantum affine space is not commutative. Considering $\theta$ generic, \ie such that if $d,d'$ have the same slope then $\langle d-d',\bullet\rangle=0$, and a ray $l$ of the form $d+ker(\langle,\rangle)$. The relative integrality Theorem from \cite[Theo A]{DavMein} gives:
\begin{align}
    \bigoplus_{d\in l}\mathcal{H}(JH_\ast\phi_W\mathcal{IC}_{\mathfrak{M}^{\theta,ss}_{Q,d}})=\Sym_{\boxtimes_\oplus}(H(B\mathbb{C}^\ast)_{vir}\otimes(\bigoplus_{d\in l}\mathcal{BPS}_{W,d}^\theta))
\end{align}
where $JH$ is the Jordan-Hölder map sending a semistable object to the associated polystable object and $\phi_W$ is the vanishing cycle functor of $\Tr(W)$ on $\mathfrak{M}^{\theta,ss}_{Q,d}$. Consider $V\in \Supp(\mathcal{BPS}_{W,d}^\theta)$. In particular, $V$ is the polystable object associated with a representation in the support of $\phi_W$, hence in the critical locus of $\Tr(W)$, and then $V$ itself is a $J_{Q,W}$-module. Suppose that $v$, which is central in $J_{Q,W}$, has at least two different eigenvalues, which we denote $\epsilon_1$ and $\epsilon_2$. We choose two disjoints open set $U_1,U_2\subset\mathbb{C}$ such that the eigenvalues of $v$ lies in $U_1\cup U_2$, $\epsilon_1\in U_1,\epsilon_2\in U_2$. Given an open set $U\subset\mathbb{C}$ we denote by $(\mathcal{M}^{\theta,ss}_{Q,d})^U$  (resp. $(\mathfrak{M}^{\theta,ss}_{Q,d})^U$) the subspace (resp. sub-stack) of representations such that $v$ has all its eigenvalues in $U$, in particular $V\in(\mathcal{M}^{\theta,ss}_{Q,d})^{U_1\cup U_2}-((\mathcal{M}^{\theta,ss}_{Q,d})^{U_1}\cup(\mathcal{M}^{\theta,ss}_{Q,d})^{U_2})$. A  critical representation $W\in \mathfrak{M}_{Q,W,d}^{U_1\cup U_2}$ splits canonically as a direct sum of representations $W_1,W_2$ where $v$ has eigenvalues respectively in $U_1,U_2$, giving:
\begin{align}
    (\mathfrak{M}^{\theta,ss}_{Q,W,d})^{U_1\cup U_2}=\bigsqcup_{d_1+d_2=d}(\mathfrak{M}^{\theta,ss}_{Q,W,d_1})^{U_1}\times (\mathfrak{M}^{\theta,ss}_{Q,W,d_2})^{U_2}
\end{align}
Remark here that necessarily at least one of the two $U_i$ is contained in $\mathbb{C}^\ast$, say $U_1$, and then contains only representations where $v$ is invertible, \ie the $d_1$ which gives nontrivial terms in the sum lies in the ray $l_0:=ker(\langle,\rangle)$, giving:
\begin{align}
    \bigsqcup_{d\in l}(\mathfrak{M}^{\theta,ss}_{Q,W,d})^{U_1\cup U_2}=\bigsqcup_{d\in l}\bigsqcup_{d_1+d_2=d}(\mathfrak{M}^{\theta,ss}_{Q,W,d_1})^{U_1}\times (\mathfrak{M}^{\theta,ss}_{Q,W,d_2})^{U_2}\nn\\
    =(\bigsqcup_{d_1\in l_0}(\mathfrak{M}^{\theta,ss}_{Q,W,d_1})^{U_1})\times(\bigsqcup_{d_2\in l}(\mathfrak{M}^{\theta,ss}_{Q,W,d_2})^{U_2})
\end{align}
By the same argument as in the proof of Lemma 4.1 of \cite{Davison2016TheIC}, it gives:
\begin{align}
    &\bigoplus_{d\in l}\mathcal{H}(JH_\ast\phi_W\mathcal{IC}_{\mathfrak{M}^{\theta,ss}_{Q,d}}|_{(\mathfrak{M}^{\theta,ss}_{Q,d})^{U_1\cup U_2}})\nn\\=&\bigoplus_{d_1\in l_0}\mathcal{H}(JH_\ast\phi_W\mathcal{IC}_{\mathfrak{M}^{\theta,ss}_{Q,d}}|_{(\mathfrak{M}^{\theta,ss}_{Q,d})^{U_1}})\boxtimes_{\oplus}\bigoplus_{d_2\in l}\mathcal{H}(JH_\ast\phi_W\mathcal{IC}_{\mathfrak{M}^{\theta,ss}_{Q,d}}|_{(\mathfrak{M}^{\theta,ss}_{Q,d})^{U_2}})
\end{align}
Applying the relative integrality Theorem of \cite{DavMein}, we obtain:
\begin{align}
    &\Sym_{\boxtimes_\oplus}(H(B\mathbb{C}^\ast)_{vir}\otimes(\bigoplus_{d\in l}\mathcal{BPS}_{W,d}^\theta|_{(\mathcal{M}^{\theta,ss}_{Q,d})^{U_1\cup U_2}})\nn\\=&\Sym_{\boxtimes_\oplus}(H(B\mathbb{C}^\ast)_{vir}\otimes(\bigoplus_{d_1\in l_0}\mathcal{BPS}_{W,d_1}^\theta|_{(\mathcal{M}^{\theta,ss}_{Q,d_1})^{U_1}}))\nn\\ &\boxtimes_{\oplus}\Sym_{\boxtimes_\oplus}(H(B\mathbb{C}^\ast)_{vir}\otimes(\bigoplus_{d_2\in l}\mathcal{BPS}_{W,d_2}^\theta|_{(\mathcal{M}^{\theta,ss}_{Q,d_2})^{U_2}}))\nn\\
    =&\Sym_{\boxtimes_\oplus}(H(B\mathbb{C}^\ast)_{vir}\otimes((\bigoplus_{d_1\in l_0}\mathcal{BPS}_{W,d_1}^\theta|_{(\mathcal{M}^{\theta,ss}_{Q,d_1})^{U_1}})\oplus(\bigoplus_{d_2\in l}\mathcal{BPS}_{W,d_2}^\theta|_{(\mathcal{M}^{\theta,ss}_{Q,d_2})^{U_2}})))
\end{align}
and then by identification one has:
\begin{align}
    \sum_{d\in l}\mathcal{BPS}_{W,d}^\theta|_{(\mathcal{M}^{\theta,ss}_{Q,d})^{U_1\cup U_2}}\simeq(\sum_{d_1\in l_0}\mathcal{BPS}_{W,d_1}^\theta|_{(\mathcal{M}^{\theta,ss}_{Q,d_1})^{U_1}})\oplus(\sum_{d_2\in l}\mathcal{BPS}_{W,d_2}^\theta|_{(\mathcal{M}^{\theta,ss}_{Q,d_2})^{U_2}})
\end{align}
We can then deduce:
\begin{align}
    \Supp(\mathcal{BPS}_{W,d}^\theta|_{(\mathcal{M}^{\theta,ss}_{Q,d})^{U_1\cup U_2}})\subset (\mathcal{M}^{\theta,ss}_{Q,d})^{U_1}\cup(\mathcal{M}^{\theta,ss}_{Q,d})^{U_2}
\end{align}
but $V\in (\mathcal{M}^{\theta,ss}_{Q,d})^{U_1\cup U_2}-((\mathcal{M}^{\theta,ss}_{Q,d})^{U_1}\cup(\mathcal{M}^{\theta,ss}_{Q,d})^{U_2})$, and so the restriction of $\mathcal{BPS}_{W,d}^\theta$ to $V$ is zero. We conclude that if $V\in \Supp(\mathcal{BPS}_{W,d}^\theta)$, then $v$ has a single eigenvalue. For a stable $J_{Q,W}$-module $W$, because $v$ is central in $J_{Q,W}$, it defines an element of $\Hom(W,W)=\mathbb{C}$ by stability, and then $v$ acts as a scalar.  $V$ is a direct sum of $\theta$-stable representations where  $v$ acts as a scalar, and $v$ has a single eigenvalue on $V$, \ie acts as a scalar on $V$. $\Box$

We can now prove the invertible/nilpotent decomposition of BPS invariants, considering $Z_I,Z_N$ two subsets of the set of sides of the toric diagram. A particular case of this result was stated and used in \cite[eq 6.18]{MP20}.\medskip

\begin{proposition}\label{nilinvdec}\medskip
    For each $Z_I,Z_N,z$ such that $z\not\in Z_I\cup Z_N$, we have:
    \begin{align}
        \Omega_\theta^{Z_I:I,Z_N:N}(x)=\Omega^{Z_I:I,Z_N,z:I}(x)+\Omega_\theta^{Z_I:I,Z_N:N,z:N}(x)
    \end{align}
\end{proposition}

Proof: Recall that the element $v^z$ is in the center of the Jacobian algebra $J_{Q,W}$, and from the Lemma \ref{lemcenter} the assumptions of Lemma \ref{supportlemma} are satisfied, \ie $v^z$ acts as a scalar on a representation in the support of $\mathcal{BPS}_{W,d}^\theta$. In particular, the support of $\mathcal{BPS}_{W,d}^\theta$ is the disjoint union of a locus where the $(v^z_i)_{i\in Q_0}$ are invertible, and a locus where they are nilpotent, \ie:
\begin{align}
    &\Supp(\mathcal{BPS}_{W,d}^\theta)\cap(\mathcal{M}^{\theta,ss}_{Q,d})^{Z_I:I,Z_N:N}\nn\\=&(\Supp(\mathcal{BPS}_{W,d}^\theta)\cap(\mathcal{M}^{\theta,ss}_{Q,d})^{Z_I:I,Z_N:N,z:I})\sqcup(\Supp(\mathcal{BPS}_{W,d}^\theta)\cap(\mathcal{M}^{\theta,ss}_{Q,d})^{Z_I:I,Z_N:N,z:N})\nn\\
    \Rightarrow &\Omega_{\theta,d}^{Z_I:I,Z_N:N}=\Omega_{\theta,d}^{Z_I:I,Z_N,z:I}+\Omega_{\theta,d}^{Z_I:I,Z_N:N,z:N}
\end{align}
where the second line holds by taking the induced long exact sequence in cohomology. The result follows by noticing that $\Omega_{\theta,d}^{Z_I:I,Z_N,z:I}$ does not depend on $\theta$, and by taking the generating series. $\Box$\medskip

\begin{remark}
The fact that, for a cycle $v^z$, the BPS invariants are the sum BPS invariants with $v^z$ invertible and BPS invariants with $v^z$ nilpotent was used in \cite[eq. 6.18]{MP20} in some specific examples of toric quivers. It was also remarked in \cite[eq 6.17]{MP20} that in those examples the invertible BPS invariants do not suffer wall crossing and are simple to compute. The formula \cite[eq. 6.18]{MP20} was proven after dimensional reduction relative to a perfect matching corresponding to a corner $p$ of the toric diagram. After this dimensional reduction, it is then possible to give an invertible/nilpotent decomposition of the cycles $v^{z},v^{z'}$ corresponding to the sides $z,z'$ adjacent to $p$, but the other cycles $v^{z''}$ vanish. Hence to provide the invertible/nilpotent decomposition on various cycles $v^z$ as in Proposition \ref{nilinvdec} one must establish this identity without doing a dimensional reduction, hence work with the formalism of vanishing cycles as done in this section.
\end{remark}

\subsection{Computation of the partially invertible part}\label{partinvsec}

\begin{proposition}\label{2partinv}
    $i)$ For any toric quiver with potential, denoting by $b$ the number of boundary points of the toric diagram and $i$ is the number of internal points of the toric diagram: 
    \begin{align}
        \Omega_{n\delta}=\mathbb{L}^{3/2}+(b-3+i)\mathbb{L}^{1/2}+i\mathbb{L}^{-1/2} for\;n\geq1
    \end{align}
    $ii)$ Consider $z,z'$ two different sides of the toric diagram, we have:
    \begin{itemize}
        \item if $z$ and $z'$ are adjacent to the same corner, then:
        \begin{align}
            \Omega^{z:I,z':I}(x)=(\mathbb{L}^{3/2}-2\mathbb{L}^{1/2}+\mathbb{L}^{-1/2})\sum_{n\geq 1}x^{n\delta}
        \end{align}
        \item otherwise:
        \begin{align}
            \Omega^{z:I,z':I}(x)=\Omega^I(x)=(\mathbb{L}^{3/2}-3\mathbb{L}^{1/2}+3\mathbb{L}^{-1/2}-\mathbb{L}^{-3/2})\sum_{n\geq 1}x^{n\delta}
        \end{align}
    \end{itemize}
\end{proposition}

Proof: $i)$ Any generic stability condition $\theta$ gives a crepant resolution $X_\theta\overset{p}{\to}X$, by taking the moduli space of $\theta$-stable $\delta$-dimensional critical representations of $(Q,W)$, as proven in \cite{MozI}[Theo 4.5]. Denote by $X^{Z:I}_\theta$ the open locus of representations such that $v^z$ is invertible for $z\in Z$, \ie $X^{Z:I}_\theta=\cap_{z\in Z}(v^{\lambda_z}p)^{-1}(\mathbb{C}^\ast)$.\medskip

As shown in \cite[Theo 4.5]{MozI}, there is an equivalence between the bounded derived category of critical representations of $(Q,W)$ and the bounded derived category of coherent sheaves on $X_\theta$. This derived equivalence restricts to a derived equivalence between the bounded derived category of critical representations such that $\sum_{i\in Q_0} v^z_i$ is invertible for $z\in Z$, and the bounded derived category of sheaves on $X^{Z:I}_\theta$. In particular, the cohomological DT/PT correspondance proven in \cite[Theo 1.1]{BriDT/PT}, and the wall crossing formula expressing the DT/PT wall crossing in terms of the BPS invariants $\Omega_{\theta,d}$, proven in \cite{Mor11} and \cite{Moz11} applies, giving:
\begin{align}
    \sum_n[(X^{Z:I}_\theta)^{[n]}]^{vir}x^{n\delta}=\Exp\left(\sum_n\Omega^{Z:I}_{\theta,n\delta}\frac{\mathbb{L}^{n/2}-\mathbb{L}^{-n/2}}{\mathbb{L}^{1/2}-\mathbb{L}^{-1/2}}x^{n\delta}\right)
\end{align}
with $(X^{Z:I}_\theta)^{[n]}$ being the Hilbert scheme of $n$ points on $X^{Z:I}_\theta$. The generating series of 
the motives of the Hilbert schemes of points on any smooth quasi-projective threefold (as $X^{Z:I}_\theta$) was computed in \cite[Theo 3.3]{BBS}:
\begin{align}
    \sum_n[(X^{Z:I}_\theta)^{[n]}]^{vir}=\Exp\left(\mathbb{L}^{-3/2}[X^{Z:I}_\theta]\sum_{n\geq 1}\frac{\mathbb{L}^{n/2}-\mathbb{L}^{-n/2}}{\mathbb{L}^{1/2}-\mathbb{L}^{-1/2}}x^{n\delta}\right)
\end{align}
\ie we can identify:
\begin{align}\label{ndelta}
    \Omega^{Z:I}_{n\delta}=\mathbb{L}^{-3/2}[X^{Z:I}_\theta]\quad for\;n\geq1
\end{align}
Here we have dropped the dependence on $\theta$  because $n\delta$ is in the kernel of the anti-symmetrized Euler form. In particular, the computation of the cohomological class $[X_\theta]$ in \cite[Lem 4.2]{MP20} gives the claimed expression for $\Omega_{n\delta}$, as explained in \cite[Remark 5.2]{MP20}.\medskip

$ii)$ Consider a subset $Z_I$ of the sides of the toric diagram containing at least two elements $z\neq z'$. For any representation in  $\mathcal{M}^{\theta,ss,Z_I:I}_{Q,W,d}$, $v^z$ and $v^{z'}$ are invertible. Consider $d$ such that $\mathcal{M}^{\theta,ss,Z_I:I}_{Q,W,d}\neq\emptyset$: the dimensions $d_i$ are constant inside the strips delimited by lines directed by $l_z$, and also inside the strips delimited by lines directed by $l_{z'}$. Since these two sets of lines intersect only at isolated points of the  torus, all the $d_i$'s must then be equal. Hence for $d\not\in\langle\delta\rangle$, $\Omega^{Z_I}_{\theta,d}$ vanish, giving:
\begin{align}
    \Omega^{Z_I:I}(x)=&\sum_{n}\Omega^{Z_I:I}_{n\delta}x^{n\delta}\nn\\
    =&\mathbb{L}^{-3/2}[X_\theta^{Z_I:I}]\sum_{n\geq 1}x^{n\delta}
\end{align}
Here we have used equation \eqref{ndelta} in the second line, considering a generic stability condition $\theta$. Recall that we have $X=Spec(\mathbb{C}[M^+])$, and then $X^{Z_I:I}=Spec(\mathbb{C}[M^+,(v^{-\lambda_z})_{z\in Z_I}]$.
There are two possible cases:
\begin{itemize}
    \item  $Z_I=\{z,z'\}$  consists of two sides of the toric diagram which are adjacent to the same corner $p$ associated to the cut $I$. In that case, the sub-semigroup of $M$ generated by $M^+$ and $-\lambda_z,-\lambda_{z'}$ is the half lattice $\{\lambda\in M|\bar{\chi_I}(\lambda)\geq 0\}$, isomorphic with $\mathbb{N}\times\mathbb{Z}^2$. This implies that $X^{Z_I:I}=(\mathbb{C}^\ast)^2\times\mathbb{C}$; $X^{Z_I:I}$, which is smooth, \ie is equal to its crepant resolution $X^{Z_I:I}_\theta$: 
    \begin{align}
        &X_\theta^{z:I,z':I}=(\mathbb{C}^\ast)^2\times\mathbb{C}\nn\\
        \Rightarrow &\Omega^{z:I,z':I}(x)=(\mathbb{L}^{3/2}-2\mathbb{L}^{1/2}+\mathbb{L}^{-1/2})\sum_{n\geq 1}x^{n\delta}
    \end{align}
    \item $Z_I$ contains two sides of the toric diagram $z,z'$ which are not on the same corner of the toric diagram. In this case the sub semigroup of $M$ generated by $L^+$ and $-\lambda_z,-\lambda_{z'}$ is the whole lattice $M$ isomorphic with $\mathbb{Z}^3$. This implies that  $X^{Z_I:I}=(\mathbb{C}^\ast)^3$, which is smooth, \ie is equal to its crepant resolution $X^{Z_I:I}_\theta$:
    \begin{align}
        &X^{Z_I:I}=X_\theta^{z:I,z':I}=X^I=(\mathbb{C}^\ast)^3\nn\\
        \Rightarrow &\Omega^{z:I,z':I}(x)=\Omega_\theta^{I}(x)=(\mathbb{L}^{3/2}-3\mathbb{L}^{1/2}+3\mathbb{L}^{-1/2}-\mathbb{L}^{-3/2})\sum_{n\geq 1}x^{n\delta}
    \end{align}
\end{itemize}
$\Box$\medskip

\begin{proposition}\label{1partinv}
Consider a side $z$ of the toric diagram:
\begin{align}
    \resizebox{1\hsize}{!}{$\Omega^{z:I}(x)=(\mathbb{L}^{3/2}+(K_z-2)\mathbb{L}^{1/2}-(K_z-1)\mathbb{L}^{-1/2})\sum_{n\geq 1}x^{n\delta}+(\mathbb{L}^{1/2}-\mathbb{L}^{-1/2})\sum_{k\neq k'}\sum_{n\geq 0}x^{n\delta+\alpha^z_{[k,k'[}}$}
\end{align}
\end{proposition}

Proof:
The main idea of the proof is to obtain an analogue of the isomorphism $X^{z:I}\cong \mathbb{C}^2/\mathbb{Z}_{K_z}\times \mathbb{C}^\ast$ at the level of the noncommutative resolution by quivers with potential. We denote by $(\bar{Q},\bar{W})$ the quiver with potential associated to the toric threefold $\mathbb{C}^2/\mathbb{Z}_{K_z}\times \mathbb{C}$: its nodes are $\bar{Q}_0=\mathbb{Z}/K_z\mathbb{Z}$, its arrows are $a_k:k\to k$, $b_k:k\to k+1$, $c_k:k+1\to k$, and its potential is:
\begin{align}
    \bar{W}=\sum_{k\in\mathbb{Z}/K_z\mathbb{Z}}(a_k c_k b_k-a_{k+1} b_k c_k)
\end{align}
We denote by $\bar{z}$ the side of the toric diagram of $\mathbb{C}^2/\mathbb{Z}_{K_z}\times \mathbb{C}$ such that $v^{\bar{z}}=\bigoplus_k a_k$. Starting with a representation $V$ of $(Q,W)$ such that $v^z$ is invertible, one can obtain a representation of $(\bar{Q},\bar{W})$ such that the $a_k$'s are invertible informally by the following procedure:
\begin{itemize}
    \item contract the strips $Q^k$ onto the single node $k$ using the fact that all the arrows of $Q^k$ are isomorphisms.
    \item Send the invertible cycles $v^z_i$, for $i\in (Q^k)_0$, to an invertible loop $a_k:k\to k$.
    \item Send the arrows of $\Zig_k:(Q^k)_0\to(Q^{k+1})_0$ to $b_k:k\to k+1$
    \item Send the arrows of $\Zag_k:(Q^{k+1})_0\to(Q^k)_0$ to $c_k:k+1\to k$
\end{itemize}
This contraction, and the corresponding operation on the toric diagrams, are illustrated in Figure \ref{contrPdP3a} for the case of the Pseudo-del Pezzo surface $PdP_{3a}$.

\begin{figure}\caption{Contraction of the quiver for $PdP_{3a}$}\label{contrPdP3a}
\centerline{\begin{tikzpicture} 
 
 \begin{scope}[scale=4/3,xshift=-1/6cm]
 \foreach \x in {0,...,2}{
      \foreach \y in {0,1}{
        \node[draw,circle,inner sep=2pt,fill] at (0+1.5*\x,0+1*\y){};
        \node[draw,circle,inner sep=2pt] at (0.5+1.5*\x,0+1*\y){};
        \node[draw,circle,inner sep=2pt,fill] at (0.75+1.5*\x,0.5+1*\y){};
        \node[draw,circle,inner sep=2pt] at (-0.25+1.5*\x,0.5+1*\y){};
        \draw (0+1.5*\x,0+1*\y) -- (0.5+1.5*\x,0+1*\y)[gray!60];
        \draw[->,>=latex,blue] (0.25+1.5*\x,-0.2+1*\y) -- (0.25+1.5*\x,0.2+1*\y);
        \draw (0.75+1.5*\x,0.5+1*\y) -- (1.25+1.5*\x,0.5+1*\y)[gray!60];
        \draw[->,>=latex,blue] (1+1.5*\x,0.3+1*\y) -- (1+1.5*\x,0.7+1*\y);
        \draw (0+1.5*\x,0+1*\y) -- (-0.25+1.5*\x,0.5+1*\y);
        \draw[->,>=latex,green] (0+1.5*\x,0.4+1*\y) -- (-0.25+1.5*\x,0.1+1*\y);
        \draw (0.75+1.5*\x,0.5+1*\y) -- (0.5+1.5*\x,1+1*\y);
        \draw[->,>=latex,green] (0.75+1.5*\x,0.9+1*\y) -- (0.5+1.5*\x,0.6+1*\y);
        \draw (0.5+1.5*\x,0+1*\y) -- (0.75+1.5*\x,0.5+1*\y);
        \draw[->,>=latex,red] (0.5+1.5*\x,0.4+1*\y) -- (0.75+1.5*\x,0.1+1*\y);
        \draw (-0.25+1.5*\x,0.5+1*\y) -- (0+1.5*\x,1+1*\y);
        \draw[->,>=latex,red] (0.5+1.5*\x-0.75,0.5+0.4+1*\y) -- (0+1.5*\x,0.5+0.1+1*\y);
        \MULTIPLY{3}{\y}{\s}
        \SUBTRACT{\s}{\x}{\so};
        \ADD{\so}{1}{\soo};
        \MODULO{\so}{6}{\sol};
        \MODULO{\soo}{6}{\sool}
        \node at (0.25+1.5*\x,0.5+1*\y) {\sol};
        \node at (1+1.5*\x,0+1*\y) {\sool};
        }}
        
        \end{scope}
        \node at  (0.4,3) {$Q^0$};
        \node at  (1.4,3) {$Q^1$};
        \node at  (2.4,3) {$Q^2$};
        \node at  (3.4,3) {$Q^0$};
        \node at  (4.4,3) {$Q^1$};
        \begin{scope}[yshift=-0.5cm]
        \node (A) at  (0.25,-1) {$0$};
        \node (B) at  (1.25,-1) {$1$};
        \node (C) at  (2.25,-1) {$2$};
        \node (D) at  (3.25,-1) {$0$};
        \node (E) at  (4.25,-1) {$1$};
        \node[red] at  (0.75,-0.5) {$b_0$};
        \node[red] at  (1.75,-0.5) {$b_1$};
        \node[red] at  (2.75,-0.5) {$b_2$};
        \node[red] at  (3.75,-0.5) {$b_0$};
        \node[green] at  (0.75,-1.5) {$c_0$};
        \node[green] at  (1.75,-1.5) {$c_1$};
        \node[green] at  (2.75,-1.5) {$c_2$};
        \node[green] at  (3.75,-1.5) {$c_0$};
        \path (A) edge [loop above,blue] node {$a_0$} (A);
        \path (B) edge [loop above,blue] node {$a_1$} (B);
        \path (C) edge [loop above,blue] node {$a_2$} (C);
        \path (D) edge [loop above,blue] node {$a_0$} (D);
        \path (E) edge [loop above,blue] node {$a_1$} (E);
        \draw[->,>=latex,red] (A) to[bend left] (B);
        \draw[->,>=latex,red] (B) to[bend left] (C);
        \draw[->,>=latex,red] (C) to[bend left] (D);
        \draw[->,>=latex,red] (D) to[bend left] (E);
        \draw[->,>=latex,green] (E) to[bend left] (D);
        \draw[->,>=latex,green] (D) to[bend left] (C);
        \draw[->,>=latex,green] (C) to[bend left] (B);
        \draw[->,>=latex,green] (B) to[bend left] (A);
        \end{scope}
        \begin{scope}[xshift=8cm,yshift=0.5cm]
\filldraw [black] (0,0) circle (2pt);
\filldraw [black] (0,1) circle (2pt);
\filldraw [black] (1,1) circle (2pt);
\filldraw [black] (-1,2) circle (2pt); 
\filldraw [black] (0,2) circle (2pt);
\filldraw [black] (1,2) circle (2pt); 
\filldraw [black] (2,2) circle (2pt);
\draw (0,0) -- (2,2) -- (-1,2) -- (0,0);
\draw[dotted] (0,1) -- (-1,2);
\draw[dotted] (0,1) -- (2,2);
\end{scope}
\begin{scope}[xshift=8cm,yshift=-3cm]
\filldraw [black] (0,1) circle (2pt);
\filldraw [black] (-1,2) circle (2pt); 
\filldraw [black] (0,2) circle (2pt);
\filldraw [black] (1,2) circle (2pt); 
\filldraw [black] (2,2) circle (2pt);
\draw (0,1) -- (2,2) -- (-1,2) -- (0,1);

\end{scope}
\end{tikzpicture}}
\end{figure}

Consider now the toric diagram associated to the toric quiver: from simple convex geometry, it contains a node which is in a line parallel with the line containing the side $z$ and which is at minimal distance from this line: for $I$ a perfect matching associated to this node, one has then $\chi_I(\lambda_z)=1$, hence each cycle $v^z_i$ contains exactly one arrow of $I$. We consider then the surjective morphism of path algebras:
\begin{align}
    \Phi:&\mathbb{C}Q\to\mathbb{C}\bar{Q}\nn\\
         &\Phi(e_i):=e_k\;\text{for}\;i\in Q^k_0\nn\\
         &\Phi(u):=a_k^{\delta_{u\in I}}\;\text{for}\;u\in Q^k_1\nn\\
         &\Phi(u):=c_kb_ka_k^{\delta_{u\in I}}\;\text{for}\;u\in J_k\nn\\
         &\Phi(u):=b_ka_k^{\delta_{u\in I}}\;\text{for}\;u\in \Zig_k\nn\\
         &\Phi(u):a_kc_k^{\delta_{u\in I}}\;\text{for}\;u\in \Zag_k
\end{align}

We evaluate then $\Phi(W)$. Consider a white (resp black) node of the perfect matching on the $k$-th zig-zag path: the corresponding cycle $w$ of the potential is of the form $cbv$ (resp $bcv$), for $v$ a path of $Q_k$ (resp $Q_{k+1}$), $b\in\Zig_k$ and $c\in\Zag_k$. If $v$ contains an arrow of $I$, then $w=c_kb_ka_k$ (resp $w=b_kc_ka_{k+1}$). Consider the cycle of $Q^k$ (resp $Q^{k+1}$) going along the $k$-th zig-zag path: it has weight $\lambda_z$, hence contains exactly one arrow of $I$, and then the $k$-th zig-zag path contains exactly one white node $w=c_kb_ka_k$ such that the arrow of $I$ adjacent to this node is in $Q^k$ and and one black node $w=b_kc_ka_{k+1}$ such that the arrow of $I$ adjacent to this node is in $Q^{k+1}$. If an arrow $b\in\Zig_k$ (resp $c\in\Zag_k$) is in $I$, then the cycles of the potential corresponding to the black and white node adjacent to this edge in the dimer model are equal to $c_kb_ka_k$, and then cancel each other. The cycles corresponding to the white and black nodes adjacent to an edge in $u\in J_k$ are both equal to $c_kb_ka_k$ and then cancel each others. Finally we obtain:
\begin{align}
    \Phi(W)=\bar{W}
\end{align}
In particular $\Phi$ pass to the quotient to a surjective morphism of Jacobi algebra:
\begin{align}
    \Phi:J_{Q,W}\to J_{\bar{Q},\bar{W}}
\end{align}
For each $k$, choose a node $i_k\in Q^k_0$, and weak paths $(u_k:i_k\to i_{k+1})=u'^{-1}bu$ (resp $(v_k:i_{k+1}\to i_k)=v'^{-1}cv$) with $u,v'$ paths of $Q^k$ and $u',v$ paths of $Q^{k+1}$, such that $\chi_I(u_k)=\chi_(v_k)=0$. We define then a morphism between localized path algebras:
\begin{align}
    \Psi:&\mathbb{C}\bar{Q}[(a_k^{-1})_k]\to\mathbb{Q}[((v^z_i)^{-1})_i]\nn\\
    &\Psi(e_k):=e_{i_k}\nn\\
    &\Psi(a_k):=v^z_{i_k}\nn\\
    &\Psi(b_k):=u_k\nn\\
    &\Psi(c_k):=v_k
\end{align}
it satisfies $\Phi\circ\Psi=\Id_{\mathbb{C}\bar{Q}[(a_k^{-1})_k]}$, and pass to the quotient to a morphism of Jacobi algebra:
\begin{align}
    \Psi:J_{\bar{Q},\bar{W}}[(a_k^{-1})_k]\to J_{Q,W}[((v^z_k)^{-1})_k]
\end{align}
which is the inverse of $\Phi$, hence $\Phi$ gives an isomorphism of localized Jacobi algebras:
\begin{align}
    \Phi:J_{Q,W}[((v^z_k)^{-1})_k]\to J_{\bar{Q},\bar{W}}[(a_k^{-1})_k]
\end{align}

For $\bar{d}\in\mathbb{N}^{\bar{Q}_0}$, denote $d=\sum_k d_k\alpha^z_k$. The morphism $\Phi$ of path algebra induces an closed embedding of stacks $\mathfrak{M}^{\bar{z}:I}_{\bar{Q},\bar{d}}\to\mathfrak{M}^{z:I}_{Q,d}$ intertwining the potential, such that the restriction to the critical locus is an isomorphism. It is then an embedding of critical charts in the sense of \cite{Joyce2013ACM}, hence:
\begin{align}
    [\mathfrak{M}^{z:I}_{Q,W,d}]^{vir}&= H_c(\mathfrak{M}^{z:I}_{Q,d},\phi_W\mathcal{IC}_{\mathfrak{M}^{z:I}_{Q,d}})\nn\\
    &=H_c(\mathfrak{M}^{\bar{z}:I}_{\bar{Q},\bar{d}},\phi_{\bar{W}}\mathcal{IC}_{\mathfrak{M}^{\bar{z}:I}_{\bar{Q},\bar{d}}})\nn\\
    &=[\mathfrak{M}^{\bar{z}:I}_{\bar{Q},\bar{W}\bar{d}}]^{vir}
\end{align}
Here the first and the last lines holds because the definition of the vanishing cycle functor is local, hence commutes with open embeddings, and the second line follows from the isomorphism of \cite[Theo 5.4]{Joycesymstab}. We obtain a relation between the generating series $\mathcal{A}$ of DT invariants of $(Q,W)$ and the generating series $\bar{\mathcal{A}}$ of DT invariants of $(\bar{Q},\bar{W})$:
\begin{align}
    \mathcal{A}^{z:I}(x)=\bar{\mathcal{A}}^{\bar{z}:I}((x^{\alpha^z_k})_k)
\end{align}
Hence, using the fact that $d$ and $\bar{d}$ are in the kernel of the anti-symmetrized Euler forms, and then that the corresponding BPS invariants $\Omega$ of $(Q,W)$ and $\bar{\Omega}$ of $(\bar{Q},\bar{W})$ does not depends on the stability parameter:
\begin{align}
    \Omega^{z:I}_d&=\bar{\Omega}^{\bar{z}:I}_{\bar{d}}\nn\\
    \Rightarrow \Omega^{z:I}(x)&=\bar{\Omega}^{\bar{z}:I}((x^{\alpha^z_k})_k)
\end{align}

Then \cite[Theo 6.1]{MozmcKay} gives:
\begin{align}\label{mckay}
    \bar{\Omega}((x^{\alpha^z_k})_k)=(\mathbb{L}^{3/2}+(K_z-1)\mathbb{L}^{1/2})\sum_{n\geq 1}x^{n\delta}+\mathbb{L}^{1/2}\sum_{k\neq k'}\sum_{n\geq 0}x^{n\delta+\alpha^z_{[k,k'[}}
\end{align}
We could compute $\bar{\Omega}^{z:I}$ by doing a method similar as in \cite{MorNag}, considering invertible/nilpotent decompositions and Jordan block decompositions. We prefer to show to extract $\bar{\Omega}^{\bar{z}:I}$ from $\bar{\Omega}$ using only formal manipulations, as an illustration of our formalism. We will prove the claim below:
\begin{align}\label{claimccast}
    \bar{\Omega}^{\bar{z}:I}(x)=\frac{[\mathbb{C}^\ast]}{[\mathbb{C}]}\bar{\Omega}(x)
\end{align}
This is related with the fact that $\bar{\Omega}^{\bar{z}:I}(x)$ (resp. $\bar{\Omega}(x)$) is the generating series of BPS invariants for a noncommutative crepant resolution of $\mathbb{C}^2/\mathbb{Z}_{K_z}\times\mathbb{C}^{\ast}$(resp. of $\mathbb{C}^2/\mathbb{Z}_{K_z}\times\mathbb{C}$). Consider a generic stability $\theta$ giving a crepant resolution $(\mathbb{C}^2/\mathbb{Z}_{K_z}\times\mathbb{C})_\theta=(\mathbb{C}^2/\mathbb{Z}_{K_z})_\theta\times\mathbb{C}$, and remark that $(\mathbb{C}^2/\mathbb{Z}_{K_z}\times\mathbb{C})^{z:I}_\theta=(\mathbb{C}^2/\mathbb{Z}_{K_z})_\theta\times\mathbb{C}^\ast$, using \eqref{ndelta}, one obtains for $n\geq 1$:
\begin{align}
    &\bar{\Omega}^{\bar{z}:I}(n\delta,y)=\mathbb{L}^{-3/2}[(\mathbb{C}^2/\mathbb{Z}_{K_z})_\theta\times\mathbb{C}^\ast]\nn\\
    &\bar{\Omega}(n\delta,y)=\mathbb{L}^{-3/2}[(\mathbb{C}^2/\mathbb{Z}_{K_z})_\theta\times\mathbb{C}]\nn\\
    \Rightarrow &\bar{\Omega}^{\bar{z}:I}(n\delta,y)=\frac{[\mathbb{C}^\ast]}{[\mathbb{C}]}\bar{\Omega}(n\delta,y)
\end{align}
We show that the relation of the claim holds also for the other dimension vectors using invertible/nilpotent decompositions and duality properties for $(\bar{Q},\bar{W})$. Denote by $\bar{z},\bar{z}',\bar{z}''$ the external edges of the toric diagram of $\mathbb{C}^2/\mathbb{Z}_{K_z}\times\mathbb{C}$ considered in the clockwise order. One has:
\begin{align}
    &\bar{\Omega}^{\bar{z}:I}(x)=\bar{\Omega}(x)-\bar{\Omega}^{\bar{z}:N}(x)\nn\\
    &\bar{\Omega}^{\bar{z}:N}(x)=\mathbb{D}(\bar{\Omega}^{\bar{z}':N,\bar{z}'':N}(x))\nn\\
    &\bar{\Omega}^{\bar{z}':N,\bar{z}'':N}(x)=\bar{\Omega}(x)-\bar{\Omega}^{\bar{z}':I,\bar{z}'':N}(x)-\bar{\Omega}^{\bar{z}'':I}(x)\nn\\
    \Rightarrow & \bar{\Omega}^{\bar{z}:I}(x)=\bar{\Omega}(x)-\mathbb{D}(\bar{\Omega}(x))+\mathbb{D}(\bar{\Omega}^{\bar{z}':I,\bar{z}'':N}(x))+\mathbb{D}(\bar{\Omega}^{\bar{z}'':I}(x))
\end{align}
where $\mathbb{D}$ is the Verdier duality, \ie the Poincaré duality for mixed Hodge structures. In the first and the third lines we have performed invertible/nilpotent decompositions using Proposition \ref{nilinvdec}, and in the second line we have used Corollary \ref{lemdual}. We have $K_{\bar{z}'}=K_{\bar{z}''}=1$, \ie the BPS invariants $\bar{\Omega}^{\bar{z}':I,\bar{z}'':N}$ and $\bar{\Omega}^{\bar{z}'':I}$ have only terms with dimension vector $n\delta$, giving:
\begin{align}
    \bar{\Omega}^{\bar{z}:I}_d=&\bar{\Omega}_d-\Sigma(\bar{\Omega}_d)\quad\forall d\not\in\langle\delta\rangle\nn\\
    =&\frac{[\mathbb{C}^\ast]}{[\mathbb{C}]}\bar{\Omega}_d\quad\forall d\not\in\langle\delta\rangle
\end{align}
The last line holds because the BPS invariants are either $0$ or $\mathbb{L}^{1/2}$ for $d\not\in\langle\delta\rangle$ from \eqref{mckay}. This ends the proof the claim \eqref{claimccast}, giving:
\begin{align}
    \Omega^{z:I}(x)&=\bar{\Omega}^{\bar{z}:I}((x^{\alpha^z_k})_k)\nn\\&\resizebox{0.9\hsize}{!}{$=(\mathbb{L}^{3/2}+(K_z-2)\mathbb{L}^{1/2}-(K_z-1)\mathbb{L}^{-1/2})\sum_{n\geq 1}x^{n\delta}+(\mathbb{L}^{1/2}-\mathbb{L}^{-1/2})\sum_{k\neq k'}\sum_{n\geq 0}x^{n\delta+\alpha^z_{[k,k'[}}$}
\end{align}
$\Box$

\subsection{Identities between partially nilpotent attractors invariants}

We  now  have all necessary ingredients to express the BPS invariants $\Omega_\theta$ in terms of BPS invariants $\Omega_\theta^{Z_N:N}$ with nilpotency constraints on a given cycle $v^z$.\medskip

\begin{proposition}\label{nilpart}
We can express, for $[z,z']$ a strict subset of the set of sides of the toric diagram:
\begin{align}
    \Omega_\theta(x)=&(\mathbb{L}^{3/2}+({\textstyle\sum_{\tilde{z}\in[z,z']}} K_{\tilde{z}}-2)\mathbb{L}^{1/2}-({\textstyle\sum_{\tilde{z}\in[z,z']}} K_{\tilde{z}}-1)\mathbb{L}^{-1/2})\sum_{n\geq 1} x^{n\delta}\nonumber\\
    &+(\mathbb{L}^{1/2}-\mathbb{L}^{-1/2})\sum_{\tilde{z}\in[z,z']}\sum_{k\neq k'}\sum_{n\geq 0} x^{n\delta+\alpha^{\tilde{z}}_{[k,k'[}}+\Omega_\theta^{[z,z']:N}(x)\label{corromega}\\
    \Omega_\theta(x)=&(\mathbb{L}^{3/2}+(b-3)\mathbb{L}^{1/2}-(b-3)\mathbb{L}^{-1/2}-\mathbb{L}^{-3/2})\sum_{n\geq 1} x^{n\delta}\nonumber\\ &+(\mathbb{L}^{1/2}-\mathbb{L}^{-1/2})\sum_{z}\sum_{k\neq k'}\sum_{n\geq 0} x^{n\delta+\alpha^z_{[k,k'[}}+\Omega_\theta^{N}(x)
\end{align}
\end{proposition}

Proof:

We use Proposition \ref{2partinv}, which says that for 
non adjacent $z_i$ and $z$, $\Omega^{z_i:I,z:I}=\Omega^I$. In particular, for $z\not\in [z_{i-1},z_{i+1}]$ we have:
\begin{align}\label{nilchain}
    &\Omega^{z_i:I,z_{i+1}:N,z:I}=0 \nn\\
    \Rightarrow &\Omega^{z_i:I,[z_{i+1},z_{i-1}[:N}=\Omega_\theta^{z_i:I,z_{i+1}:N}
\end{align}
We have also:
\begin{align}\label{nilchain2}
    \Omega^{z_i:I,[z_{i+1},z_i[:N}=&\Omega^{z_i:I,[z_{i+1},z_{i-1}[:N}-\Omega^{z_{i-1}:I,z_i:I,[z_{i+1},z_{i-1}[:N}\nn \\
    =&\Omega^{z_i:I,z_{i+1}:N}-\Omega^{z_{i-1}:I,z_i:I}+\Omega^I
\end{align}
Graphically, the two equations \eqref{nilchain} and \eqref{nilchain2} can be written:

\centerline{
\begin{tikzpicture} [scale=0.5]
\begin{scope}[xshift=4cm]
\draw[->,>=latex] (0,0) -- (-1,0);
\draw[->,>=latex] (0,0) -- (0.5,1);
\draw[->,>=latex] (0,0) -- (-0.5,1);
\draw[->,>=latex] (0,0) -- (1,0);
\draw[->,>=latex] (0,0) -- (0.5,-1);
\draw[->,>=latex] (0,0) -- (-0.5,-1);
\node at (-0.5,1) [above]{$I$};
\node at (0.5,1) [above]{$N$};
\end{scope}
\node at (2.2,0) {$=$};

\begin{scope}
\draw[->,>=latex] (0,0) -- (-1,0);
\draw[->,>=latex] (0,0) -- (0.5,1);
\draw[->,>=latex] (0,0) -- (-0.5,1);
\draw[->,>=latex] (0,0) -- (1,0);
\draw[->,>=latex] (0,0) -- (0.5,-1);
\draw[->,>=latex] (0,0) -- (-0.5,-1);
\node at (-0.5,1) [above]{$I$};
\node at (0.5,1) [above]{$N$};
\node at (1,0) [right]{$N$};
\node at (0.5,-1) [below]{$N$};
\node at (-0.5,-1) [below]{$N$};
\end{scope}
\end{tikzpicture}}
\centerline{
\begin{tikzpicture}[scale=0.5] 
\begin{scope}
\draw[->,>=latex] (0,0) -- (-1,0);
\draw[->,>=latex] (0,0) -- (0.5,1);
\draw[->,>=latex] (0,0) -- (-0.5,1);
\draw[->,>=latex] (0,0) -- (1,0);
\draw[->,>=latex] (0,0) -- (0.5,-1);
\draw[->,>=latex] (0,0) -- (-0.5,-1);
\node at (-1,0) [left]{$N$};
\node at (0.5,1) [above]{$N$};
\node at (-0.5,1) [above]{$I$};
\node at (1,0) [right]{$N$};
\node at (0.5,-1) [below]{$N$};
\node at (-0.5,-1) [below]{$N$};
\end{scope}
\node at (2.2,0) {$=$};

\begin{scope}[xshift=4cm]
\draw[->,>=latex] (0,0) -- (-1,0);
\draw[->,>=latex] (0,0) -- (0.5,1);
\draw[->,>=latex] (0,0) -- (-0.5,1);
\draw[->,>=latex] (0,0) -- (1,0);
\draw[->,>=latex] (0,0) -- (0.5,-1);
\draw[->,>=latex] (0,0) -- (-0.5,-1);
\node at (-0.5,1) [above]{$I$};
\node at (0.5,1) [above]{$N$};
\end{scope}
\node at (6,0) {$-$};

\begin{scope}[xshift=8cm]
\draw[->,>=latex] (0,0) -- (-1,0);
\draw[->,>=latex] (0,0) -- (0.5,1);
\draw[->,>=latex] (0,0) -- (-0.5,1);
\draw[->,>=latex] (0,0) -- (1,0);
\draw[->,>=latex] (0,0) -- (0.5,-1);
\draw[->,>=latex] (0,0) -- (-0.5,-1);
\node at (-0.5,1) [above]{$I$};
\node at (-1,0) [above]{$I$};
\end{scope}
\node at (10,0) {$+$};

\begin{scope}[xshift=12cm]
\draw[->,>=latex] (0,0) -- (-1,0);
\draw[->,>=latex] (0,0) -- (0.5,1);
\draw[->,>=latex] (0,0) -- (-0.5,1);
\draw[->,>=latex] (0,0) -- (1,0);
\draw[->,>=latex] (0,0) -- (0.5,-1);
\draw[->,>=latex] (0,0) -- (-0.5,-1);
\node at (-1,0) [left]{$I$};
\node at (0.5,1) [above]{$I$};
\node at (-0.5,1) [above]{$I$};
\node at (1,0) [right]{$I$};
\node at (0.5,-1) [below]{$I$};
\node at (-0.5,-1) [below]{$I$};
\end{scope}
\end{tikzpicture}}
We can combine the formulas of Propositions \ref{2partinv} and \ref{1partinv}:
\begin{align}\label{partIN}
    \Omega^{z_l:I,z_{l+1}:N}=&\Omega^{z_l:I}-\Omega^{z_l:I,z_{l+1}:I}\nn\\
    =&(\mathbb{L}^{1/2}-\mathbb{L}^{-1/2})\left(K_z\sum_{n\geq 1}x^{n\delta}+\sum_{k\neq k'}\sum_{n\geq 0}x^{\alpha^z_{[k,k'[}+n\delta}\right)
\end{align}
We  decompose successively, denoting for convenience $z=z_i,z'=z_j$:
\begin{align}
    \Omega_\theta=&\Omega^{z_j:I}+\sum_{z_l\in[z_i,z_j[}\Omega^{z_l:I,]z_l,z_j]:N}+\Omega_\theta^{[z_i,z_j]:N} \nn\\
    =&\Omega^{z_j:I,z_{j+1}:I}+\Omega^{z_j:I,z_{j+1}:N}+\sum_{z_l\in[z_i,z_j[}\Omega^{z_l:I,z_{l+1}:N}+\Omega_\theta^{[z_i,z_j]:N}\nn\\
    =&(\mathbb{L}^{3/2}+({\textstyle\sum_{z\in[z_i,z_j]}} K_z-2)\mathbb{L}^{1/2}-({\textstyle\sum_{z\in[z_i,z_j]}} K_z-1)\mathbb{L}^{-1/2})\sum_{n\geq 1} x^{n\delta}\nonumber\\
    &+(\mathbb{L}^{1/2}-\mathbb{L}^{-1/2})\sum_{z\in[z_i,z_j]}\sum_{k\neq k'}\sum_{n\geq 0} x^{n\delta+\alpha^z_{[k,k'[}}+\Omega_\theta^{[z_i,z_j]:N}(x)
\end{align}
Here in the second line we have used \eqref{nilchain}, and we have used \eqref{partIN} and Proposition \ref{2partinv} in the last line. The manipulations above can be represented graphically as:

\centerline{
\begin{tikzpicture} [scale=0.4]
\begin{scope}
\draw[->,>=latex] (0,0) -- (-1,0);
\draw[->,>=latex] (0,0) -- (0.5,1);
\draw[->,>=latex] (0,0) -- (-0.5,1);
\draw[->,>=latex] (0,0) -- (1,0);
\draw[->,>=latex] (0,0) -- (0.5,-1);
\draw[->,>=latex] (0,0) -- (-0.5,-1);
\end{scope}
\node at (2.5,0) {$=$};

\begin{scope}[xshift=5cm]
\draw[->,>=latex] (0,0) -- (-1,0);
\draw[->,>=latex] (0,0) -- (0.5,1);
\draw[->,>=latex] (0,0) -- (-0.5,1);
\draw[->,>=latex] (0,0) -- (1,0);
\draw[->,>=latex] (0,0) -- (0.5,-1);
\draw[->,>=latex] (0,0) -- (-0.5,-1);
\node at (1,0) [right]{$I$};
\end{scope}
\node at (7.5,0) {$+$};

\begin{scope}[xshift=10cm]
\draw[->,>=latex] (0,0) -- (-1,0);
\draw[->,>=latex] (0,0) -- (0.5,1);
\draw[->,>=latex] (0,0) -- (-0.5,1);
\draw[->,>=latex] (0,0) -- (1,0);
\draw[->,>=latex] (0,0) -- (0.5,-1);
\draw[->,>=latex] (0,0) -- (-0.5,-1);
\node at (0.5,1) [above]{$I$};
\node at (1,0) [right]{$N$};
\end{scope}
\node at (12.5,0) {$+$};

\begin{scope}[xshift=15cm]
\draw[->,>=latex] (0,0) -- (-1,0);
\draw[->,>=latex] (0,0) -- (0.5,1);
\draw[->,>=latex] (0,0) -- (-0.5,1);
\draw[->,>=latex] (0,0) -- (1,0);
\draw[->,>=latex] (0,0) -- (0.5,-1);
\draw[->,>=latex] (0,0) -- (-0.5,-1);
\node at (-0.5,1) [above]{$I$};
\node at (0.5,1) [above]{$N$};
\node at (1,0) [right]{$N$};
\end{scope}
\node at (17.5,0) {$+$};

\begin{scope}[xshift=20cm]
\draw[->,>=latex] (0,0) -- (-1,0);
\draw[->,>=latex] (0,0) -- (0.5,1);
\draw[->,>=latex] (0,0) -- (-0.5,1);
\draw[->,>=latex] (0,0) -- (1,0);
\draw[->,>=latex] (0,0) -- (0.5,-1);
\draw[->,>=latex] (0,0) -- (-0.5,-1);
\node at (-1,0) [left]{$I$};
\node at (-0.5,1) [above]{$N$};
\node at (0.5,1) [above]{$N$};
\node at (1,0) [right]{$N$};
\end{scope}
\node at (22.5,0) {$+$};
\begin{scope}[xshift=25cm]
\draw[->,>=latex] (0,0) -- (-1,0);
\draw[->,>=latex] (0,0) -- (0.5,1);
\draw[->,>=latex] (0,0) -- (-0.5,1);
\draw[->,>=latex] (0,0) -- (1,0);
\draw[->,>=latex] (0,0) -- (0.5,-1);
\draw[->,>=latex] (0,0) -- (-0.5,-1);
\node at (-1,0) [left]{$N$};
\node at (-0.5,1) [above]{$N$};
\node at (0.5,1) [above]{$N$};
\node at (1,0) [right]{$N$};
\end{scope}
\node at (2.5,-4) {$=$};

\begin{scope}[xshift=5cm,yshift=-4cm]
\draw[->,>=latex] (0,0) -- (-1,0);
\draw[->,>=latex] (0,0) -- (0.5,1);
\draw[->,>=latex] (0,0) -- (-0.5,1);
\draw[->,>=latex] (0,0) -- (1,0);
\draw[->,>=latex] (0,0) -- (0.5,-1);
\draw[->,>=latex] (0,0) -- (-0.5,-1);
\node at (1,0) [right]{$I$};
\node at (0.5,-1) [below]{$I$};
\end{scope}
\node at (7.5,-4) {$+$};

\begin{scope}[xshift=10cm,yshift=-4cm]
\draw[->,>=latex] (0,0) -- (-1,0);
\draw[->,>=latex] (0,0) -- (0.5,1);
\draw[->,>=latex] (0,0) -- (-0.5,1);
\draw[->,>=latex] (0,0) -- (1,0);
\draw[->,>=latex] (0,0) -- (0.5,-1);
\draw[->,>=latex] (0,0) -- (-0.5,-1);
\node at (1,0) [right]{$I$};
\node at (0.5,-1) [below]{$N$};
\end{scope}

\node at (12.5,-4) {$+$};

\begin{scope}[xshift=15cm,yshift=-4cm]
\draw[->,>=latex] (0,0) -- (-1,0);
\draw[->,>=latex] (0,0) -- (0.5,1);
\draw[->,>=latex] (0,0) -- (-0.5,1);
\draw[->,>=latex] (0,0) -- (1,0);
\draw[->,>=latex] (0,0) -- (0.5,-1);
\draw[->,>=latex] (0,0) -- (-0.5,-1);
\node at (0.5,1) [above]{$I$};
\node at (1,0) [right]{$N$};
\end{scope}
\node at (17.5,-4) {$+$};

\begin{scope}[xshift=20cm,yshift=-4cm]
\draw[->,>=latex] (0,0) -- (-1,0);
\draw[->,>=latex] (0,0) -- (0.5,1);
\draw[->,>=latex] (0,0) -- (-0.5,1);
\draw[->,>=latex] (0,0) -- (1,0);
\draw[->,>=latex] (0,0) -- (0.5,-1);
\draw[->,>=latex] (0,0) -- (-0.5,-1);
\node at (-0.5,1) [above]{$I$};
\node at (0.5,1) [above]{$N$};
\end{scope}
\node at (22.5,-4) {$+$};

\begin{scope}[xshift=25cm,yshift=-4cm]
\draw[->,>=latex] (0,0) -- (-1,0);
\draw[->,>=latex] (0,0) -- (0.5,1);
\draw[->,>=latex] (0,0) -- (-0.5,1);
\draw[->,>=latex] (0,0) -- (1,0);
\draw[->,>=latex] (0,0) -- (0.5,-1);
\draw[->,>=latex] (0,0) -- (-0.5,-1);
\node at (-1,0) [left]{$I$};
\node at (-0.5,1) [above]{$N$};
\end{scope}
\node at (27.5,-4) {$+$};
\begin{scope}[xshift=30cm,yshift=-4cm]
\draw[->,>=latex] (0,0) -- (-1,0);
\draw[->,>=latex] (0,0) -- (0.5,1);
\draw[->,>=latex] (0,0) -- (-0.5,1);
\draw[->,>=latex] (0,0) -- (1,0);
\draw[->,>=latex] (0,0) -- (0.5,-1);
\draw[->,>=latex] (0,0) -- (-0.5,-1);
\node at (-1,0) [left]{$N$};
\node at (-0.5,1) [above]{$N$};
\node at (0.5,1) [above]{$N$};
\node at (1,0) [right]{$N$};
\end{scope}
\end{tikzpicture}}

Similarly we can decompose :
\begin{align}
    \Omega_\theta=&\Omega^{z_{n-1}:I}+\sum_{z_l\in]z_0,z_{n-1}[}\Omega^{z_l:I,]z_l,z_{n-1}]:N,}+\Omega^{z_0:I,]z_0,z_{n-1}]:N}+\Omega_\theta^{N}\nn\\
    =&\Omega^{z_{n-1}:I,z_0:I}+\Omega^{z_{n-1}:I,z_0:N}+\sum_{z_l\in]z_0,z_{n-1}[}\Omega^{z_l:I,z_{l+1}:N}+\Omega^{z_0:I,z_1:N}-\Omega^{z_{n-1}:I,z_0:I}+\Omega^I+\Omega_\theta^{N}\nn\\
    =&\sum_{z_l\in[z_0,z_{n-1}]}\Omega^{z_l:I,z_{l+1}:N}+\Omega^I+\Omega_\theta^{N}\nn\\
    =&(\mathbb{L}^{3/2}+(b-3)\mathbb{L}^{1/2}-(b-3)\mathbb{L}^{-1/2}-\mathbb{L}^{-3/2})\sum_{n\geq 1} x^{n\delta}\nonumber\\ &+(\mathbb{L}^{1/2}-\mathbb{L}^{-1/2})\sum_{z}\sum_{k\neq k'}\sum_{n\geq 0} x^{n\delta+\alpha^z_{[k,k'[}}+\Omega_\theta^{N}(x)
\end{align}
Here we have used \eqref{nilchain} and \eqref{nilchain2} in the second line, we have simplified in the third line, and \eqref{partIN} and the formulas of Propositions \ref{2partinv} in the last line, recalling $b=\sum_z K_z$. Graphically, the manipulations above corresponds to:

\centerline{
\begin{tikzpicture} [scale=0.5]
\begin{scope}
\draw[->,>=latex] (0,0) -- (-1,0);
\draw[->,>=latex] (0,0) -- (0,1);
\draw[->,>=latex] (0,0) -- (1,0);
\draw[->,>=latex] (0,0) -- (0,-1);
\end{scope}
\node at (2.5,0) {$=$};
\begin{scope}[xshift=5cm]
\draw[->,>=latex] (0,0) -- (-1,0);
\draw[->,>=latex] (0,0) -- (0,1);
\draw[->,>=latex] (0,0) -- (1,0);
\draw[->,>=latex] (0,0) -- (0,-1);
\node at (0,-1) [below]{$I$};
\end{scope}
\node at (7.5,0) {$+$};
\begin{scope}[xshift=10cm]
\draw[->,>=latex] (0,0) -- (-1,0);
\draw[->,>=latex] (0,0) -- (0,1);
\draw[->,>=latex] (0,0) -- (1,0);
\draw[->,>=latex] (0,0) -- (0,-1);
\node at (1,0) [right]{$I$};
\node at (0,-1) [below]{$N$};
\end{scope}
\node at (12.5,0) {$+$};
\begin{scope}[xshift=15cm]
\draw[->,>=latex] (0,0) -- (-1,0);
\draw[->,>=latex] (0,0) -- (0,1);
\draw[->,>=latex] (0,0) -- (1,0);
\draw[->,>=latex] (0,0) -- (0,-1);
\node at (0,1) [above]{$I$};
\node at (1,0) [right]{$N$};
\node at (0,-1) [below]{$N$};
\end{scope}
\node at (17.5,0) {$+$};
\begin{scope}[xshift=20cm]
\draw[->,>=latex] (0,0) -- (-1,0);
\draw[->,>=latex] (0,0) -- (0,1);
\draw[->,>=latex] (0,0) -- (1,0);
\draw[->,>=latex] (0,0) -- (0,-1);
\node at (-1,0) [left]{$I$};
\node at (0,1) [above]{$N$};
\node at (1,0) [right]{$N$};
\node at (0,-1) [below]{$N$};
\end{scope}
\node at (22.5,0) {$+$};
\begin{scope}[xshift=25cm]
\draw[->,>=latex] (0,0) -- (-1,0);
\draw[->,>=latex] (0,0) -- (0,1);
\draw[->,>=latex] (0,0) -- (1,0);
\draw[->,>=latex] (0,0) -- (0,-1);
\node at (-1,0) [left]{$N$};
\node at (0,1) [above]{$N$};
\node at (1,0) [right]{$N$};
\node at (0,-1) [below]{$N$};
\end{scope}

\node at (2.5,-4) {$=$};
\begin{scope}[xshift=5cm,yshift=-4cm]
\draw[->,>=latex] (0,0) -- (-1,0);
\draw[->,>=latex] (0,0) -- (0,1);
\draw[->,>=latex] (0,0) -- (1,0);
\draw[->,>=latex] (0,0) -- (0,-1);
\node at (-1,0) [left]{$I$};
\node at (0,-1) [below]{$I$};
\draw[gray!60] (-1,-1) -- (1,1);
\draw[gray!60] (-1,1) -- (1,-1);
\end{scope}
\node at (7.5,-4) {$+$};
\begin{scope}[xshift=10cm,yshift=-4cm]
\draw[->,>=latex] (0,0) -- (-1,0);
\draw[->,>=latex] (0,0) -- (0,1);
\draw[->,>=latex] (0,0) -- (1,0);
\draw[->,>=latex] (0,0) -- (0,-1);
\node at (-1,0) [left]{$N$};
\node at (0,-1) [below]{$I$};
\end{scope}
\node at (12.5,-4) {$+$};
\begin{scope}[xshift=15cm,yshift=-4cm]
\draw[->,>=latex] (0,0) -- (-1,0);
\draw[->,>=latex] (0,0) -- (0,1);
\draw[->,>=latex] (0,0) -- (1,0);
\draw[->,>=latex] (0,0) -- (0,-1);
\node at (1,0) [right]{$I$};
\node at (0,-1) [below]{$N$};
\end{scope}
\node at (17.5,-4) {$+$};
\begin{scope}[xshift=20cm,yshift=-4cm]
\draw[->,>=latex] (0,0) -- (-1,0);
\draw[->,>=latex] (0,0) -- (0,1);
\draw[->,>=latex] (0,0) -- (1,0);
\draw[->,>=latex] (0,0) -- (0,-1);
\node at (0,1) [above]{$I$};
\node at (1,0) [right]{$N$};
\end{scope}
\node at (2.5,-8) {$+$};
\begin{scope}[xshift=5cm,yshift=-8cm]
\draw[->,>=latex] (0,0) -- (-1,0);
\draw[->,>=latex] (0,0) -- (0,1);
\draw[->,>=latex] (0,0) -- (1,0);
\draw[->,>=latex] (0,0) -- (0,-1);
\node at (-1,0) [left]{$I$};
\node at (0,1) [above]{$N$};
\end{scope}
\node at (7.5,-8) {$-$};
\begin{scope}[xshift=10cm,yshift=-8cm]
\draw[->,>=latex] (0,0) -- (-1,0);
\draw[->,>=latex] (0,0) -- (0,1);
\draw[->,>=latex] (0,0) -- (1,0);
\draw[->,>=latex] (0,0) -- (0,-1);
\node at (-1,0) [left]{$I$};
\node at (0,-1) [below]{$I$};
\draw[gray!60] (-1,-1) -- (1,1);
\draw[gray!60] (-1,1) -- (1,-1);
\end{scope}
\node at (12.5,-8) {$+$};
\begin{scope}[xshift=15cm,yshift=-8cm]
\draw[->,>=latex] (0,0) -- (-1,0);
\draw[->,>=latex] (0,0) -- (0,1);
\draw[->,>=latex] (0,0) -- (1,0);
\draw[->,>=latex] (0,0) -- (0,-1);
\node at (-1,0) [left]{$I$};
\node at (0,1) [above]{$I$};
\node at (1,0) [right]{$I$};
\node at (0,-1) [below]{$I$};
\end{scope}
\node at (17.5,-8) {$+$};
\begin{scope}[xshift=20cm,yshift=-8cm]
\draw[->,>=latex] (0,0) -- (-1,0);
\draw[->,>=latex] (0,0) -- (0,1);
\draw[->,>=latex] (0,0) -- (1,0);
\draw[->,>=latex] (0,0) -- (0,-1);
\node at (-1,0) [left]{$N$};
\node at (0,1) [above]{$N$};
\node at (1,0) [right]{$N$};
\node at (0,-1) [below]{$N$};
\end{scope}
\end{tikzpicture}}
Here the two crossed terms cancel. Notice that the correction between the partially nilpotent invariants lie in the center of the quantum affine space, \ie they are not subject to wall crossing: in particular, the same relations hold for the BPS invariants $\Omega_\theta$ at any generic stability $\theta$, an then also for attractor invariants $\Omega_\ast$.
$\Box$

Using the duality result of Corollary \ref{lemdual}, we are able to derive a universal formula expressing  BPS invariants up to an unknown self-Poincaré dual contribution:\medskip

\begin{theorem}\label{theoasym}
\begin{align}
    \Omega_\theta(x)=(\mathbb{L}^{3/2}+(b-3+i)\mathbb{L}^{1/2}+i\mathbb{L}^{-1/2})\sum_{n\geq 1} x^{n\delta}+\mathbb{L}^{1/2}\sum_{z}\sum_{k\neq k'}\sum_{n\geq 0} x^{n\delta+\alpha^z_{[k,k'[}}+\Omega^{sym}_\theta(x)
\end{align}
with $\Omega^{sym}_\theta(x)$ self Poincaré dual, and supported on dimension vectors $d\not\in\langle\delta\rangle$. The same formula holds for attractor invariants.
\end{theorem}

Proof: We have from Corollary \ref{lemdual}:
\begin{align}
    &\Omega_\theta(x)=\mathbb{D}(\Omega_\theta^N(x))
\end{align}
We have then, using the formula of Proposition \ref{nilpart}:
\begin{align}
    &\Omega_\theta(x)-\mathbb{D}(\Omega_\theta(x))=(\mathbb{L}^{3/2}+(b-3)\mathbb{L}^{1/2}-(b-3)\mathbb{L}^{-1/2}-\mathbb{L}^{-3/2})\sum_{n\geq 1} x^{n\delta}\nn\\
    &+(\mathbb{L}^{1/2}-\mathbb{L}^{-1/2})\sum_{z}\sum_{k\neq k'}\sum_{n\geq 0} x^{n\delta+\alpha^z_{[k,k'[}})
\end{align}
Hence, introducing:
\begin{align}
    &\resizebox{0.95\hsize}{!}{$\Omega^{sym}_\theta(x):=\Omega_\theta(x)-\left((\mathbb{L}^{3/2}+(b-3+i)\mathbb{L}^{1/2}+i\mathbb{L}^{-1/2})\sum_{n\geq 1} x^{n\delta}+\mathbb{L}^{1/2}\sum_{z}\sum_{k\neq k'}\sum_{n\geq 0} x^{n\delta+\alpha^z_{[k,k'[}}\right)$} \nn \\
    \Rightarrow &\Omega^{sym}_\theta(x)=\mathbb{D}(\Omega^{sym}_\theta(x))
\end{align}
From $i)$ of Proposition \ref{2partinv}, one further obtains that $\Omega^{sym}_\theta$ is supported on dimension vectors $d\not\in\langle\delta\rangle$. The same property follows for $\Omega_\ast(x)$, by noticing that $\Omega_{\ast,d}=\Omega_{\theta_d,d}$ for $\theta_d$ a specific stability parameter.
$\Box$

For local curves, \ie symmetric quivers corresponding to toric diagrams without interior lattice points, the BPS invariants (which does not depend on the stability $\theta$, because the quantum affine space is commutative in this case) have been computed explicitly. We will check as an illustration the compatibility of those results with our formula in Section \ref{localcurves}. It appears that the only contribution to the symmetric part $\Omega^{sym}(x)$ come from dimension vectors $d$ with $\Omega_d=1$.
\medskip

For toric diagrams with $i\geq1$ interior lattice points, the symmetric part $\Omega^{sym}_\theta(x)$ can be quite complicated, and in particular it is subject to wall crossing. The attractor invariants are expected to be simpler than BPS invariants for generic $\theta$. The simple representations, with dimension vectors $e_i,i\in Q_0$, always contribute to the attractor invariants, with $\Omega_{\ast,e_i}=1$ because there are no self 1-cycles in this case. A natural question is then whether there exist other dimension vectors for which the attractor invariants have a non-zero symmetric part $\Omega^{sym}_{\ast,d}$. We conjecture, based on evidence collected in \cite{BMP20} and on computations in \cite{MP20} recalled in Section \ref{seccompdiv}, that such dimension vectors do not exist:\medskip

\begin{conjecture}\label{conj}
For toric diagram with $i\geq 1$ internal lattice points, the attractor invariants are given by:
   \begin{align}
       \Omega_\ast(x)=\sum_i x_i+(\mathbb{L}^{3/2}+(b-3+i)\mathbb{L}^{1/2}+i\mathbb{L}^{-1/2})\sum_{n\geq 1} x^{n\delta}+\mathbb{L}^{1/2}\sum_{z}\sum_{k\neq k'}\sum_{n\geq 0} x^{n\delta+\alpha^z_{[k,k'[}}
    \end{align}
\end{conjecture}

\section{Toric localization for framed quivers with potential}

\subsection{Torus fixed variety and attracting variety}

Consider a one dimensional torus $\mathbb{C}^\ast$ acting on a variety $X$. We consider as in \cite{Braden2002HyperbolicLO}, \cite{drinfeld2013algebraic} and \cite{Ric16} the hyperbolic localization diagram:
\[\begin{tikzcd}
 X & X^\pm\arrow[l,swap,"\eta^\pm"]\arrow[r,"p^\pm"] & X^0
\end{tikzcd}\]
Here $X^0$ denotes the closed sub-variety of $\mathbb{C}^\ast$-fixed points, and $X^\pm$ the attracting (resp. repelling) variety, \ie the disjoint union of the components of $X$ flowing to a $\mathbb{C}^\ast$-fixed component when $t\to 0$ (resp. $t\to\infty$), $\eta^\pm$ gives the disjoint union of the closed embeddings of those components, and $p^\pm$ gives the projection to the $\mathbb{C}^\ast$-fixed component. The functors of constructible complexes $(p^+)_!(\eta^+)^\ast:D^b_c(X)\to D^b_c(X^0)$ is called the hyperbolic localization functor.\medskip

We consider the torus $(\mathbb{C}^\ast)^{(Q_f)_1}$ acting on $\mathbb{C}Q_f$, and therefore also on $\mathfrak{P}_f$, by scaling the arrows of $Q_f$.  We consider the subtorus $T$ leaving invariant the relations of the potential $W_f$, hence such that $W_f$ is homogeneous, with weight denoted by $\kappa$: its action on $\mathcal{M}_{Q_f,d}$ restricts to an action on $\mathcal{M}_{Q_f,W_f,d}$. The gauge torus $T_G=(\mathbb{C}^\ast)^{Q_0}$ acts on $(\mathbb{C}^\ast)^{(Q_f)_1}$ by adjunction $(t_a)_{(a:i\to j)\in (Q_f)_1}\mapsto(t_it_at_j^{-1})_{(a:i\to j)\in (Q_f)_1}$, where we denote $t_\infty=1$. The action of $T$ on $\mathcal{M}_{Q_f,d}$ and $\mathcal{M}_{Q_f,W_f,d}$ descends to an action of $T/T_G$. The scaling of the framing arrow $p:\infty\to i$ $p$ can be cancelled by the action of $T_G$, and for the $D4$ brane framing, the condition that the cycle $paq$ has weight $\kappa$ determines the weight of the relation arrow $q:j\to\infty$. Hence the torus acting on $\mathcal{M}_{Q_f,d}$ by leaving the potential homogeneous is the three dimensional torus $\mathbb{T}^3$ with weight lattice $M$, and the subtorus leaving the potential invariant is the two dimensional subtorus $\mathbb{T}^2$ with weight lattice $M/\kappa\mathbb{Z}=L$.\medskip

To study toric localization we will the consider one dimensional subtorus $\mathbb{C}^\ast\subset\mathbb{T}^2$ leaving the potential invariant. This data is called a choice of slope in K-theoretic Donaldson-Thomas theory, because such a torus is determined by the data of line passing through the origin of $L$, separating this lattice into two half planes $L^{>0}$ (resp $L^{<0}$), the half space of positive (resp negative) weights. In the following lemma, we consider the more general case of a subtorus $\mathbb{C}^\ast\subset\mathbb{T}^3$, which is then determined by the separation of $M$ into the half spaces $M^{>0}$ (resp $M^{<0}$) of positive (resp negative) weights:\medskip

\begin{lemma}\label{lemattr}
    $i)$ For $D$ the divisor corresponding to the corner $p$ of the toric diagram lying between the two sides $z,z'$, if $\lambda_z\in M^{>0}$ and $\lambda_{z'}\in M^{>0}$: 
    \begin{align}
        \mathcal{M}^+_{Q_D,W_D,d}=\mathcal{M}_{Q_D,W_D,d}
    \end{align}
    $ii)$ If $\lambda_{\tilde{z}}\in M^{<0}\iff\tilde{z}\in[z,z']$, then:
    \begin{align}
        \mathcal{M}^{\theta,ss,+}_{Q,W,d}&=\mathcal{M}^{\theta,ss,[z,z'[:N}_{Q,W,d}\nn\\
        \mathcal{M}^+_{Q_i,W,d}&=\mathcal{M}^{[z,z']:N}_{Q_i,W,d}
    \end{align}\medskip
\end{lemma}

Proof: Consider the set $C$ of cycles of $Q$ (resp $Q_f$) of length less than $\sum_{i\in Q_0}d_i$, and consider the map:
\begin{align}
    Tr:&\mathcal{M}^{\theta,ss}_{Q,d}\to\mathbb{C}^C\nn\\
    &V\hookrightarrow(Tr(w))_{w\in C}
\end{align}
and the same map for $\mathcal{M}_{Q_f,d}$. From general geometric invariant theory (see \cite{king_moduli}), these maps are projective, hence:
\begin{align}
    \mathcal{M}^{\theta,ss,\pm}_{Q,d}&=Tr^{-1}(\mathbb{C}^C)^\pm\nn\\
    \mathcal{M}^\pm_{Q_f,d}&=Tr^{-1}(\mathbb{C}^C)^\pm
\end{align}
\ie the attracting (resp repelling) variety is the variety of representations where the cycles in $M^{<0}$ are nilpotents.\medskip

From \cite[Corrolary 3.6]{MozI}, $M^+$, the lattice of weights of cycles of $Q$, is saturated in the cone $M^+_\mathbb{Q}$ which is the convex hull of the rays with direction $\lambda_{\tilde{z}}$, for $\tilde{z}$ a side of the toric diagram. For $w\in M^+$, there are then sides of the toric diagram $z_i$ and integers $n,n_i\in\mathbb{N}^\ast$ such that $w^n=\prod_i (v^z)^{n_i}$. In particular, $w\in M^{<0}$ if and only if there is a $z_i$ such that $\lambda_{z_i}\in M^{<0}$, and $w$ is nilpotent in a critical representation if and only if there is a $z_i$ which is nilpotent. Then in a critical representation the property "$w\in M^{<0}\iff w$ is nilpotent" is equivalent to "$\lambda_z\in M^{<0}\iff v^z$ is nilpotent". Hence if $l_{\tilde{z}}\in M^{<0}\iff\tilde{z}\in[z,z']$, one has $ii)$:
\begin{align}
        \mathcal{M}^{\theta,ss,+}_{Q,W,d}&=\mathcal{M}^{\theta,ss,[z,z']:N}_{Q,W,d}\nn\\
        \mathcal{M}^+_{Q_i,W,d}&=\mathcal{M}^{[z,z']:N}_{Q_i,W,d}
    \end{align}
Because each arrow of $I$ acts trivially on $P_D$, only the sub-cone $M^+_\mathbb{Q}\cap\bar{\chi_I}^{-1}(0)$ acts non-trivially in a framed representation of $(Q_D,W_D)$, hence $v^{z''}$ acts trivially for $z''\neq z,z'$, from which we deduce $i)$.$\Box$\medskip

\begin{corollary}\label{lemdual}
For any dimension vector $d$ and stability parameter $\theta$, one has:
\begin{align}
    \Omega^{[z,z'[:N}_{\theta,d}=\mathbb{D}(\Omega^{[z',z[:N}_{\theta,d})
\end{align}
Here $\mathbb{D}$ denotes the Poincaré duality at the level of monodromic mixed Hodge structures.
\end{corollary}

Proof: Consider a $\mathbb{C}^\ast$ action scaling all the arrows of the quiver, such that $\lambda_{\tilde{z}}\in M^{<0}\iff \tilde{z}\in [z,z'[$ (such an action exists because $M^+$ is convex). Hence from Corollary \ref{lemdual}:
\begin{align}
    \mathcal{M}^{\theta,ss,+}_{Q,d}&=\mathcal{M}^{\theta,ss,[z,z'[:N}_{Q,d}\nn\\
    \mathcal{M}^{\theta,ss,-}_{Q,d}&=\mathcal{M}^{\theta,ss,[z',z[:N}_{Q,d}
\end{align}
Then denoting by $q:\mathcal{M}^{\theta,ss,0}_{Q,d}\to \ast$ the projection to a point, one obtains:
\begin{align}
    \mathbb{D}(\Omega^{[z',z[:N}_{\theta,d})&=\mathbb{D}H_c(\mathcal{M}^{\theta,ss,-}_{Q,d},\phi_W\mathcal{IC}_{\mathcal{M}^{\theta,ss}_{Q,d}})\nn\\
    &=\mathbb{D}q_!(p^-)_!(\eta^-)^\ast\phi_W\mathcal{IC}_{\mathcal{M}^{\theta,ss}_{Q,d}}\nn\\
    &=q_!(p^+)_!(\eta^+)^\ast\mathbb{D}\phi_W\mathcal{IC}_{\mathcal{M}^{\theta,ss}_{Q,d}}\nn\\
    &=q_!(p^+)_!(\eta^+)^\ast\phi_W\mathcal{IC}_{\mathcal{M}^{\theta,ss}_{Q,d}}\nn\\
    &=H_c(\mathcal{M}^{\theta,ss,+}_{Q,d},\phi_W\mathcal{IC}_{\mathcal{M}^{\theta,ss}_{Q,d}})\nn\\
    &=\Omega^{[z,z'[:N}_{\theta,d}
\end{align}
Here we have used the definitions of the BPS invariants and the characterization of the attracting and repelling varieties in the first and last lines, the fact that $q_!=q_\ast$ because $q$ is projective, and Braden's contraction lemma of \cite{Braden2002HyperbolicLO} in its form \cite[Theo B]{Ric16} using the fact that $\phi_W\mathcal{IC}_{\mathcal{M}^{\theta,ss}_{Q,d}}$ is $\mathbb{C}^\ast$-equivariant in the third line, and the self-duality of the vanishing cycles functor and the intersection complex in the fourth line.$\Box$\medskip

We will now describe the fixed points of the action of $\mathbb{T}^2$, the torus leaving the potential invariant, on $\mathcal{M}_{Q_f,W_f,d}$. The action of the torus $\mathbb{T}^3$ scaling the arrows of the quiver by leaving the potential invariant induces a $\Lambda$-grading on $J_{Q_f,W_f}$, and $P_f$ has a $\Lambda$-grading as an $J_{Q_f,W_f}$-module, \ie 
\begin{align}
    a.(P_f)_\lambda\subset (P_f)_{\lambda+\wt(a)}\ ,\quad \mbox{for}\;\lambda\in\Lambda
\end{align}
A path which does not vanish in $P_i$ (resp.  $P_D$) is of the form $vp$, with $v$ a path of $Q$, and two paths with he same $\Lambda$-weights agree in $J_{Q,W}$, then $(P_i)_\lambda$ (resp. $(P_D)_\lambda$) is at most one dimensional for $\lambda\in\Lambda$. We define the Empty Room (ERC) Configuration $\Delta_f$ as the subset of $\Lambda$ such that $(P_f)_\lambda$ is not empty, hence one dimensional. One calls the elements of $\lambda\in\Delta_f$ such that $d_1(\lambda)=i$ the atoms of color $i$ of the ERC. One denotes $\lambda\leq \mu$ for $\lambda,\mu\in\Lambda$ if there exist $v\in \mathbb{C}Q_f$ such that $\mu=\lambda+\wt(v)$. The relation $\leq$ is manifestly reflexive and transitive. If $\lambda\leq\mu\leq\lambda$, then  there are paths $v,w\in\mathbb{C}Q_f$ such that $\wt(v)+\wt(w)=0$, \ie $wv$ has $\Lambda$ weight $0$, then from \cite[Prop 4.8]{MR} $wv$ is trivial in $\mathbb{C}Q$, and then $w:j\to i$ is trivial in $\mathbb{C}Q_f$, giving $\lambda=\mu$: $\leq$ is then anti-symmetric. Thus $\leq$ defines a poset structure on $\Delta_f$.\medskip

We denote by $\Pi_f$ the set of finite ideals of $\Delta_f$, \ie subsets $\pi\subset \Delta_f$ such that $x\leq y,y\in \pi\Rightarrow x\in\pi$. For $\pi\in\Pi_f$ we denote by $d_\pi\in\mathbb{N}^{Q_0}$ the dimension vector such that $(d_\pi)_i$ gives the number of atoms of color $i$ in $\pi$. Those posets can be visualized in three dimension by choosing a quasi-inverse $\Lambda\to M$ of the embedding $M\hookrightarrow \Lambda$. In the D6 brane framing case, $\Delta_i$ can be seen as a three-dimensional pyramid with atom $i$ on the top which is the dual of the toric fan of $X$, with each atom being obtained by a path of the quiver. In the D4-brane framing case, the action of the relation arrow $q:j\to\infty$ and of all the arrows of $I$ are trivial on $P_D$: consider the quiver with relations $(Q_I,\partial_IW)$, where one has removed the arrows of $I$ and enforced the relations $\partial_a W$ for $a\in I$. Hence $\Delta_D$ can be seen as a facet of the pyramid $\Delta_i$, which is the dual of the ray of the toric fan supporting the corner of the toric diagram corresponding to $d$, obtained by considering only paths of $(Q_I,\partial_IW)$ starting at $i$. \medskip

\begin{lemma}\label{lemfix}
    The $\mathbb{T}^2$-fixed variety $\mathcal{M}^{\mathbb{T}^2}_{Q_f,W_f,d}$ is an union of isolated points, which are in natural bijection with the set of $d$-dimensional pyramids in $\Pi_f$.
\end{lemma}

Proof: As in \cite[Theo 2.4]{MR}, f-cyclic representations fixed by a torus actions are quotients of $P_f$ by an ideal $\rho$ which i homogeneous under this torus action, hence the elements of $\mathcal{M}^{\mathbb{T}^2}_{Q_f,W_f,d}$ are the quotients $P_f/\rho$ for $\rho$ $\Lambda_f/\mathbb{Z}\kappa$-homogeneous.

It was shown in \cite[Rem 4.10]{MR} that any path $v:i\to \iota$ of $Q$ can be written in $\mathbb{C}Q/\partial W$ as $w^nv_0$, for $v_0:i\to \iota$ a minimal path of the unframed quiver, $w:\iota\to \iota$ an arbitrary cycle of the potential $W$, and $n\in\mathbb{N}$. Recall that any cycle of $W$ contains an arrow of $I$, and then has a trivial action on $P_D$. Then for $\iota\in Q_0$, $l\in\Lambda/\mathbb{Z}\kappa$ (resp. $l\in\Lambda/\mathbb{Z}\kappa$) such that $(P_i)_l$ (resp. $(P_D)_l$) is not empty, one has:
\begin{align}
    &(P_D)_l=\langle v_0p\rangle\nn\\
    &(P_i)_l=\langle (w^nv_0p)_{n\in\mathbb{N}}\rangle
\end{align}
with $v_0:i\to \iota$ a minimal path of the periodic quiver and $w:\iota\to \iota$ a cycle of $W$. Consider $\rho=\bigoplus_l\rho_l$ a $\Lambda/\mathbb{Z}\kappa$ homogeneous sub-module of $P_i$ (resp. of $P_D$) with finite codimension. For a 'D4 brane' framing $\rho_l$ is automatically $\mathbb{Z}\kappa$-homogeneous, and then $\rho=\bigoplus_l\rho_l$ is $\Lambda$-homogeneous. For $z$ a side of the toric diagram, $Tr((v^z)^n)$ is scaled by $\mathbb{T}^2$, hence vanish on a $\mathbb{T}^2$-fixed point, \ie $v^z$ is nilpotent on a $\mathbb{T}^2$-fixed point. Because $M^+$ is saturated in $M_\mathbb{Q}$, there is an $n\in\mathbb{N}^\ast$ such that $w^n$ can be expressed as a product of $v^z$, then $w$ is itself nilpotent on a $\mathbb{T}^2$-fixed point. In particular $\rho_l=\langle (w^nv_0p)_{n\geq N}\rangle$ for $N\in\mathbb{N}$, hence it is $\kappa\mathbb{Z}$-homogeneous, therefore $\rho$,is $\Lambda$-homogeneous. Hence:
\begin{align}
    \mathcal{M}^{\mathbb{T}^2}_{Q_f,W_f,d}=\mathcal{M}^{\mathbb{T}^3}_{Q_f,W_f,d}
\end{align}
and $\mathbb{T}^2$-fixed points corresponds to quotients $P_f/\rho$ with $\rho$ $\Lambda$-homogeneous. A $\Lambda$-homogeneous sub-module $\rho$  with codimension $(1,d)$ of $P_f$ is then a sum of graded components of $P_f$ (recall that $(P_f)_\lambda$ is at most one dimensional):
\begin{align}
    \rho=\bigoplus_{\lambda\not\in\pi}(P_f)_\lambda
\end{align}
The condition that $\rho$ is a sub-module is equivalent, by construction of the poset structure on $\Delta_f$, to the condition that $\pi$ is an ideal of $\Delta_f$, and we have $d_\pi=d$. There is then a natural bijection between the set of $\Lambda_f$-homogeneous sub-module of $P_f$ with finite codimension (which by $iii)$ in the assumption is equal to the set of  $\Lambda/\mathbb{Z}\kappa$-homogeneous sub-module of $P_f$ with finite codimension) and $\Pi_f$. $\Box$\medskip

\subsection{The tangent-obstruction complex}

We have seen that representations in $\mathcal{M}^{\mathbb{T}^3}_{Q_f,d}$ are $\Lambda$-graded. For $i\in Q_0$, there is a tautological sheaf $V_i$ on the moduli space of representations, whose stalk at each points corresponding to a representation $V$ is the vector space $V_i$ at the node $i$. The restriction of these tautological sheaves on $\mathcal{M}^{\mathbb{T}^3}_{Q_f,d}$ are then $\Lambda$-graded. The restriction $T_{\mathcal{M}_{Q_f,d}}|_{\mathcal{M}^{\mathbb{T}^3}_{Q_f,d}}$ of the tangent space of $\mathcal{M}_{Q_f,d}$ at $\mathcal{M}^{\mathbb{T}^3}_{Q_f,d}$ as a $\mathbb{T}^3$-equivariant structure. Denoting by $t^a$ the $\mathbb{T}^3$-equivariant line bundle with weight $a$, the $\mathbb{T}^3$-equivariant tangent space on $\mathcal{M}^{\mathbb{T}^3}_{Q_f,\pi}$ is the cokernel of the map of fiber bundles:
\begin{align}
    S^0_\pi\overset{\delta_0}{\to}S^1_\pi\nn\\
\end{align}
Where:
\begin{itemize}
\item The fiber bundle $S^0_\pi$ is the space of infinitesimal gauge transformations $\delta g_i$ (we denote for convenience $\delta g_i=0$ for i a framing node):
\begin{align}
    S^0_\pi=\bigoplus_{i\in Q_0} \Hom_\mathbb{C}(V_{i},V_{i})
    \label{S0}
\end{align}
\item $S^1_\pi$ is the space of infinitesimal deformations of the arrows $(\delta a)$:
\begin{align}
S^1_\pi=\bigoplus_{(a:i\to j)\in (Q_f)_1} \Hom_\mathbb{C}(V_{i},V_{j})\otimes t^a\label{S1}
\end{align}
\item The differential $\delta_0$ is the linearization of gauge transformations (taking care of the fact that framing nodes are not gauged):
\begin{align}
\delta_0:(\delta g_i)_{i\in Q_0}\mapsto(\delta g_j a-a\delta g_i)_{(a:i\to j)\in (Q_f)_1}
\end{align}
\end{itemize}
Notice that this complex is in fact $M$-graded.\medskip

For a torus $\mathbb{C}^\ast\subset\mathbb{T}^2$ leaving the potential invariant, the attracting and repelling behaviour of the action near the fixed component can be studied using this complex. The signed number of weights in $L^{>0}$ (resp in $L^0$, resp in $L^{<0}$) in this complex gives the number $d^+_\pi$ (resp $d^0_\pi$, resp $d^-_\pi$) of contracting (resp invariant, resp repelling) weights in $T_{\mathcal{M}_{Q_f,d}}|_{\mathcal{M}^{\mathbb{T}^3}_{Q_f,d}}$. We define then, following the notations used in K-theoretic Donaldson-Thomas theory:
\begin{align}
    \Index^s_\pi:=d^+_\pi-d^-_\pi
\end{align}
here we have insisted on the dependency on the slope $s$. Notice that the $\mathbb{T}^3$-equivariant structure of $S^0_\pi$ is self-dual, hence the number of contracting and repelling weights in $S^0_\pi$ are equal: to compute $\Index^s_\pi$, it suffice then to compute the difference between the number of contracting and repelling weights in $S^1_\pi$.\medskip

\begin{remark}
The weights of the cycles given by the subtorus $\mathbb{C}^\ast\subset\mathbb{T}^2$ must be integers, hence the line separating $L^0$ must be directed by an element of the lattice $L$. The general procedure to obtain a localization result is, for a given dimension vector, to choose a slope generic for this dimension vector, \ie such that the only weight of $L^0$ appearing in the tangent complex is $0$, and then such that the fixed points of $\mathbb{C}^\ast$ are the fixed points of $\mathbb{T}^2$. Then rigorously to compute the generating series of framed invariants one must consider a family of slope $(s_d)_d$, one for each dimension vector. To avoid heavy notations, we consider slopes $s$ with irrational coefficients, hence such that $L^0=\{0\}$, which we call generic slopes. For a given dimension vector $d$, we establish then the localization result by choosing a rational slope $s_d$ generic for $d$ and approximating $s$, \ie such that all the weights appearing in the tangent space, and also the weights $l_z$ of the cycles $v^z$, have the same contracting and repelling behaviour under the slopes $s$ and $s_d$.
\end{remark}

$\mathcal{M}_{Q_f,W_f,d}$ is the critical locus of the potential $Tr(W_f)$ inside the smooth scheme $\mathcal{M}_{Q_f,d}$. The derivative of $Tr(W_f)$ gives then a duality between the tangent directions and the obstructions of $\mathcal{M}_{Q_f,W_f,d}$, it is then a  $[-1]$-shifted symplectic scheme in the language of derived geometry. The Hessian of $Tr(W_f)$ defines then the tangent-obstruction complex:
\begin{align}
    0\to T_{\mathcal{M}_{Q_f,d}}\to T^\ast_{\mathcal{M}_{Q_f,d}}\to 0
\end{align}
One of the main idea of derived geometry is to replace the tangent space, which behaves not very well for singular spaces, by the tangent-obstruction complex, which behaves here far better. The tangent-obstruction complex of $\mathcal{M}_{Q_f,W_f,d}$ restricted to $\mathcal{M}^0_{Q_f,W_f,d}$ has also a $\mathbb{T}^3$-equivariant structure, and the obstructions spaces and tangent spaces are dual as $\mathbb{T}^2$-equivariant fiber bundles (but not as $\mathbb{T}^3$-equivariant fiber bundles). Hence $d^-_\pi$, the number of repelling weights in the tangent space is also the number of contracting weights in the obstruction space. Then $\Index^s_\pi$ is the signed number of contracting weights in the tangent obstruction complex.

\subsection{Derived Białynicki-Birula decomposition}

\begin{theorem}\label{theolocD}
    $i)$ For $D$ a non-compact divisor of $X$, corresponding to the corner $p$ of the toric diagram lying between the two sides $z,z'$, and a generic slope $s$ such that $l_z,l_{z'}\in L^{>0}$ (such slopes always exist, because the angle between $l_z$ and $l_{z'}$ is smaller than $\pi$), we have:
    \begin{align}
        Z_D(x)=\sum_{\pi\in\Pi_D}\mathbb{L}^{\Index^s_\pi/2}x^{d_\pi}
    \end{align}\\
    $ii)$ For a generic slope $s$ such that $l_{\tilde{z}}\in L^{<0}\iff\tilde{z}\in[z,z']$, one has:
    \begin{align}
    Z_i(x)&=S_{-i}[\Exp\left(\textstyle\sum_d\Delta^{s}\Omega_d\frac{\mathbb{L}^{d_i}-1}{\mathbb{L}^{1/2}-\mathbb{L}^{-1/2}}x^d\right)]\sum_{\pi\in\Pi_i}\mathbb{L}^{\Index^s_\pi/2}x^{d_\pi}
    \end{align}
    Using the correction term:
    \begin{align}
        \Delta^s\Omega(x)=&(\mathbb{L}^{3/2}+({\textstyle\sum_{\tilde{z}\in[z,z']}} K_{\tilde{z}}-2)\mathbb{L}^{1/2}-({\textstyle\sum_{\tilde{z}\in[z,z']}} K_{\tilde{z}}-1)\mathbb{L}^{-1/2})\sum_{n\geq 1} x^{n\delta}\nonumber\\
    &+(\mathbb{L}^{1/2}-\mathbb{L}^{-1/2})\sum_{\tilde{z}\in[z,z']}\sum_{k\neq k'}\sum_{n\geq 0} x^{n\delta+\alpha^{\tilde{z}}_{[k,k'[}}
    \end{align}
\end{theorem}

Proof: Consider the correspondence coming from the $\mathbb{C}^\ast$ action on the smooth scheme $\mathcal{M}_{Q_f,d}$:
\[\begin{tikzcd}
    \mathcal{M}_{Q_f,d}& \mathcal{M}^+_{Q_f,d}\arrow[l,swap,"p^+"]\arrow[r,"\eta^+"] &\mathcal{M}^0_{Q_f,d}
\end{tikzcd}\]
Consider the connected component $\mathcal{M}^0_{Q_f,\pi}$ of $\mathcal{M}^0_{Q_f,d}$ containing the $\mathbb{C}^\ast$-fixed representation corresponding to the pyramid $\pi\in\Pi_f$. The tangent space of $\mathcal{M}_{Q_f,d}$ at $\mathcal{M}^0_{Q_f,\pi}$ has $d^+_\pi$ contracting weights, $d^0_\pi$ invariant weights and $d^-_\pi$ repelling weights, with:
\begin{align}
    \dim(\mathcal{M}_{Q_f,d})&=d^+_\pi+d^0_\pi+d^-_\pi\nn\\
    \dim(\mathcal{M}^0_{Q_f,\pi})&=d^0_\pi\nn\\
    \Index^s_\pi&=d^+_\pi-d^-_\pi
\end{align}
then from Białynicki-Birula \cite{BiaynickiBirula1973SomeTO}, $\mathcal{M}^+_{Q_f,d}$ is a disjoint union of affine fiber bundles $p^+_\pi:\mathcal{M}^+_{Q_f,\pi}\to \mathcal{M}^0_{Q_f,\pi}$ of dimension $d^+_\pi$, and the fixed components are smooth of dimension $d^0_\pi$, hence using the hyperbolic localization functor:
\begin{align}\label{BB}
    (p^+)_!(\eta^+)^\ast\mathcal{IC}_{\mathcal{M}_{Q_f,d}}&=\bigoplus_{\pi\in\Pi_f|d_\pi=d}(p^+_\pi)_!\mathbb{Q}_{\mathcal{M}^+_{Q_f,\pi}}[d^+_\pi+d^0_\pi+d^-_\pi]\nn\\
    &=\bigoplus_{\pi\in\Pi_f|d_\pi=d}\mathbb{Q}_{\mathcal{M}^0_{Q_f,d}}[-d^+_\pi+d^0_\pi+d^-_\pi]\nn\\
    &=\bigoplus_{\pi\in\Pi_f|d_\pi=d}\mathbb{L}^{\Index^s_\pi/2}\mathcal{IC}_{\mathcal{M}^0_{Q_f,\pi}}
\end{align}
Here we have used the fact that $\mathcal{M}_{Q_f,d}$ and $\mathcal{M}^0_{Q_f,\pi}$ are smooth in the first and last line, and the fact that $p^+$ is an affine fiber bundle in the second line. This isomorphism lift to an isomorphism of mixed Hodge modules. We denote now by $\phi_{W^0}$ the vanishing cycles functor of the restriction $\Tr(W)|_{\mathcal{M}^0_{Q_f,d}}$. The natural functoriality of the vanishing cycle functor gives a morphism $(p^+)_!(\eta^+)^\ast\phi_W\to\phi_{W^0}(p^+)_!(\eta^+)^\ast$, which lifts to a morphism of monodromic mixed Hodge modules. Then, because $\mathcal{IC}_{\mathcal{M}_{Q_f,d}}$ is $\mathbb{C}^\ast$-equivariant, \cite[theo 3.3]{Ric16} gives that:
\begin{align}\label{hyplocvan}
    (p^+)_!(\eta^+)^\ast\phi_W\mathcal{IC}_{\mathcal{M}_{Q_f,d}}\to\phi_{W^0}(p^+)_!(\eta^+)^\ast\mathcal{IC}_{\mathcal{M}_{Q_f,d}}
\end{align}
is an isomorphism at the level of perverse sheaves. Because the fact of being an isomorphism can be checked at the level of the underlying perverse sheaves, it is also an isomorphism of monodromic mixed Hodge modules. Hence one obtains, denoting by $q:\mathcal{M}_{Q_f,d}^0\to\ast$ the projection to a point:
\begin{align}
    [\mathcal{M}^+_{Q_f,d}]^{vir}&=H_c(\mathcal{M}^+_{Q_f,d},\phi_W\mathcal{IC}_{\mathcal{M}_{Q_f,d}})\nn\\
    &=q_!(p^+)_!(\eta^+)^\ast\phi_W\mathcal{IC}_{\mathcal{M}_{Q_f,d}}\nn\\
    &=q_!\phi_{W^0}(p^+)_!(\eta^+)^\ast\mathcal{IC}_{\mathcal{M}_{Q_f,d}}\nn\\
    &=\bigoplus_{\pi\in\Pi_f|d_\pi=d}\mathbb{L}^{\Index^s_\pi/2}q_!\phi_{W^0}\mathcal{IC}_{\mathcal{M}^0_{Q_f,\pi}}\nn\\
    &=\bigoplus_{\pi\in\Pi_f|d_\pi=d}\mathbb{L}^{\Index^s_\pi/2}H_c(\mathcal{M}^0_{Q_f,\pi},\phi_{W^0}\mathcal{IC}_{\mathcal{M}^0_{Q_f,\pi}})
\end{align}
Here we have used \eqref{hyplocvan} in the third line and \eqref{BB} in the fourth line. The perverse sheaf $\phi_{W^0}\mathcal{IC}_{\mathcal{M}^0_{Q_f,\pi}}$ is supported on the critical locus of $\Tr(W)|_{\mathcal{M}^0_{Q_f,\pi}}$, which is just a single point, the representation associated to $\pi$. Hence the cohomology of the vanishing cycle on this point is just the Milnor number of this point, which is $1$ from \cite[Prop 3.3]{Behrend08symmetricobstruction}, \ie:
\begin{align}\label{Kform}
    [\mathcal{M}^+_{Q_f,d}]^{vir}=\bigoplus_{\pi\in\Pi_f|d_\pi=d}\mathbb{L}^{\Index^s_\pi/2}
\end{align}
Now using the Lemma \ref{lemattr}, one obtains in the case $i)$:
\begin{align}
    Z_D(x):&=\sum_d[\mathcal{M}^+_{Q_D,d}]^{vir}x^d\nn\\
    &=\bigoplus_{\pi\in\Pi_f}\mathbb{L}^{\Index^s_\pi/2}x^{d_\pi}
\end{align}
And in the case $ii)$:
\begin{align}
    Z_i^{[z,z']:N}(x):&=\sum_d[\mathcal{M}^+_{Q_i,d}]^{vir}x^d\nn\\
    &=\bigoplus_{\pi\in\Pi_f}\mathbb{L}^{\Index^s_\pi/2}x^{d_\pi}
\end{align}
Because $\sum_d\Delta^s\Omega_dx^d=\Omega_\theta(x)-\Omega^{[z,z']:N}_\theta(x)$ lies in the center of the quantum affine space, one has:
\begin{align}
    \mathcal{A}(x)=\Exp(\sum_d\frac{\Delta^s\Omega_d}{\mathbb{L}^{1/2}-\mathbb{L}^{-1/2}}x^d)\mathcal{A}^{[z,z']:N}(x)
\end{align}
And using Proposition \ref{nilpart}:
\begin{align}
    Z_i(x)=&S_i(\mathcal{A}(x))S_{-i}(\mathcal{A}(x)^{-1})\nn\\
    =&S_i(\Exp(\textstyle\sum_d\frac{\Delta^s\Omega_d}{\mathbb{L}^{1/2}-\mathbb{L}^{-1/2}}x^d)\mathcal{A}^{[z,z']:N}(x))S_{-i}((\Exp(\textstyle\sum_d\frac{\Delta^s\Omega_d}{\mathbb{L}^{1/2}-\mathbb{L}^{-1/2}}x^d)\mathcal{A}^{[z,z']:N}(x))^{-1})\nn\\
    =&S_i(\Exp(\textstyle\sum_d\frac{\Delta^s\Omega_d}{\mathbb{L}^{1/2}-\mathbb{L}^{-1/2}}x^d))S_{-i}(\Exp(-\textstyle\sum_d\frac{\Delta^s\Omega_d}{\mathbb{L}^{1/2}-\mathbb{L}^{-1/2}}x^d))S_i(\mathcal{A}^{[z,z']:N}(x))S_{-i}(\mathcal{A}^{[z,z']:N}(x)^{-1})\nn\\
    =&S_{-i}[\Exp\left(\textstyle\sum_d\Delta^{s}\Omega_d\frac{\mathbb{L}^{d_i}-1}{\mathbb{L}^{1/2}-\mathbb{L}^{-1/2}}x^d\right)]Z^{[z,z']:N}_i(x)
\end{align}
Here we have used once more in the third line the fact that $\sum_d\Delta^s\Omega_dx^d$ lies in the center of the quantum affine space. $\Box$
\medskip

\subsection{Link with K-theoretic computations}\label{secKtheo}

The localization formula \eqref{Kform} is very similar to the localization formula for K-theoretic invariants defined in \cite{NekOk}. Those invariants are defined for projective spaces with symmetric obstruction theories, hence in particular for projective critical locus of a potential on a smooth space, and are expected to give the $\chi_y$ genus of the Hodge polynomial coming from cohomological invariants. Consider a moduli space $M$ which is the critical locus of a potential, with a $\mathbb{C}^\ast$-action leaving the potential invariant (the choice of such an action is called a choice of slope $s$ in K-theoretic theory). We denote the attracting variety by $M^+$ and the fixed components of this $\mathbb{C}^\ast$-action by $M^0_\pi$ for $\pi\in\Pi$, and denotes as before by $\Index^s_\pi$ the signed number of contracting weight in the restriction of the tangent-obstruction complex of $M$ at $M^0_\pi$. Because $M$ is projective, the attracting variety is $M$, hence the same reasoning that leads to \eqref{Kform} gives:
\begin{align}\label{locnonproj}
    [M^+]^{vir}=\bigoplus_{\pi\in\Pi}\mathbb{L}^{\Index^s_\pi/2}[M^0_\pi]^{vir}
\end{align}
In fact, the author proved in \cite{Descombes:2022cpc} that this formula holds also when $M$ is a [-1]-shifted symplectic scheme or stack, \ie is locally described as the critical locus of a potential. If $M$ is projective, $M^+=M$, and then taking the $\chi_y$-genus gives:
\begin{align}\label{locproj}
    \chi_y([M]^{vir})=\sum_{\pi\in\Pi}(-y)^{\Index^s_\pi} \chi_y([M^0_\pi]^{vir})
\end{align}
It is exactly the localization formula of K-theoretic invariants from \cite[Sec 8.3]{NekOk}. It can be applied to cases of framed quiver with potential when the moduli space is projective, for example when the Empty Room Configuration has a finite number of atoms, as in \cite{Cir19}\medskip

It was proposed to define K-theoretic invariants by the formula \eqref{locproj} when $M^0$ is projective, but $M$ is potentially non-projective. However this definition depends on the choice of slope, \ie on the choice of $\mathbb{C}^\ast$ action. This dependency is explained by the formula \eqref{locnonproj}, because the attracting variety depends on the slope. For toric Calabi-Yau threefolds, the dependency of K-theoretic invariants of framed sheaves was studied in \cite[Prop 3.3]{Arb}, and it was established that the K-theoretic invariants change only when a toric coordinates becomes attracting or repelling. Framed sheaves corresponds with D6-framed representations of the toric quiver with potential, and the cycles $v^z$ are scaled as the toric coordinates of the quiver. Hence \cite[Prop 3.3]{Arb} is coherent with the formula \eqref{locnonproj}, because the attracting variety changes only when a weight $\lambda_z$ becomes attracting or repelling. Hence Theorem \ref{theolocD} explains the discrepancy between K-theoretic and cohomological/motivic invariants observed for toric quivers, as we will check in several cases in Section \ref{localcurves}.

\section{Examples of toric quivers}\

In this section, we compare our results with the known BPS invariants for local curves, \ie when the number $i$ of internal lattice points in the toric diagram vanishes, and spell out our results and conjecture for local toric surfaces, \ie for $i=1$.

\subsection{Local curves}\label{localcurves}

In those cases, the quantum affine space is commutative, there is no wall crossing, \ie the BPS invariants are independent of $\theta$.  The generating series of cohomological invariants were explicitly computed. Moreover, depending on the toric diagram, there can be 'preferred slope' as introduced in \cite{IKV}, for which there is many cancellations in $\Index^s_\pi$: this index becomes then a sum of simple contributions for each atom of $\pi$, and the localization formula gives a closed expression. We show then the agreement between the cohomological computations and the localization computations of \cite{IKV}, corrected as in Theorem \ref{theolocD}.
\medskip

$\bullet$ $\mathbb{C}^3$\medskip

The toric diagram of $\mathbb{C}^3$ is given by Figure \ref{torc3}. We label by $z,z',z''$ the three edges of the toric diagram in the clockwise order, and consider a slope $s$ such that $z\in L^{<0}$ and $z',z''\in L^{>0}$. Then according to the discussion in Section \ref{seccompdiv}, the refined generating series of framed invariant computed from K-theoretic localization in \cite[sec 8.3]{NekOk} using the slope $s$ (resp $-s$) correspond in our formalism to the generating series $Z_i^{z:N}$ (resp $Z_i^{z':N,z'':N}$). Then equation the result \cite[sec 8.3]{NekOk} can be expressed as:
\begin{align}
    &\Omega^{z:N}(x)=\mathbb{L}^{1/2}\sum_{n\geq 1}x^{n\delta}\nn\\
    &\Omega^{z':N,z'':N}(x)=\mathbb{L}^{-1/2}\sum_{n\geq 1}x^{n\delta}
\end{align}
Using the correction given in Proposition \ref{nilpart}, one obtains for any of these two choices of slope:
\begin{align}\label{C3check}
    \Omega(x)=\mathbb{L}^{3/2}\sum_{n\geq 1}x^{n\delta}
\end{align}
in perfect agreement  with the result of \cite[Theo 2.7]{BBS}. Note that the symmetric part  vanishes in this case.
\medskip

$\bullet$ $\mathbb{C}^3/(\mathbb{Z}_2\times\mathbb{Z}_2)$\medskip

\begin{figure}
\caption{Toric diagram and brane tiling for $\mathbb{C}^3/(\mathbb{Z}_2\times\mathbb{Z}_2)$}\label{tordi}
\medskip
\centerline{\hfill
\begin{tikzpicture}[rotate=45, yshift=2cm]
\filldraw [black] (0,0) circle (2pt);
\draw (0,0) node[below]{$p_0$};
\filldraw [black] (1,0) circle (2pt);
\filldraw [black] (2,0) circle (2pt);
\draw (2,0) node[right]{$p_2$}; 
\filldraw [black] (0,1) circle (2pt);
\filldraw [black] (1,1) circle (2pt);
\filldraw [black] (0,2) circle (2pt);
\draw (0,2) node[left]{$p_1$};
\draw (0,0) -- (2,0); 
\draw (2,0) -- (0,2);
\draw (0,0) -- (0,2);
\draw (1,1) -- (1,0); 
\draw (1,0) -- (0,1);
\draw (1,1) -- (0,1);
\end{tikzpicture}
\hfill
\begin{tikzpicture} 
 \foreach \x in {0,...,2}{
      \foreach \y in {0,...,2}{
        \node[draw,circle,inner sep=2pt,fill] at (0+1.5*\x,0+1*\y){};
        \node[draw,circle,inner sep=2pt] at (0.5+1.5*\x,0+1*\y){};
        \node[draw,circle,inner sep=2pt,fill] at (0.75+1.5*\x,0.5+1*\y){};
        \node[draw,circle,inner sep=2pt] at (-0.25+1.5*\x,0.5+1*\y){};
        \draw (0+1.5*\x,0+1*\y) -- (0.5+1.5*\x,0+1*\y)[gray!60];
        \draw (0.75+1.5*\x,0.5+1*\y) -- (1.25+1.5*\x,0.5+1*\y)[gray!60];
        \draw (0+1.5*\x,0+1*\y) -- (-0.25+1.5*\x,0.5+1*\y)[green];
        \draw (0.75+1.5*\x,0.5+1*\y) -- (0.5+1.5*\x,1+1*\y)[green];
        \draw (0.5+1.5*\x,0+1*\y) -- (0.75+1.5*\x,0.5+1*\y)[red];
        \draw (-0.25+1.5*\x,0.5+1*\y) -- (0+1.5*\x,1+1*\y)[red];
        \ADD{\x}{\y}{\s}
        \MODULO{\s}{2}{\so};
        \ADD{\so}{2}{\sol};
        \node at (0.25+1.5*\x,0.5+1*\y) {\so};
        \node at (1+1.5*\x,0+1*\y) {\sol};
        }}
\end{tikzpicture}
\hfill}
\end{figure}
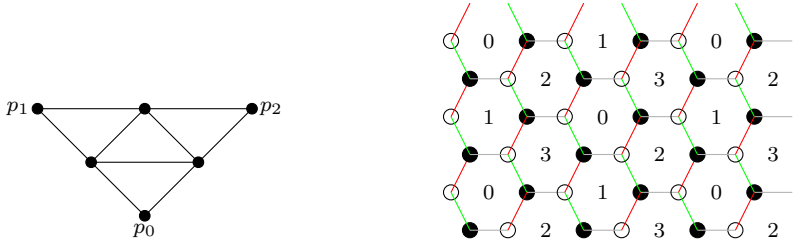
\medskip

The toric diagram and perfect matchings are represented in the Figure \ref{tordi}. Each external edge of the toric diagram has two subdivisions. In green and red, we have written the zig-zag path corresponding with the edge $z_{3/2}$ between $p_1$ and $p_2$. It divides the quiver into two sub-quivers $Q^0$ and $Q^1$, with nodes $\{0,1\}$ and $\{2,3\}$. The two other external edges give zig-zag paths that are similar but rotated by an angle $\pm2\pi/3$, dividing the quiver respectively into sub-quivers with nodes $\{0,3\}$ and $\{1,2\}$, resp. $\{0,2\}$ and $\{1,3\}$. 
Our computation of the anti-symmetric part of the attractor invariants gives (denoting $x_{i_1i_2...i_r}=x_{i_1}x_{i_2}...x_{i_r}$):
\begin{align}
    \Omega(x)=((\mathbb{L}^{3/2}+3\mathbb{L}^{1/2})x^\delta+\mathbb{L}^{1/2}(x_{01}+x_{23}+x_{03}+x_{12}+x_{02}+x_{13}))\sum_{n\geq 0}x^{n\delta}+\Omega^{sym}(x)
\end{align}
It agrees with the result of \cite[Remark 5.2]{MR21}:
\begin{align}
    \Omega(x)=&((\mathbb{L}^{3/2}+3\mathbb{L}^{1/2})x^\delta+\mathbb{L}^{1/2}(x_{01}+x_{23}+x_{03}+x_{12}+x_{02}+x_{13})\nn\\
    &+1(x_0+x_1+x_2+x_3+x_{123}+x_{230}+x_{301}+x_{012}))\sum_{n\geq 0}x^{n\delta}
\end{align}
\medskip

$\bullet$ Other small crepant resolutions:

The other toric small crepant resolutions are resolutions of the zero locus of $XY-Z^{N_0}W^{N_1}$ in $\mathbb{C}^4$. The corresponding toric diagram has a trapezoidal shape with height $1$, a lower edge of length $N_0$, and upper edge of length $N_1$. A noncommutative resolution of this threefold is determined by a triangulation $\sigma$ of the toric diagram. The construction of the corresponding quiver and brane tiling is described in \cite[sec 1.1]{Nag11}. We enumerate triangles by $T_i$ from the right to the left, for $i\in \mathcal{I}=\mathbb{Z}_N$, for $N=N_0+N_1$ ( in particular $b=N+2$), cyclically identifying the right external edge of the toric diagram with the left external edge of the toric diagram. The triangulation defines a bijection:
\begin{align}
    \sigma=(\sigma_x,\sigma_y):I_N=\{0,...,N-1\}\to (I_{N_0}\times \{0\})\cup(I_{N_1}\times \{1\})
\end{align}
We define:
\begin{align}
    \mathcal{J}=\{i\in \mathcal{I}|\sigma_y(i)=\sigma_y(i+1)\}
\end{align}
which enumerates $i\in \mathcal{I}$ such that triangles $T_i$ and $T_{i+1}$ have adjoint horizontal edges (we consider triangles $T_{N-1}$ and $T_0$ for $i=N-1$.
\medskip

We construct then a quiver with nodes $\mathcal{I}$, a pair of bidirectional arrows between successive nodes $i,i+1$, and an edge loop at nodes of $\mathcal{J}$. The corresponding brane tiling is obtained by stacking up layers of the form:

\medskip
\centerline{
\begin{tikzpicture} 
 \foreach \x in {0,...,5}{
    \draw (\x,0)--(\x+0.5,0.25);
    \draw (\x,0)--(\x+0.5,-0.25);
    \draw (\x+1,0)--(\x+0.5,0.25);
    \draw (\x+1,0)--(\x+0.5,-0.25);
    \node at (\x+0.5,0) {$i$};
    }
\end{tikzpicture}
}
if $i\in \mathcal{I}-\mathcal{J}$, and of the form:

\medskip
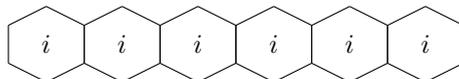
\centerline{\begin{tikzpicture} 
 \foreach \x in {0,...,5}{
    \draw (\x,0)--(\x,0.5);
    \draw (\x,0)--(\x+0.5,-0.25);
    \draw (\x+1,0)--(\x+1,0.5);
    \draw (\x+1,0)--(\x+0.5,-0.25);
    \draw (\x,0.5)--(\x+0.5,0.75);
    \draw (\x+1,0.5)--(\x+0.5,0.75);
    \node at (\x+0.5,0.25) {$i$};
    }
\end{tikzpicture}}

if $i\in \mathcal{I}$.

The zig-zag paths corresponding to the lower (resp. upper) edge of the toric diagram, denoted $z^0$ (resp. $z^1$) are given by the border between two successive layers $i-1,i$ such that $\sigma_y(T_i)=0$ (resp. $\sigma_y(T_i)=1$). 
\medskip

\begin{figure}
\caption{Example of triangulation corresponding to a small crepant resolution}\label{trig}
\centerline{\begin{tikzpicture} 
\filldraw [black] (0,0) circle (2pt);
\draw (0,0) node[left]{$p_0$};
\filldraw [black] (2,0) circle (2pt);
\filldraw [black] (4,0) circle (2pt);
\filldraw [black] (6,0) circle (2pt);
\filldraw [black] (8,0) circle (2pt);
\draw (8,0) node[right]{$p_3$}; 
\filldraw [black] (0,2) circle (2pt);
\draw (0,2) node[left]{$p_1$}; 
\filldraw [black] (2,2) circle (2pt);
\filldraw [black] (4,2) circle (2pt);
\draw (4,2) node[right]{$p_2$};
\draw (0,0) -- (8,0);
\draw (0,2) -- (4,2);
\draw (0,0) -- (0,2);
\draw (0,0) -- (2,2);
\draw (2,0) -- (2,2);
\draw (2,0) -- (4,2);
\draw (4,0) -- (4,2);
\draw (6,0) -- (4,2);
\draw (8,0) -- (4,2);
\draw (0,1) node[left]{$5$};
\draw (1,1) node[left]{$4$};
\draw (2,1) node[left]{$3$};
\draw (3,1) node[left]{$2$};
\draw (4,1) node[left]{$1$};
\draw (5,1) node[left]{$0$};
\draw (6,1) node[left]{$5$};
\draw (0.5,1.5) node {$T_5$};
\draw (1.5,0.5) node {$T_4$};
\draw (2.5,1.5) node {$T_3$};
\draw (3.5,0.5) node {$T_2$};
\draw (4.5,0.5) node {$T_1$};
\draw (6.5,0.5) node {$T_0$};
\end{tikzpicture}}
\end{figure}
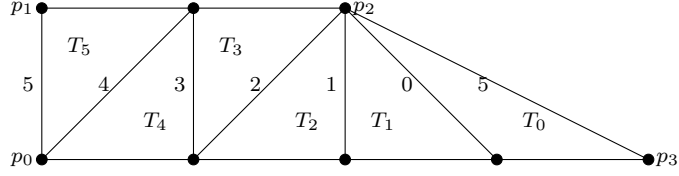

\begin{figure}\caption{Corresponding brane tiling}\label{bratrig}\centerline{
\begin{tikzpicture} 
 \foreach \x in {0,...,5}{
    \node[draw,circle,inner sep=2pt,fill] at (\x+0.5,-0.25) {};
    \draw (\x,0) -- (\x+0.5,-0.25)[blue];
    \draw (\x+1,0) -- (\x+0.5,-0.25)[blue];
    \node[draw,circle,inner sep=2pt] at (\x,0) {};
    \draw (\x,0) -- (\x+0.5,0.25)[red];
    \draw (\x+1,0) -- (\x+0.5,0.25)[red];
    \node[draw,circle,inner sep=2pt,fill] at (\x+0.5,0.25) {};
    \draw (\x,0.5) -- (\x+0.5,0.25)[blue];
    \draw (\x+1,0.5) -- (\x+0.5,0.25)[blue];
    \node[draw,circle,inner sep=2pt] at (\x,0.5) {};
    \draw (\x,0.5) -- (\x+0.5,0.75)[red];
    \draw (\x+1,0.5) -- (\x+0.5,0.75)[red];
    \node[draw,circle,inner sep=2pt,fill] at (\x+0.5,0.75) {};
    \draw (\x,1) -- (\x+0.5,0.75)[blue];
    \draw (\x+1,1) -- (\x+0.5,0.75)[blue];
    \node[draw,circle,inner sep=2pt] at (\x,1) {};
    \draw (\x,1) -- (\x,1.5)[gray!60];
    \node[draw,circle,inner sep=2pt,fill] at (\x,1.5) {};
    \draw (\x,1.5) -- (\x+0.5,1.75)[blue];
    \draw (\x+1,1.5) -- (\x+0.5,1.75)[blue];
    \node[draw,circle,inner sep=2pt] at (\x+0.5,1.75) {};
    \draw (\x+0.5,1.75) -- (\x+0.5,2.25)[gray!60];
    \node[draw,circle,inner sep=2pt,fill] at (\x+0.5,2.25) {};
    \draw (\x+0.5,2.25) -- (\x+1,2.5)[blue];
    \draw (\x+0.5,2.25) -- (\x,2.5)[blue];
    \node[draw,circle,inner sep=2pt] at (\x,2.5) {};
    \draw (\x,2.5) -- (\x+0.5,2.75)[red];
    \draw (\x+1,2.5) -- (\x+0.5,2.75)[red];
    \node[draw,circle,inner sep=2pt,fill] at (\x+0.5,2.75) {};
    \node at (\x+0.5,0) {$5$};
    \node at (\x,0.25) {$4$};
    \node at (\x+0.5,0.5) {$3$};
    \node at (\x,0.75) {$2$};
    \node at (\x+0.5,1.25) {$1$};
    \node at (\x,2) {$0$};
    \node at (\x+0.5,2.5) {$5$};
    }
\end{tikzpicture}
}
\end{figure}
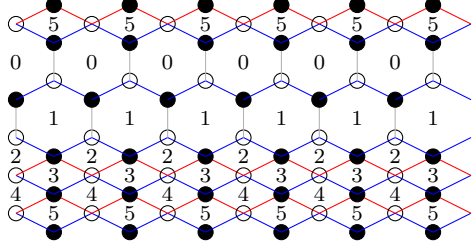
\medskip

As an example, consider the triangulation of Figure \ref{trig}, for $N_0=4$, $N_1=2$, $N=6$: We have $\sigma=((3,0),(2,0),(1,0),(1,1),(0,0),(0,1))$, $\mathcal{I}=\mathbb{Z}_6$, $\mathcal{J}=\{0,1\}$. The corresponding brane tiling is given by Figure \ref{bratrig}, where we have drawn  in red the zig-zag paths $z^1=z_{3/2}$ corresponding to the upper edge  between $p_1$ and $p_2$, and in blue the the zig-zag paths $z^0=z_{-1/2}$ corresponding to the lower edge between $p_3$ and $p_0$.
\medskip

We can use our evaluation of the anti-symmetric part of the attractor invariants. To this aim we must find all the dimension vectors that are of the form $\alpha^z_{[k,k'[}$, for $z$ a side of the toric diagram and $k\neq k'\in\mathbb{Z}/K_z\mathbb{Z}$. The left and right external edges of the toric diagram have only one subdivision, \ie they do not give such roots. According to our description of zig-zag paths corresponding to the above and below side of the toric diagram, such a dimension vector is of the form $d=e_j+e_{j+1}+...+e_{l}$, $l+1\neq j$, (in particular it is in the set of real root $\Delta^{re}_+$), such that the corresponding layers of the brane tiling lie between two zig-zag paths of $z^0$ or two zig-zag paths of $z^1$, \ie such that $\sigma_y(T_j)=\sigma_y(T_{l+1})$. According to our definition of $\mathcal{J}$, this holds if and only if $\sum_{i\not\in \mathcal{J}}d_i$ is even. We obtain:
\begin{align}
    \Omega(x)=(\mathbb{L}^{3/2}+(N-1)\mathbb{L}^{1/2})\sum_{d\in\Delta^{im}_+} x^d+\mathbb{L}^{1/2}\sum_{d\in\Delta^{re}_+|\sum_{i\not\in \mathcal{J}}d_i\: even}x^d+\Omega^{sym}(x)
\end{align}

This is in agreement with \cite[Theo 0.1]{MorNag}, upon adding a suitable symmetric correction:
\begin{align}
    \Omega(x)=(\mathbb{L}^{3/2}+(N-1)\mathbb{L}^{1/2})\sum_{d\in\Delta^{im}_+} x^d+\mathbb{L}^{1/2}\sum_{d\in\Delta^{re}_+|\sum_{i\not\in \mathcal{J}}d_i\: even}x^d+1\sum_{d\in\Delta^{re}_+|\sum_{i\not\in \mathcal{J}}d_i\: odd}x^d
\end{align}\medskip

$\bullet$ The conifold\medskip

It is a particular case of small crepant resolutions for $N_0=N_1=1$. Then \cite[Theo 1]{MMNS} gives:
\begin{align}
    \Omega(x)&=(\mathbb{L}^{3/2}+\mathbb{L}^{1/2})\sum_{n\geq1}x^{n\delta}+\sum_{n\geq 0}(x_0+x_1)x^{n\delta}
\end{align}
We label by $z_1,z_2,z_3,z_4$ the edges of the toric diagram. When one choose a generic slope, one can consider up to a circular permutation of the sides of the diagram that $z\in L^{<0}\iff z\in[z_1,z_2]$. By toric localization, one access then to the BPS invariants $\Omega^{[z_1,z_2]:N}(x)$. Every edge has only one subdivision, hence the full and partially invertible invariants are related by:
\begin{align}
    \Omega(x)&=(\mathbb{L}^{3/2}-\mathbb{L}^{-1/2})\sum_{n\geq1}x^{n\delta}+\Omega^{[z_1,z_2]:N}(x)\nn\\
    \implies\Omega^{[z_1,z_2]:N}(x)&=(\mathbb{L}^{1/2}+\mathbb{L}^{-1/2})\sum_{n\geq1}x^{n\delta}+\sum_{n\geq 0}(x_0+x_1)x^{n\delta}\label{eqcon}
\end{align}
It is in agreement with the computations of the refined topological vertex in \cite[Sec 5.1]{IKV}, as exposed in \cite[Sec 4.3]{MMNS}. Namely, in \cite[Sec 5.1]{IKV} the refined topological vertex is used to compute the 'PT partition function' which gives access to the invariants $\Omega^{[z_1,z_2]:N}_{(d_1,d_2)}$ for $d_1>d_2$, hence using symmetries for all the invariants for $d\not\in\mathbb{N}\delta$, and those invariants are in agreement with \eqref{eqcon} (and coincide in fact with the motivic invariants $\Omega_d$). The invariants $\Omega^{[z_1,z_2]:N}_{n\delta}$ and $\Omega_{n\delta}$ are different: it explains the observation of \cite[p. 2]{MMNS} hat the motivic and refined computations agree 'up to a subtlety involving the Hilbert scheme of points', and the ambiguity in defining the refinement of the MacMahon function in \cite{Dimofte:2009bv}.\medskip

$\bullet$ $\mathbb{C}^2/(\mathbb{Z}/2\mathbb{Z})\times\mathbb{C}$\medskip

It is a particular case of small toric crepant resolution for $N_0=2$ and $N_1=0$. Then \cite[Theo 0.1]{MMNS} gives:
\begin{align}
  \Omega(x)=(\mathbb{L}^{3/2}+\mathbb{L}^{1/2})\sum_{n\geq1}x^{n\delta}+\mathbb{L}^{1/2}\sum_{n\geq 0}(x_0+x_1)x^{n\delta}
\end{align}
We label by $z,z',z''$ the edges of the toric diagram, $z$ being the edge with two subdivisions. The slope in \cite[Sec 5.3, Fig. 6)b)]{IKV} gives $z''\in L^{<0}$ and $z,z'\in L^{>0}$, hence the refined topological gives access to the BPS invariants $\Omega^{[z_1,z_2]:N}$. The edge $z''$ has only one subdivision, hence:
\begin{align}
    \Omega(x)&=(\mathbb{L}^{3/2}-\mathbb{L}^{1/2})\sum_{n\geq1}x^{n\delta}+\Omega^{[z_1,z_2]:N}(x)\nn\\
    \implies\Omega^{[z_1,z_2]:N}(x)&=2\mathbb{L}^{1/2}\sum_{n\geq1}x^{n\delta}+\mathbb{L}^{1/2}\sum_{n\geq 0}(x_0+x_1)x^{n\delta}
\end{align}
It is in agreement with the computations of the refined topological vertex in \cite[Sec 5.3]{IKV}.\medskip

\subsection{Toric threefolds with one compact divisors}\label{seccompdiv}\medskip

$\bullet$ Canonical bundle over toric Fano surfaces:
\medskip

In this case, the toric diagram has one internal lattice point and the only points on the boundary are the corners, \ie external edges have only one subdivision. Our result gives then:
\begin{align}
    \Omega_\ast(x)=(\mathbb{L}^{3/2}+(b-2)\mathbb{L}^{1/2}+\mathbb{L}^{-1/2})\sum_{n\geq 1}x^{n\delta}+\Omega^{sym}_\ast(x)
\end{align}
The arguments of \cite{BMP20}, and explicit computations for small dimension vectors done in \cite[sec 6]{MP20}, support our conjectural formula \ref{conj}:
\begin{align}
    \Omega_\ast(x)=\sum_i x_i+(\mathbb{L}^{3/2}+(b-2)\mathbb{L}^{1/2}+\mathbb{L}^{-1/2})\sum_{n\geq 1}x^{n\delta}
\end{align}

$\bullet$ Canonical bundle over toric weak Fano surfaces:\medskip

In those cases, the toric diagram has one internal lattice point, and its external edges can have various number of subdivisions. For completeness, we will give here our conjectural formula \ref{conj} (which is proven up to a symmetric correction) for those various geometries, using the notations of \cite{HanSeo} (when there is some misprint in this reference, we use the label of the nodes given in the brane tiling):
\medskip

$\ast$ $F_2$ (model 13 of \cite{HanSeo})
\begin{align}
    \Omega_\ast=\sum_ix_i+((\mathbb{L}^{3/2}+2\mathbb{L}^{1/2}+\mathbb{L}^{-1/2})x^\delta+\mathbb{L}^{1/2}(x_{13}+x_{24}))\sum_{n\geq 0}x^{n\delta}
\end{align}

$\ast$ $PdP_2$ (model 11 of \cite{HanSeo})
\begin{align}
    \Omega_\ast=\sum_ix_i+((\mathbb{L}^{3/2}+3\mathbb{L}^{1/2}+\mathbb{L}^{-1/2})x^\delta+\mathbb{L}^{1/2}(x_{12}+x_{345}))\sum_{n\geq 0}x^{n\delta}
\end{align}

$\ast$ $PdP_{3b}$ (model 9 of \cite{HanSeo})
\begin{align}
    &phase\:a:\quad \Omega_\ast=\sum_ix_i+((\mathbb{L}^{3/2}+4\mathbb{L}^{1/2}+\mathbb{L}^{-1/2})x^\delta+\mathbb{L}^{1/2}(x_{12}+x_{3456}))\sum_{n\geq 0}x^{n\delta}\nn\\
    &phase\:b:\quad \Omega_\ast=\sum_ix_i+((\mathbb{L}^{3/2}+4\mathbb{L}^{1/2}+\mathbb{L}^{-1/2})x^\delta+\mathbb{L}^{1/2}(x_{126}+x_{345}))\sum_{n\geq 0}x^{n\delta}\nn\\
    &phase\:c:\quad \Omega_\ast=\sum_ix_i+((\mathbb{L}^{3/2}+4\mathbb{L}^{1/2}+\mathbb{L}^{-1/2})x^\delta+\mathbb{L}^{1/2}(x_{1246}+x_{35}))\sum_{n\geq 0}x^{n\delta}\nn\\
\end{align}

$\ast$ $PdP_{3c}$ (model 8 of \cite{HanSeo})
\begin{align}
    &phase\:a:\quad \Omega_\ast=\sum_ix_i+((\mathbb{L}^{3/2}+4\mathbb{L}^{1/2}+\mathbb{L}^{-1/2})x^\delta+\mathbb{L}^{1/2}(x_{126}+x_{345}+x_{1234}+x_{56}))\sum_{n\geq 0}x^{n\delta}\nn\\
    &phase\:b:\quad \Omega_\ast=\sum_ix_i+((\mathbb{L}^{3/2}+4\mathbb{L}^{1/2}+\mathbb{L}^{-1/2})x^\delta+\mathbb{L}^{1/2}(x_{156}+x_{234}+x_{1345}+x_{26}))\sum_{n\geq 0}x^{n\delta}
\end{align}

$\ast$ $PdP_{3a}$ (model 7 of \cite{HanSeo})
\begin{align}
    \Omega_\ast=&\sum_ix_i+((\mathbb{L}^{3/2}+4\mathbb{L}^{1/2}+\mathbb{L}^{-1/2})x^\delta+\mathbb{L}^{1/2}(x_{124}+x_{356}+x_{15}+x_{34}+x_{26}\nn\\&+x_{15}x_{34}+x_{34}x_{26}+x_{26}x_{15}))\sum_{n\geq 0}x^{n\delta}
\end{align}

$\ast$ $PdP_{4a}$ (model 6 of \cite{HanSeo})
\begin{align}
    &phase\:a:\quad \Omega_\ast=\sum_ix_i+((\mathbb{L}^{3/2}=5\mathbb{L}^{1/2}+\mathbb{L}^{-1/2})x^\delta+\mathbb{L}^{1/2}(x_{137}+x_{2456}+x_{1345}+x_{267}))\sum_{n\geq 0}x^{n\delta}\nn\\
    &phase\:b:\quad \Omega_\ast=\sum_ix_i+((\mathbb{L}^{3/2}+5\mathbb{L}^{1/2}+\mathbb{L}^{-1/2})x^\delta+\mathbb{L}^{1/2}(x_{12345}+x_{67}+x_{1237}+x_{456}))\sum_{n\geq 0}x^{n\delta}\nn\\
    &phase\:c:\quad \Omega_\ast=\sum_ix_i+((\mathbb{L}^{3/2}+5\mathbb{L}^{1/2}+\mathbb{L}^{-1/2})x^\delta+\mathbb{L}^{1/2}(x_{167}+x_{2345}+x_{1456}+x_{237}))\sum_{n\geq 0}x^{n\delta}\nn\\
\end{align}

$\ast$ $PdP_{4b}$ (model 5 of \cite{HanSeo})
\begin{align}
     \Omega_\ast=&\sum_ix_i+((\mathbb{L}^{3/2}+5\mathbb{L}^{1/2}+\mathbb{L}^{-1/2})x^\delta+\mathbb{L}^{1/2}(x_{1234}+x_{567}+x_{17}+x_{26}+x_{345}+x_{17}x_{26}\nn\\&+x_{26}x_{345}+x_{345}x_{17}))\sum_{n\geq 0}x^{n\delta}
\end{align}

$\ast$ $PdP_{5}$ (model 4 of \cite{HanSeo})
\begin{align}
     phase\:a:\quad\Omega_\ast=&\sum_ix_i+((\mathbb{L}^{3/2}+6\mathbb{L}^{1/2}+\mathbb{L}^{-1/2})x^\delta+\mathbb{L}^{1/2}(x_{1234}+x_{5678}+x_{1638}+x_{2745}\nn\\&+x_{1674}+x_{2385}+x_{1278}+x_{3456}))\sum_{n\geq 0}x^{n\delta}\nn\\
     phase\:b:\quad\Omega_\ast=&\sum_ix_i+((\mathbb{L}^{3/2}+6\mathbb{L}^{1/2}+\mathbb{L}^{-1/2})x^\delta+\mathbb{L}^{1/2}(x_{46}+x_{123578}+x_{28}+x_{134567}\nn\\&+x_{1568}+x_{2347}+x_{1245}+x_{3678}))\sum_{n\geq 0}x^{n\delta}\nn\\
     phase\:c:\quad\Omega_\ast=&\sum_ix_i+((\mathbb{L}^{3/2}+6\mathbb{L}^{1/2}+\mathbb{L}^{-1/2})x^\delta+\mathbb{L}^{1/2}(x_{234}+x_{15678}+x_{368}+x_{12457}\nn\\&+x_{278}+x_{13456}+x_{467}+x_{12358}))\sum_{n\geq 0}x^{n\delta}\nn\\
     phase\:d:\quad\Omega_\ast=&\sum_ix_i+((\mathbb{L}^{3/2}+6\mathbb{L}^{1/2}+\mathbb{L}^{-1/2})x^\delta+\mathbb{L}^{1/2}(x_{1678}+x_{2345}+x_{1247}+x_{3568}\nn\\&+x_{1346}+x_{2578}+x_{1238}+x_{4567}))\sum_{n\geq 0}x^{n\delta}
\end{align}

$\ast$ $L_{1,3,1}/\mathbb{Z}_2 (0,1,1,1)$ (model 3 of \cite{HanSeo})
\begin{align}
    phase\:a:\quad \Omega_\ast=&\sum_ix_i+((\mathbb{L}^{3/2}+6\mathbb{L}^{1/2}+\mathbb{L}^{-1/2})x^\delta+\mathbb{L}^{1/2}(x_{1453}+x_{2786}+x_{1756}+x_{2483}\nn\\&+x_{18}+x_{37}+x_{2456}+x_{18}x_{37}+x_{37}x_{2456}+x_{2456}x_{18}))\sum_{n\geq 0}x^{n\delta}\nn\\
    phase\:b:\quad \Omega_\ast=&\sum_ix_i+((\mathbb{L}^{3/2}+6\mathbb{L}^{1/2}+\mathbb{L}^{-1/2})x^\delta+\mathbb{L}^{1/2}(x_{14253}+x_{678}+x_{17256}+x_{348}\nn\\&+x_{18}+x_{237}+x_{456}+x_{18}x_{237}+x_{237}x_{456}+x_{456}x_{18}))\sum_{n\geq 0}x^{n\delta}
\end{align}

$\mathbb{C}^3/(\mathbb{Z}_4\times\mathbb{Z}_2)  (1, 0, 3)(0, 1, 1)$ (model 2 of \cite{HanSeo})
\begin{align}
     \Omega_\ast=&\sum_ix_i+((\mathbb{L}^{3/2}+6\mathbb{L}^{1/2}+\mathbb{L}^{-1/2})x^\delta+\mathbb{L}^{1/2}(x_{1458}+x_{2367}+x_{2864}+x_{1753}+x_{12}+x_{34}+x_{56}+x_{78}\nn\\&+x_{12}x_{34}+x_{34}x_{56}+x_{56}x_{78}+x_{78}x_{12}+x_{12}x_{34}x_{56}+x_{34}x_{56}x_{78}+x_{56}x_{78}x_{12}+x_{78}x_{12}x_{34}))\sum_{n\geq 0}x^{n\delta}
\end{align}

$\ast$ $\mathbb{C}^3/(\mathbb{Z}_3\times\mathbb{Z}_3)  (1, 0, 2)(0, 1, 2)$ (model 1 of \cite{HanSeo})
\begin{align}
     \Omega_\ast=&\sum_ix_i+((\mathbb{L}^{3/2}+7\mathbb{L}^{1/2}+\mathbb{L}^{-1/2})x^\delta+\mathbb{L}^{1/2}(x_{153}+x_{678}+x_{294}+x_{153}x_{678}+x_{678}x_{294}\nn\\&+x_{294}x_{153}+x_{189}+x_{237}+x_{456}+x_{189}x_{237}+x_{237}x_{456}+x_{456}x_{189}+x_{126}\nn\\&+x_{597}+x_{348}+x_{126}x_{597}+x_{597}x_{348}+x_{348}x_{126}))\sum_{n\geq 0}x^{n\delta}
\end{align}

\bibliography{ref}
\bibliographystyle{utphys}
\end{document}